\documentclass[12pt,amsfonts]{article}

\usepackage{bm}

\usepackage{latexsym}
\usepackage{amsfonts,amssymb}
\usepackage{mathrsfs} 
\usepackage{graphicx}
\usepackage{pdfpages}
\usepackage{multicol}

\usepackage[T1]{fontenc}
\usepackage{lmodern}
\usepackage{esint}

\usepackage{tikz}
\usetikzlibrary{math} 
\usetikzlibrary{calc} 

\newtheorem{theorem}{Theorem}[section]
\newtheorem{lemma}[theorem]{Lemma}%[section]
\newtheorem{remark} %[theorem] 
{Remark}%[section]
\newtheorem{corollary}[theorem]{Corollary}%[section]
\newtheorem{proposition}[theorem]{Proposition}%[section]

\newtheorem{definition}[theorem]{Definition}
\newcommand{\br}[1]{\begin{remark}\label{#1}\rm}
\newcommand{\er}{\end{remark}}
\newcommand{\bd}[1]{\begin{definition}\label{#1}\rm}
\newcommand{\ed}{\end{definition}}
\newcommand{\bt}[1]{\begin{theorem}\label{#1}}
\newcommand{\et}{\end{theorem}}
\newcommand{\bprop}[1]{\begin{proposition}\label{#1}}
\newcommand{\eprop}{\end{proposition}}
\newcommand{\bcor}[1]{\begin{corollary}\label{#1}}
\newcommand{\ecor}{\end{corollary}}

\newcommand{\lra}{\longrightarrow}

\newcommand{\stack}[2]{\raisebox{-2pt} 
{\renewcommand{\arraystretch}{.01} 
\begin{tabular}{c} 
$#2$\\$\scriptscriptstyle #1$ 
\end{tabular} 
}}

\renewcommand{\l}{{\it leb}}

\newcommand{\vp}{\varphi}
\newcommand{\ve}{\varepsilon}
\newcommand{\nid}{\noindent}
\newcommand{\qed}{\hfill$\Box$} 

\def\1{\, {\rm I}\mskip-10mu 1}

\renewcommand{\t}[1]{\tilde{#1}} 

\newcommand{\bmu}{\mbox{\boldmath${\mu}$}}
\newcommand{\bnu}{\mbox{\boldmath${\nu}$}} 
\newcommand{\sbmu}{\mbox{\tiny\boldmath${\mu}$}} 
\newcommand{\sbnu}{\mbox{\tiny\boldmath${\nu}$}}

\begin{document}
\title{Knudsen Type Group for Time in $\mathbb{R}$ and Related Boltzmann Type 
Equations} 
\par
\author{J\"org-Uwe L\"obus
%%\footnotemark[1]
\\ Matematiska institutionen \\ 
Link\"opings universitet \\ 
SE-581 83 Link\"oping \\ 
Sverige 
}
\date{}
\maketitle
{\footnotesize
\noindent
\begin{quote}
{\bf Abstract}
%\footnotetext[1]{\noindent } 
We consider certain Boltzmann type equations on a bounded physical and 
a bounded velocity space under the presence of both, reflective as well 
as diffusive boundary conditions. We introduce conditions on the shape of 
the physical space and on the relation between the reflective and the 
diffusive part in the boundary conditions such that the associated Knudsen 
type semigroup can be extended to time $t\in\mathbb{R}$. Furthermore, we 
provide conditions under which there exists a unique global solution to a 
Boltzmann type equation for time $t\ge 0$ or for time $t\in [\tau_0,\infty)$ 
for some $\tau_0<0$ which is independent of the initial value at time 0. 
Depending on the collision kernel, $\tau_0$ can be arbitrarily small. 
\noindent 

{\bf AMS subject classification (2020)} primary 35Qxx, secondary 76P05

\noindent
{\bf Keywords}  Boltzmann type equation, Knudsen type group, spectral 
analysis 

\end{quote}
}
%\bigskip
%\noindent (Running title: ) 

\section{Introduction}\label{sec:1}
\setcounter{equation}{0} 

Boltzmann type equations provide a mathematical description of rarefied 
gases in vessels. In the paper we consider such equations on a bounded 
physical space $\Omega$ and a bounded velocity space $V$. The boundary 
conditions include both, a specular reflective and a diffusive part. The 
goal is to formulate conditions on the parameters of the equation, the 
shape of the physical space, and the boundary conditions guaranteeing 
existence and uniqueness of the solution to a Boltzmann type equation 
for time $t\ge 0$ or time $t\in [\tau_0,\infty)$. Here $\tau_0<0$ is 
independent of the initial value at time 0. In particular we are 
concerned with the {\it integrated} or {\it mild} form 
\begin{eqnarray*}
\quad p(r,v,t)=S(t-\tau)\, p(r,v,\tau)+\lambda\int_\tau^t S(t-s)\, Q(p,p) 
\, (r,v,s)\, ds\, ,  
\end{eqnarray*} 
$(r,v,t)\in\Omega\times V\times [\tau,\infty)$,where $\tau_0\le\tau\le 0$, 
endowed with the boundary conditions. Here $Q$ denotes the collision 
operator and $S$ stands for the associated Knudsen type transport 
(semi)group. In addition, the initial value $p_0=p(\cdot,\cdot,0)$ is 
uniformly bounded from above a.e.~on $\Omega\times V$. Moreover, $p_0=p 
(\cdot, \cdot,0)$ is assumed to be uniformly positively bounded from below 
and above a.e.~on $\Omega\times V$. The parameter $\lambda>0$ is meaningful 
for the construction of preliminary local solutions and controls $\tau_0 
\equiv\tau_0(\lambda)$ by $\lim_{\lambda\to 0}\tau_0(\lambda)=-\infty$. 

More specifically, the collision operator involves spatial smearing. The 
collision kernel is bounded and in a certain sense decaying. The boundary 
conditions for the Knudsen type transport (semi)group and hence also for the 
Boltzmann type equation are explained in detail in the beginning of Subsection 
\ref{sec:1:2}. In particular, there is a certain balance between the specular 
reflective part and the diffusive part, see Proposition 1 and Theorem 2. The 
velocity redistribution density of the diffusive part is bounded from below as 
well as above and otherwise non-specified.
\medskip

The central part of the paper is the analysis of the Knudsen type transport 
(semi)group $S$. The proof of its existence for time $t\in\mathbb{R}$ is 
prepared in Section \ref{sec:3}. In Subsection \ref{sec:3:1}, a.e.~uniform 
positive boundedness from below and above of $S(t)\, p_0(r,v)$ on $(r,v,t)
\in\Omega\times V\times [0,\infty)$ is demonstrated. In Subsection 
\ref{sec:3:2}, certain explicit representations of $S(t)\, p_0$, $t\ge 0$, 
are provided and strong continuity in $L^1(\Omega\times V)$ is proved. The 
actual proof of the existence of $S(t)$ for $t\in\mathbb{R}$, is carried out 
in Section \ref{sec:4} using spectral theoretical methods and facts from the 
theory of mathematical billiards. Important for this proof is the particular 
shape of the physical boundary $\partial\Omega$ and the form of the boundary 
conditions. 

In Section \ref{sec:5} we are then concerned with existence and uniqueness 
of solutions to Boltzmann type equations, first on some small time interval 
$[0,T]$, see Subsection \ref{sec:5:1}. Using this and a certain representation 
of the Boltzmann type equation along the paths of transport, existence and 
uniqueness is shown in Subsection \ref{sec:5:2} for time $t\ge 0$. As a 
consequence of the results of Sections \ref{sec:3} and \ref{sec:4} we derive 
then existence and uniqueness of solutions to Boltzmann type equations for 
time in $[\tau_0,\infty)$. Denoting by $s\le 0$ a time for which $p(\cdot, 
\cdot ,s)$ is uniformly positively bounded from below a.e.~on $\Omega\times 
V$, these solutions are uniformly positively bounded from below and above 
a.e.~on $\Omega\times V\times [s,T]$ for any $T>0$. 

\subsection{Relation to existing literature}\label{sec:1:1} 

Significant progress in the analysis of {\it collisionless gases} in bounded 
domains has been achieved over the past decade. For bounded $C^1$-domains 
$\Omega$ endowed with certain composite boundary conditions, containing a 
specular reflective as well as a diffusive part, in \cite{LMR20} spectral 
methods are used in order to examine the associate transport semigroup. 
As a result, conditions can be formulated such that a unique invariant 
density exists. It is also demonstrated that if no invariant density exists 
then the total mass of the transport semigroup may concentrate near certain 
sets of zero measure as $t\to\infty$. 

On the other hand, the optimal rate of convergence $1/(1+t)^{d-}$ to the 
stationary density has been established in e.g. \cite{AG11}, \cite{KLT13}, 
and \cite{Be20}. In the latest paper \cite{Be20} this rate is even shown 
to hold for bounded $C^2$-domains $\Omega$ under the presence of certain 
composite boundary conditions, consisting of a specular reflective as well 
as a diffusive part. 

We mention that the setup of the present paper is not compatible with 
these recent developments since the global condition (vii) below fails for 
everywhere differentiable boundaries $\partial\Omega$, cf. the instructive 
figures 6 and 7 below. 

For the Knudsen type transport (semi)group $S(t)$ considered in the present 
paper, the convergence to a unique stationary density $\bar{g}$ as $t\to 
\infty$ is demonstrated. Moreover $c\le S(t)p_0\le C$ a.e.~for all $t\in 
[t_0,\infty)$ and some $0<c\le C<\infty$ is shown, whenever 
$S(t_0)p_0\in L^\infty(\Omega\times V)$, $S(t_0)p_0\ge 0$ a.e. and $\|1/ 
S(t_0)p_0\|_{L^\infty(\Omega\times V)}<\infty$. In addition it is proved 
that the stationary density $\bar{g}$ is a.e.~uniformly positively bounded 
from below and above. 

\medskip

In \cite{CPW98} and \cite{Lo18} the {\it particle redistribution} after collision 
with a physical boundary is diffusive with a {\it non-specified redistribution} 
density. This comes close to the situation of the present paper. 
\medskip 

The physical relevance of a {\it boundary condition containing a specular 
reflective fraction and a diffusive fraction} has already been pointed out 
by C. Maxwell, see the Appendix to \cite{Ma1879}. Since then this type of 
boundary conditions, also called {\it Maxwell boundary conditions}, has 
drawn the attention of many researchers. The state of the art is documented 
in \cite{S-R09}, Theorem 2.3.4 and its proof and in \cite{Mi10}. Even 
extensions of Maxwell boundary conditions have been properly examined, most 
recently in \cite{Ch20}. In the just cited references the redistribution 
density in the diffusive part of the boundary conditions is a certain wall 
Maxwellian while in the present paper this density remains non-specified. 
We emphasize that in the present paper the particular boundary conditions 
together with the shape of the boundary imply unique extensions of the 
Knudsen type transport semigroup to time $t\in\mathbb{R}$ as well as the 
solution to the Boltzmann type equation to time $t\in [\tau_0,\infty)$. 
\medskip

A recurring problem is the existence of {\it lower bounds on solutions} 
to different types of Boltzmann equations. There are several results for 
the spatially homogeneous case. For example, in \cite{PW97} a certain 
time independent Maxwellian lower bound has been derived and in 
\cite{Fo01} Malliavin calculus arguments have been applied to establish 
a strictly positive function bound in a certain weak sense. The papers 
\cite{Mo05}, \cite{Br15-1}, \cite{Br15-2}, and \cite{IMS20} are 
concerned with several types of function lower bounds on solutions to 
non-homogeneous Boltzmann equations. In particular, they refer to 
different types of boundary conditions, such as periodic, specular 
reflective, and Maxwellian diffusive ones. For a result concerning a 
stationary solution to a Boltzmann equation on a bounded physical space 
with general diffusive boundary conditions, see \cite{Lo18}. 
\medskip

{\it Boltzmann's $H$-theorem} is a crucial tool in order to prove 
irreversibility of solutions to the Boltzmann equation. As for example in 
\cite{BGSS18}, the proof of a corresponding $H$-theorem is possible in a 
number of situations, in particular when the velocity space $V$ is $\{v 
\in\mathbb{R}^d:|v|>0\}$ and the physical space $\Omega$ is $\mathbb{R}^d$ 
or the $d$-dimensional torus. On the other hand, for bounded $\Omega$ or $V$, 
certain types of boundary conditions may not be compatible to known proofs 
of the $H$-theorem. For example, the proof in \cite{RV08}, Theorem 1 of 
Section 1.1.2 cannot be applied or modified in order to cope with diffusive 
boundary conditions as used in the present paper since boundary conditions 
containing a diffusive component do not conserve energy, in general. For the 
characterization of the entropy under the presence of diffusive boundary 
conditions, there is a modified $H$-functional, $H=\iint p(\cdot,\cdot, t) 
\log p(\cdot,\cdot, t)\, dv\, dr+\beta_w\iint |v|^2 p(\cdot ,\cdot, t)\, dv 
\, dr$. Here $\beta_w$ is the inverse temperature. Under certain isothermal 
diffusive boundary conditions this $H$ decreases in time, see (4.1), (4.2) 
of Chapter 9 in \cite{CIP94}. For further investigations in terms of relative 
entropy, we refer to \cite{S-R09}. 
\medskip

The papers \cite{BC02} an \cite{WZ06} provide {\it eternal solutions}, i.e.  
non-stationary non-negative solutions to different versions of the Boltzmann 
equation for $t\in\mathbb{R}$. The infinite energy solution constructed 
in \cite{BC02}, has been numerically examined in \cite{RW05}, Sections 4.4 
and B.7. Besides these known eternal solutions, 
existence for time $t\in\mathbb{R}$ in form of {\it stationarity} is  
well-established for the Boltzmann equation. In particular, we would like to 
draw attention to the paper \cite{EM20} since it reviews several approaches 
to stationary solutions of the Boltzmann equation which are different from 
the Maxwell-Boltzmann equilibrium. The construction of a unique stationary 
solution is also an objective in \cite{CPW98} where a certain Boltzmann 
equation satisfying diffusive boundary conditions with a non-specified 
redistribution density is analyzed. It has been shown in this paper that 
for bounded $\Omega$, $V=\{v\in\mathbb{R}^d:|v|>v_0>0\}$, and appropriate 
initial configurations, the distribution of a suitable $N$-particle system 
converges in a certain sense to the unique stationary solution of the 
Boltzmann equation considered.

\subsection{Main results}\label{sec:1:2}

For a certain class of physical spaces which includes certain convex polygons 
in dimension $d=2$ or a certain class of polyhedrons in $d=3$ let us consider 
the Knudsen type semigroup $S(t)$, $t\ge 0$. It is introduced for 
time $t\ge 0$ by the solution to the initial boundary value problem  
\begin{eqnarray*}
\left(\frac{d}{dt}+v\circ\nabla_r\right)(S(t)p_0)(r,v)=0\quad \mbox{\rm on} 
\quad (r,v,t)\in \Omega\times V\times [0,\infty)\, ,  
\end{eqnarray*} 
$S(0)p_0=p_0$, and the boundary conditions 
\begin{eqnarray*} 
(S(t)p_0)(r,v)=\omega\, (S(t)p_0)(r,R_r(v))+(1-\omega)J(r,t)(S(\cdot)p_0)M 
(r,v)\, ,\quad t>0,  
\end{eqnarray*} 
for all $(r,v)\in\partial^{(1)}\Omega\times V$ with $v\circ n(r)\le 0$. 
Here $n(r)$ is the outer normal at the point $r$ belonging to the boundary 
part $\partial^{(1)}\Omega$ obtained from $\partial\Omega$ by removing all 
vertices and edges if $d=3$. The symbol ``$\circ$" denotes the inner product 
in $\mathbb{R}^d$. Furthermore, $R_r(v):=v-2v\circ n(r)\cdot n(r)$ for $(r,v) 
\in\partial^{(1)}\Omega\times V$ indicates reflection of the velocity $v$ at  
$r\in\partial^{(1)}\Omega$ and  
\begin{eqnarray*} 
J(r,t)(S(\cdot)p_0):=\int_{v\circ n(r)\ge 0}v\circ n(r)\, S(t)p_0(r,v)\, dv 
\, .  
\end{eqnarray*} 
The function $M$ on $\{(r,v):r\in\partial^{(1)}\Omega,\ v\in V,\ v\circ n(r) 
\le 0\}$ quantifies the diffusive part of the boundary conditions; it is 
positive, continuous, and uniformly bounded from below and above. The constant 
$\omega\in (0,1)$ controls the relation between reflective and diffusive 
boundary conditions. 
\bigskip\bigskip\bigskip\medskip

\fbox{ 
 
\begin{tikzpicture}[scale=3.2]

\useasboundingbox (0.06,0.023) rectangle (1.94,0.975);  
        
        \draw[color=green,dashed,-<] (1.2,1) -- ($(1.2,1) !0.5cm! (0,0.92)$);
        \draw[color=green,dashed,-<] (1.2,1) -- ($(1.2,1) !0.5cm! (0,0.75)$);
        \draw[color=green,dashed,-<] (1.2,1) -- ($(1.2,1) !0.5cm! (0,0.6)$);
        \draw[color=green,dashed,-<] (1.2,1) -- ($(1.2,1) !0.5cm! (0,0.4)$);
        \draw[color=green,dashed,-<] (1.2,1) -- ($(1.2,1) !0.5cm! (0,0.2)$);
        \draw[color=green,dashed,-<] (1.2,1) -- ($(1.2,1) !0.5cm! (0.1,0)$);
        \draw[color=green,dashed,-<] (1.2,1) -- ($(1.2,1) !0.5cm! (0.4,0)$);
        \draw[color=green,dashed,-<] (1.2,1) -- ($(1.2,1) !0.5cm! (0.7,0)$);
        \draw[color=green,dashed,-<] (1.2,1) -- ($(1.2,1) !0.5cm! (1,0)$);
        \draw[color=green,dashed,-<] (1.2,1) -- ($(1.2,1) !0.5cm! (1.3,0)$);
        \draw[color=green,dashed,-<] (1.2,1) -- ($(1.2,1) !0.5cm! (1.6,0)$);
        \draw[color=green,dashed,-<] (1.2,1) -- ($(1.2,1) !0.5cm! (1.9,0)$);
        \draw[color=green,dashed,-<] (1.2,1) -- ($(1.2,1) !0.5cm! (2,0.2)$);
        \draw[color=green,dashed,-<] (1.2,1) -- ($(1.2,1) !0.5cm! (2,0.4)$); 
        \draw[color=green,dashed,-<] (1.2,1) -- ($(1.2,1) !0.5cm! (2,0.6)$);
        \draw[color=green,dashed,-<] (1.2,1) -- ($(1.2,1) !0.5cm! (2,0.8)$); 
        \draw[color=green,dashed,-<] (1.2,1) -- ($(1.2,1) !0.5cm! (2,0.95)$);
        
        \node[color=black] at (0.85,1.15) {\footnotesize $S(t)p_0$ at $\quad\quad(r,v)$}; 
        
        \node[color=black] at (1.2,0.35) {\footnotesize $\int v\circ n(r)\, S(t)p_0(r,v)\, dv$}; 
      
        \filldraw[white] (1.19,0.93) circle (1.8pt); 
        
        \node[color=black] at (1.19,0.93) {\footnotesize $r$}; 

        \draw[color=black,thick,->] ($(-0.01,0.8) !0.6cm! (1.2,1)$) -- (1.2,1); 

\end{tikzpicture} 

} 
\hspace{1.5cm} 
\fbox{

\begin{tikzpicture}[scale=3.2]

\useasboundingbox (0.06,0.023) rectangle (1.94,0.975);
        
        \draw[color=red,dashed,->] (1.2,1) -- ($(1.2,1) !0.5cm! (0,0.92)$);
        \draw[color=red,dashed,->] (1.2,1) -- ($(1.2,1) !0.5cm! (0,0.75)$);
        \draw[color=red,dashed,->] (1.2,1) -- ($(1.2,1) !0.5cm! (0,0.6)$);
        \draw[color=red,dashed,->] (1.2,1) -- ($(1.2,1) !0.5cm! (0,0.4)$);
        \draw[color=red,dashed,->] (1.2,1) -- ($(1.2,1) !0.5cm! (0,0.2)$);
        \draw[color=red,dashed,->] (1.2,1) -- ($(1.2,1) !0.5cm! (0.1,0)$);
        \draw[color=red,dashed,->] (1.2,1) -- ($(1.2,1) !0.5cm! (0.4,0)$);
        \draw[color=red,dashed,->] (1.2,1) -- ($(1.2,1) !0.5cm! (0.7,0)$);
        \draw[color=red,dashed,->] (1.2,1) -- ($(1.2,1) !0.5cm! (1,0)$);
        \draw[color=red,dashed,->] (1.2,1) -- ($(1.2,1) !0.5cm! (1.3,0)$);
        \draw[color=red,dashed,->] (1.2,1) -- ($(1.2,1) !0.5cm! (1.6,0)$);
        \draw[color=red,dashed,->] (1.2,1) -- ($(1.2,1) !0.5cm! (1.9,0)$);
        \draw[color=red,dashed,->] (1.2,1) -- ($(1.2,1) !0.5cm! (2,0.2)$);
        \draw[color=red,dashed,->] (1.2,1) -- ($(1.2,1) !0.5cm! (2,0.4)$); 
        \draw[color=red,dashed,->] (1.2,1) -- ($(1.2,1) !0.5cm! (2,0.6)$);
        \draw[color=red,dashed,->] (1.2,1) -- ($(1.2,1) !0.5cm! (2,0.8)$); 
        \draw[color=red,dashed,->] (1.2,1) -- ($(1.2,1) !0.5cm! (2,0.95)$);
        
        \node[color=blue] at (1.4,1.15) {\footnotesize$(r,v)$};
        \node[color=black] at (0.5,1.15) {\footnotesize $S(t)p_0$ at $\quad\quad(r,R_v(v))$}; 
       
        \node[color=black] at (1.2,0.35) {\footnotesize $J(r,t)(S(\cdot)p_0)\cdot M(r,v)$};

        \filldraw[white] (1.19,0.93) circle (1.8pt);
        \node[color=black] at (1.19,0.93) {\footnotesize $r$}; 

        \draw[color=blue,thick,->] (1.2,1) -- ($(1.2,1) !0.6cm! (2.01,0.867)$); 
        \draw[color=black,thick,->] ($(-0.01,0.8) !0.6cm! (1.2,1)$) -- (1.2,1); 

\end{tikzpicture} 

 } 

\bigskip

\noindent
{\small{{\bf Fig.~1} Case ${v\circ n(r)\ge 0}$, $\quad r\in\partial\Omega$ 
\hspace{2.1cm} {\bf Fig.~2} Case ${v\circ n(r)\le 0}$, $\quad r\in\partial 
\Omega$ \\ \\ 
The directions of the arrows represent the directions of the velocities $v$. 
The lengths of the arrows do not represent $|v|$. Fig.~1: $\bm{v\circ  
n(r)\ge 0}:$ The particle densities directed from $r\in\partial\Omega$ to the 
outside of $\Omega$ at time $t$ (green) are integrated to $J(r,t)(S(\cdot)p_0) 
=\int v\circ n(r)\, S(t)p_0(r,v)\, dv$ where the integral ranges over all $v$ 
with $v\circ n(r)\ge 0$. The integral $J$ controls to the diffusive part in the 
boundary conditions. In addition, $S(t)p_0(r,v)$ (black) contributes to the 
reflective part in the boundary conditions. Fig.~2: $\bm 
{v\circ n(r)\le 0}:$ There are two sources of the particle densities directed 
from $r\in\partial\Omega$ to the inside of $\Omega$ at time $t$. With probability 
$\omega$ these particle densities come from specular reflection (blue). Moreover, 
with probability $1-\omega$ they emerge from $J(r,t)(S(\cdot)p_0)$ (red) and are 
distributed over $\{v:v\circ n(r)\le 0\}$ according to the density $M(r,v)$ with 
respect to $|v\circ n(r)|\, dv$.}} 

\bigskip\medskip

The crucial technical part of the paper is the following spectral property. 
\medskip 

\noindent
{\bf Proposition 1} {\it (Corollary \ref{Corollary4.4} (c) below) 
Under certain conditions on the shape of $\partial\Omega$ there exist 
$k_0\in\mathbb{N}$ and $m_1>-\infty$ such that for $2^{-\frac{1}{k_0}} 
<\omega<1$ and $t>0$, the resolvent set of $S(t)$ contains the set 
$\left\{\lambda=e^{t\mu}:\mathfrak{Re}\mu<m_1\right\}\cup\left\{\lambda 
=e^{t\mu}:\mathfrak{Re}\mu>0\right\}$. }
\medskip 

\noindent 
The proof of Proposition 1 includes explicit non-standard solutions to 
certain Banach space valued differential equations. By means of 
Proposition 1 we have proved the subsequent theorem. Here the condition 
$2^{-\frac{1}{k_0}}<\omega<1$ is sufficient for the existence of an 
extension of the Knudsen type transport semigroup to time $t\in\mathbb 
{R}$. It remains an open question whether this condition, or some 
modification, is necessary and sufficient for the Knudsen type transport 
semigroup to admit an extension to $t\in\mathbb{R}$.
\medskip 

\noindent
{\bf Theorem 2} {\it (Theorem \ref{Theorem4.5} below) Suppose the 
conditions of Proposition 1 are satisfied. There exists $k_0\in 
\mathbb{N}$ such that for $2^{-\frac{1}{k_0}}<\omega<1$, the semigroup 
$S(t)$, $t\ge 0$, extends to a strongly continuous group in $L^1(\Omega 
\times V)$ which we will denote by $S(t)$, $t\in\mathbb{R}$. } 
\medskip 

\noindent
The following boundedness property of the Knudsen-type group is the next 
major result. It is an important technical tool throughout the whole 
paper and of independent interest. 
\medskip 

\noindent
{\bf Proposition 3} {\it (Lemma \ref{Lemma3.1} and Corollary 
\ref{Corollary4.6} (a) below) Let $p_0\in L^1(\Omega\times V)$ and $t_0 
\in\mathbb{R}$ such that $S(t_0)p_0\in L^\infty(\Omega\times V)$. Then 
there are finite real numbers $p_{0,{\rm min}}$ and $p_{0,{\rm max}}$ 
such that 
\begin{eqnarray*}
p_{0,{\rm min}}\le S(t)\, p_0\le p_{0,{\rm max}}\quad\mbox{\rm a.e. on } 
\Omega\times V
\end{eqnarray*} 
for all $t\ge t_0$. In particular, if $S(t_0)p_0\ge 0$ a.e. and $\|1/ 
S(t_0)p_0\|_{L^\infty(\Omega\times V)}<\infty$ then there exists $p_{ 
0,{\rm min}}>0$.  A corresponding statement holds for time $t\ge 0$.}
\medskip 

\noindent
Let us  turn to Boltzmann type equations for time $t\ge 0$ or time $t\in 
[\tau_0,\infty)$. The following theorem presents the result on the unique 
existence of global solutions to their integrated (mild) as well as their 
{\it differentiated (classical)} versions, where the terms {\it mild} and 
{\it classical} are used as in the standard reference \cite{Pa83}. It has 
been proved under reasonable assumptions on the collision operator $Q$. 
Theorem 2 and Proposition 3 are essential for the proof of the subsequent
Theorem 4. 
\medskip

\noindent 
{\bf Theorem 4} {\it (Theorem \ref{Theorem5.6}, Corollary 
\ref{Corollary5.7}, and Theorem \ref{Theorem5.8} below) Suppose that the 
conditions of Proposition 1 are satisfied. There exists $\tau_0\equiv\tau_0 
(\lambda)<0$ with $\lim_{\lambda\to 0}\tau_0(\lambda)=-\infty$ such that 
the following holds. \\ 
(a) Let $\|p_0\|_{L^1(\Omega\times V)}=1$ and suppose that there are 
constants $0<c\le C<\infty$ with $c\le p_0\le C$ a.e. on $\Omega\times V$. 
Then there is a unique map $[\tau_0,\infty)\ni t\mapsto p_t(p_0)\equiv p 
(\cdot ,\cdot ,t)\in L^1(\Omega\times V)$ such that for every $\tau_0\le\tau 
\le 0$ 
\begin{eqnarray}\label{1.1} 
p(r,v,t)=S(t-\tau)\, p(r,v,\tau)+\lambda\int_\tau^t S(t-s)\, Q(p,p)\, (r,v,s) 
\, ds 
\end{eqnarray} 
a.e.~on $\Omega\times V\times [\tau,\infty)$ and $p(\cdot ,\cdot ,0)=p_0$. 
The solution $p\equiv p(p_0)$ to (\ref{1.1}) has the following properties. 
\begin{itemize} 
\item[(1)] The map $[\tau_0,\infty)\ni t\mapsto p(\cdot ,\cdot ,t)\in L^1 
(\Omega\times V)$ is continuous with respect to the topology in $L^1(\Omega 
\times V)$. 
\item[(2')] We have $\|p(\cdot ,\cdot ,t)\|_{L^1(\Omega\times V)}\le K< 
\infty$, $t\ge\tau_0$ and $\|p(\cdot ,\cdot ,t)\|_{L^1(\Omega\times V)}=1$, 
$t\ge 0$. 
\item[(3')] If $p(\cdot,\cdot ,s)\ge c'$ a.e. on $\Omega\times V$ for some 
$c'>0$ and some $s\in [\tau_0,0]$ then there exists a strictly decreasing 
positive function $[s,\infty)\ni t\mapsto c_t\equiv c_t(c',s)$ such that 
$p(\cdot ,\cdot ,t)\ge c_t$ a.e. on $\Omega\times V$. 
\item[(4')] For $s$ introduced in (3') and $t\ge s$ we have 
\begin{eqnarray*}
\|p(\cdot ,\cdot ,t)\|_{L^\infty(\Omega\times V)}\le a\cdot\exp\left\{\lambda 
b\cdot |t|\vphantom{l^1}\right\}\, ,    
\end{eqnarray*}
for some constants $a,b>0$. 
\end{itemize}
If (\ref{1.1}) is modified by replacing $p$ with $p/(1\vee\|p/K\|_{L^1(\Omega 
\times V)})$, then the modified equation has a unique continuous solution 
$\mathbb{R}\ni t\mapsto \t p_t(p_0)\equiv\t p(\cdot,\cdot ,t)\in L^1(\Omega 
\times V)$ with $\t p_t(p_0)=p_t(p_0)$ for all $t\ge\tau_0$. The following 
holds.  
\begin{itemize} 
\item[(2'')] For all $t\in\mathbb{R}$ we have $\int_\Omega\int_V \t 
p(\cdot,\cdot ,t)\, dv\, dr =1$ and for all $t\le 0$ it holds that 
\begin{eqnarray*}
\left\|\t p_t(p_0)\right\|_{L^1(\Omega\times V)}\le\exp\left\{-t\cdot 
\left(-m_1+a_1\lambda K\vphantom{l^1}\right)\right\}
\end{eqnarray*}
where $a_1$ is some constant depending on the collision kernel and $m_1 
<0$ is the constant introduced in Proposition 1. The constant $K$ is the 
one of (2'). It is related to $a_1,\tau_0,m_1$ by $-\tau_0=\log(K)/(-m_1 
+a_1\lambda K)$. 
\end{itemize} 
(b) Let $p_0\in D(A)$ with $\|p_0\|_{L^1(\Omega\times V)}=1$ and suppose 
that there are constants $0<c\le C<\infty$ with $c\le p_0\le C$ a.e. on 
$\Omega\times V$. Then there is a unique map $[\tau_0,\infty)\ni t\mapsto 
p_t(p_0)\equiv p(\cdot,\cdot ,t)\in D(A)$ which is continuous on $[\tau_0, 
\infty)$ and continuously differentiable on $(\tau_0,\infty)$ with respect 
to the topology of $L^1(\Omega\times V)$ such that 
\begin{eqnarray*}
\frac{d}{dt}\, p(\cdot ,\cdot ,t)=Ap(\cdot ,\cdot ,t)+\lambda\, Q(p,p)\, 
(\cdot ,\cdot ,t)\, ,\quad t\ge\tau_0, 
\end{eqnarray*} 
and $p(\cdot ,\cdot, 0)=p_0$. Here, $d/dt$ denotes the derivative in $L^1 
(\Omega\times V)$. At $t=\tau_0$ it is the right derivative. This map 
coincides with the solution to the equation (\ref{1.1}) if there $p(\cdot 
,\cdot, 0)=p_0\in D(A)$. We have properties (2')-(4') and (2'') of part (a).
}
\medskip

\noindent
The time $\tau_0\equiv\tau_0(\lambda)$ and the constants $a$ and $b$ in 
(4') are explicitly determined.

\section{Preliminaries}\label{sec:2}
\setcounter{equation}{0} 

Let $d\in\{2,3\}$. For a bounded set $\Gamma\subseteq\mathbb{R}^d$ let 
$\overline{\Gamma}$ denote its closure in $\mathbb{R}^d$. Let $\Omega 
\subset\mathbb{R}^d$ be a bounded domain, called the {\it physical 
space}. Let us suppose that $\partial\Omega=\bigcup_{i=1}^{n_\partial} 
\overline{\Gamma_i}$ where the $\Gamma_i$ are $(d-1)$-dimensional 
manifolds which are smooth up to the boundary and, at the same time, 
open sets in the topology of $\partial\Omega$, $i\in \{1,\ldots,n_ 
\partial\}$. Let us furthermore assume that for $i\neq j$ the intersection 
$\overline {\Gamma_i}\cap\overline{\Gamma_j}$ either is the empty set or 
a smooth closed $(d-2)$-dimensional manifold which, by definition, for 
$d=2$ is just a single point in $\mathbb{R}^2$. 

In case of $d=3$, for any two neighboring $\Gamma_i$ and $\Gamma_j$ 
let us call $x\in\overline{\Gamma_i}\cap\overline{\Gamma_j}$ an {\it 
end point} if for any sufficiently small open ball $B_x$ in $\mathbb 
{R}^d$ with center $x$ the set $\partial B_x\cap (\overline{\Gamma_i} 
\cap\overline{\Gamma_j})$ consists of exactly one point. For the boundary 
$\partial\Omega$ let us also assume that there is $\xi\in (0,\pi)$ such 
that for any two neighboring $\Gamma_i$ and $\Gamma_j$ and any $x\in 
\overline{\Gamma_i}\cap\overline{\Gamma_j}$, not being an end point if 
$d=3$, we have the following. For the angle $\xi(i,j;x)$ between the 
rays $R_i$ and $R_j$ from $x$ tangential to $\Gamma_i$ and $\Gamma_j$ 
respectively, and orthogonal to $\overline{\Gamma_i}\cap\overline 
{\Gamma_j}$ if $d=3$, it holds that 
\begin{eqnarray}\label{2.1}
\xi<\xi(i,j;x)<2\pi-\xi\, . 
\end{eqnarray}
The latter guarantees compatibility with condition (vii) below and Assumption 
A in the proof of Theorem 2.1 in \cite{CPW98}, as we will discuss in Remark 
\ref{3} below. 
\bigskip\bigskip

\begin{tikzpicture}[scale=0.7] 
\draw[thick,rounded corners=50pt] (4,0.5)--(7,0.5)--(10,0);
\draw[thick,rounded corners=20pt] (4,0.5)--(3.7,1.5)--(3,2.5);
\draw[thick,rounded corners=25pt] (10,0)--(10.8,1)--(12,2); 
\draw[thick,rounded corners=25pt] (3,2.5)--(5,2.9)--(7.2,3.5);
\draw[thick,rounded corners=25pt] (7.2,3.5)--(9.5,2.6)--(12,2);
\draw[thick,rounded corners=25pt] (7.2,3.5)--(7.2,4.3);
\draw[thick,rounded corners=25pt] (3,2.5)--(3.1,2.9);
\draw[thick,rounded corners=25pt] (12,2)--(12.3,2.6);
\draw[thick,dashed,rounded corners=25pt,color=red] (10,0)--(10.3,1.5)--(11,3); 
\draw[thick,dashed,rounded corners=25pt,color=red] (4,0.5)--(4.5,2)--(4,3); 
\draw[thick,dashed,rounded corners=25pt,color=red] (4,3)--(5.5,3.3)--(7,4);
\draw[thick,dashed,rounded corners=25pt,color=red] (7,4)--(9,3.3)--(11,3);
\draw[thick,dashed,rounded corners=25pt,color=red] (7,4)--(7,5);
\draw[thick,dashed,rounded corners=15pt,color=red] (4,3)--(3.85,3.5)--(3.4,4);
\draw[thick,dashed,rounded corners=15pt,color=red] (11,3)--(11.4,3.5)--(12,4);
\node at (7.5,1.5) {$\Gamma_1$};
\node[color=red] at (8,2.5) {\footnotesize$\Gamma_2$};
\node[color=red] at (11.5,2.5) {\footnotesize$\Gamma_3$};
\draw[color=white] (7,0) arc (90:70:3.5cm and 3.5cm);
\end{tikzpicture} 
\hspace{1.5cm} 
\begin{tikzpicture}[scale=0.7]
\draw[thick] (-0.04,0) arc (90:70:18cm and 18cm);
\draw[thick] (0.5,-2) arc (90:70:18cm and 18cm);
\draw[thick,color=red] (1.5,1) arc (90:70:18cm and 18cm);
\draw[thick] (6.12,-1.1) arc (30:0:4cm and 4cm);
\draw[thick] (0.0,0) arc (30:0:4cm and 4cm);
\draw[thick,color=red] (7.65,-0.06) arc (144.2:180:5.2cm and 5.2cm);
\draw[thick,color=red] (1.52,1) arc (143.8:156.5:5.2cm and 5.2cm);
\draw[thick,color=red,dashed] (1.52,1) arc (144.2:180:5.2cm and 5.2cm);
\draw[thick,color=blue,->] (4.57,1.2) arc (115:134:6cm and 6cm);
\draw[thick,color=red,dashed] (3.5,-2.26)--(5,2.5);
\draw[thick,color=red,->] (4.08,-0.42)--(5,2.5); 
\draw[thick,color=black,dashed] (3.5,-2.26)--(2.5,2);
\draw[ultra thick,color=black,->] (3.04,-0.27)--(2.5,2);
\node at (6,-2) {$\Gamma_i$};
\node[color=red] at (6,0) {\footnotesize$\Gamma_j$};
\node at (3.5,-2.7) {\footnotesize$x$};
\node[color=black] at (2,1.7) {$R_i$};
\node[color=red] at (5.5,2.1) {\footnotesize$R_j$};
\node[color=blue] at (3.7,1.4) {\footnotesize$\xi(i,j;x)$};
\end{tikzpicture} 
\bigskip 

\noindent
{\small{{\bf Fig.~3}  \hspace{7.5cm}{\bf Fig.~4} 
\\ \\ 
Black indicates the visible faces and edges, red the non-visible ones. Fig.~3: 
The principal shape of $\partial\Omega$ is displayed in this figure. Fig.~4: 
(\ref{2.1}) is illustrated. }} 
\bigskip

Let $V:=\{v\in \mathbb{R}^d:0<v_{min}<|v|<v_{max}<\infty\}$ be the 
{\it velocity space} and let $\lambda>0$. Denote by $n(r)$ the outer 
normal at 
\begin{eqnarray*} 
r\in\partial^{(1)}\Omega:=\bigcup_{i=1}^{n_{\partial}}\Gamma_i\, ,  
\end{eqnarray*} 
indicate the inner product in $\mathbb{R}^d$ by ``$\circ$", and let 
$R_r(v):=v-2v\circ n(r)\cdot n(r)$ for $(r,v)\in\partial^{(1)}\Omega 
\times V$. 
\medskip 

For $(r,v,t)\in\Omega\times V\times [0,\infty)$, consider the {\it 
Boltzmann type equation} 
\begin{eqnarray}\label{2.2} 
\frac{d}{dt}\, p(r,v,t)=-v\circ\nabla_rp(r,v,t)+\lambda Q(p,p)\, (r,v,t)
\end{eqnarray} 
with boundary conditions 
\begin{eqnarray}\label{2.3} 
p(r,v,t)=\omega\, p(r,R_r(v),t) +(1-\omega)J(r,t)(p)M(r,v)\, ,\quad r 
\in\partial^{(1)}\Omega,\ v\circ n(r)\le 0, 
\end{eqnarray} 
for some $\omega\in [0,1)$ and initial probability density $p(0,\cdot,\cdot) 
:=p_0$ on $\Omega\times V$. Consider also its {\it integrated (mild) version} 
\begin{eqnarray}\label{2.4} 
p(r,v,t)=S(t)\, p_0(r,v)+\lambda\int_0^t S(t-s)\, Q(p,p)\, (r,v,s) 
\, ds\, . 
\end{eqnarray} 

The following global conditions (i)-(vii) on the terms in 
(\ref{2.2})-(\ref{2.4}) and the shape of $\Omega$ will be used in the 
paper. In the preliminary paragraph of each section we list those global 
conditions which we suppose in the section.
\begin{itemize} 
\item[(i)] For all $t\ge 0$ and all $r\in\partial^{(1)}\Omega$, the 
function $J$ is given by 
\begin{eqnarray*} 
J(r,t)(p)=\int_{v\circ n(r)\ge 0}v\circ n(r)\, p(r,v,t)\, dv\, .  
\end{eqnarray*} 
\item[(ii)] The real function $M$ on $\{(r,v):r\in\partial^{(1)}\Omega,
\ v\in V,\ v\circ n(r)\le 0\}$ has positive lower and upper bounds 
$M_{\rm min}$ and $M_{\rm max}$, is continuous on every $\Gamma_i$, $i\in 
\{1,\ldots ,n_\partial\}$, and satisfies
\begin{eqnarray*} 
\int_{v\circ n(r)\le 0}|v\circ n(r)|\, M(r,v)\, dv=1\, ,\quad r\in 
\partial^{(1)}\Omega. 
\end{eqnarray*} 
\item[(iii)] $S(t)$, $t\ge 0$, is given by the solution to the initial 
boundary value problem 
\begin{eqnarray}\label{2.5}
\left(\frac{d}{dt}+v\circ\nabla_r\right)(S(t)p_0)(r,v)=0\quad \mbox 
{\rm on}\quad (r,v,t)\in \Omega\times V\times [0,\infty)\, ,  
\end{eqnarray} 
$S(0)p_0=p_0$, and 
\begin{eqnarray*} 
(S(t)p_0)(r,v)=\omega\, (S(t)p_0)(r,R_r(v))+(1-\omega)J(r,t)(S(\cdot)p_0) 
M(r,v)\, ,\quad t>0,  
\end{eqnarray*} 
for all $(r,v)\in\partial^{(1)}\Omega\times V$ with $v\circ n(r)\le 0$. 
We call $S(t)$, $t\ge 0$, {\it Knudsen type semigroup}. If $\omega=0$ we 
indicate this in the notation by $S_0(t)$, $t\ge 0$. We  follow 
\cite{CPW98} and call $S_0(t)$, $t\ge 0$, just {\it Knudsen semigroup}. 
\item[(iv)] Denoting by $\chi$ the indicator function and setting $p:=0$ 
as well as $q:=0 $ on $\Omega\times (\mathbb{R}^d\setminus V)\times 
[0,\infty)$, the {\it collision operator} $Q$ is given by 
\begin{eqnarray*}
Q(p,q)(r,v,t)&&\hspace{-.5cm}=\frac12\int_{\Omega}\int_V\int_{S_+^{d-1}} 
B(v,v_1,e)h_\gamma(r,y)\chi_{\{(v^\ast ,v_1^\ast)\in V\times V\}}\times \\ 
&&\hspace{-1.0cm}\times\left(\left(p(r,v^\ast ,t)q(y,v_1^\ast ,t)-p(r,v,t) 
q(y,v_1,t)\vphantom{l^1}\right)\vphantom{\dot{f}}\right.\vphantom{\int} \\ 
&&\hspace{-.7cm}\left.+\left(q(r,v^\ast ,t)p(y,v_1^\ast ,t)-q(r,v,t)p(y, 
v_1,t)\vphantom{l^1}\right)\vphantom{\dot{f}}\right)\, de\, dv_1\, dy\, .  
\vphantom{\int} 
\end{eqnarray*} 
Here $S^{d-1}$ is the unit sphere. Moreover, $S^{d-1}_+\equiv S^{d-1}_+(v- 
v_1):=\{e\in S^{d-1}:e\circ(v-v_1)>0\}$, $v^\ast :=v-e\circ (v-v_1)\, e$, 
$v_1^\ast:=v_1+e\circ (v-v_1)\, e$ for $e\in S^{d-1}_+$ as well as $v,v_1 
\in V$, and $de$ refers to the normalized surface measure on $S^{d-1}_+$.  
\item[(v)] The {\it collision kernel} $B$ is assumed to be non-negative 
and continuous on $V\times V\times S^{d-1}$, symmetric in $v$ and $v_1$, 
and satisfying $B(v^\ast,v_1^\ast,e)=B(v,v_1,e)$ for all $v,v_1\in V$ and 
$e\in S^{d-1}_+$ for which $(v^\ast ,v_1^\ast)\in V\times V$. Furthermore, 
there is a constant $b\in (0,\infty)$ such that with $\varphi(v_1) \equiv 
\varphi(v;v_1,e) :=-v+2v_1+(|v-v_1|^2/e\circ(v-v_1))\, e$, defined for 
$(v,v_1,e)\in \mathbb{R}^d \times \mathbb{R}^d \times S^{d-1}$, we have 
\begin{eqnarray*}
B(v,\varphi(v_1),e)\le 2^{1-d}b\cdot\left(\frac{v-v_1}{|v-v_1|}\circ e 
\right)^2
\end{eqnarray*} 
for all $(v,v_1,e)\in V \times V \times S^{d-1}_+$ such that $\varphi(v_1) 
\in V$.

\bigskip 

\begin{multicols}{2} 

\begin{tikzpicture}[scale=2.7]
\draw[color=black,scale=0.6] (0,1) node [left] {$v$} -- (0,0) node [left,above] {$v_1\quad \ $} -- (0,-1) node [left] {$-v+2v_1$}; 
\draw[color=red,scale=0.6,->] (0,-1) -- (1,-0.5) node [below] {$e$}; 
\draw[color=black,scale=0.6] (1,-0.5) -- (2,0) node [right,above] {$\varphi(v_1)$}; 
\draw[color=red,scale=0.6] (-1,0) -- (3,0) node [below,color=black] {$G$};
\node at (1.6,-0.75) {\small{\bf Fig.~5}};
\end{tikzpicture}

\vspace{0.5cm} 

{\small{
Fig.~5: This is the minimal requirement on $B$ in order to establish (2.8) 
below. The estimate (2.8) will imply boundedness of global solutions. 
Let $G$ denote the hyperplane through $v_1$ orthogonal to $v-v_1$. 
For fixed $v\in V$, $B(v,\varphi(v_1),e)$ is decaying at least quadratically 
with respect to $|\varphi(v_1)-v_1|$, from $v_1$ along any direction on $G$. 
}} 
\end{multicols} 
\end{itemize}
\begin{itemize}
\item[(vi)] $h_\gamma$ is a continuous function on $\overline{\Omega 
\times\Omega}$ which is non-negative and symmetric, and vanishes for $|r-y| 
\ge\gamma>0$.  
\end{itemize} 

Define $\sigma\equiv\sigma(e):=(v^\ast-v_1^\ast)/|v^\ast-v_1^\ast|$ and 
$\t B(v,v_1,\sigma):=B(v,v_1,e)\cdot\chi_{\{(v^\ast,v_1^\ast)\in V\times V\}} 
(v,v_1,e)$. Let $d\sigma$ denote the normalized surface measure on $S^{d-1}$ 
relative to the variable $\sigma$ and keep in mind that $de$ refers to the 
normalized surface measure on $S^{d-1}_+\equiv S^{d-1}_+(v-v_1)$ relative to 
the variable $e$. For given $v,v_1\in \mathbb{R}^d$ with $v\neq v_1$, the map  
$S^{d-1}_+\equiv S^{d-1}_+(v-v_1)\ni e\mapsto\sigma\in S^{d-1}\setminus\{(v- 
v_1)/|v-v_1|\}$ is bijective. Let us fix $v\in \mathbb{R}^d$ as well as $\sigma 
\in S^{d-1}$ and define, in addition to $\vp$ introduced in condition (v),  
$\psi(w)\equiv\psi(v,w,\sigma):=\vp(v,w,\sigma)+v-w$ where $w\in\{u\in {\mathbb 
{R}^d}:\sigma\circ(v-u)>0\}$. We note that $\vp$ is injective. Moreover, it 
holds that $v^\ast=\psi(v_1^\ast)\equiv\psi(v,v_1^\ast,\sigma)$ and $v_1=\vp(v, 
v_1^\ast,\sigma)$, the latter with Jacobian determinant 
\begin{eqnarray*} 
\left|\frac{dv_1}{dv_1^\ast}\right|=\frac{2^{d-1}|v-v_1^\ast|^2}{(\sigma 
\circ (v-v_1^\ast))^2}\, ,\quad d=2,3.
\end{eqnarray*} 

For the next calculation, we  take into consideration that, for given 
$\sigma\in S^{d-1}$ and $v\in V$, it holds that 
\begin{eqnarray}\label{2.6}
\left\{v_1^\ast=\frac{v+v_1}{2}-\frac{|v-v_1|}{2}\, \sigma:v_1\in\mathbb 
{R}^d,\, v_1\neq v^\ast\right\}=\{v_1^\ast\in \mathbb{R}^d:(v-v_1^\ast) 
\circ\sigma>0\}. 
\end{eqnarray} 
Furthermore, we note that, for fixed $\sigma\in S^{d-1}$ and $v\in V$,  
$\left\{\vp(v_1^\ast):v_1\in V\right\}$ is not necessarily a subset of $V$. 
Introduce 
\begin{eqnarray*} 
\t B_J(v,v_1,\sigma):=\t B(v,v_1,\sigma)\left|\frac{de}{d\sigma}\right| 
(v,v_1,\sigma) 
\end{eqnarray*} 
and extend $\t B_J$ from $V\times V\times S^{d-1}$ to $V\times \mathbb{R}^d 
\times S^{d-1}$ by zero. It turns out that with 
\begin{eqnarray*}
p_\gamma(r,v_1,t):=\int_{y\in\Omega}p(y,v_1,t)\, h_\gamma(r,y)\, dy 
\end{eqnarray*} 
it holds that 
\begin{eqnarray}\label{2.7}
&&\hspace{-.5cm}Q^+(p,p)\, (r,v,t):=\int_V\int_{S_+^{d-1}}B(v,v_1,e)\chi_{ 
\{(v^\ast,v_1^\ast)\in V\times V\}}p(r,v^\ast ,t)p_\gamma(r,v_1^\ast ,t)\, 
de\, dv_1\nonumber \\ 
&&\hspace{.5cm}=\int_{S^{d-1}}\int_{\{v_1^\ast:v_1\in V\}}\left(\frac{ 
2^{d-1}\t B_J(v,\vp(v_1^\ast),\sigma)\, |v-v_1^\ast|^2}{(\sigma\circ(v- 
v_1^\ast))^2}\right)p(r,\psi(v_1^\ast),t)\, p_\gamma(r,v_1^\ast,t)\, d 
v_1^\ast\, d\sigma\nonumber \\ 
&&\hspace{.5cm}=\int_V\int_{S^{d-1}_+}\left(\frac{2^{d-1}\t B_J(v,\vp(v, 
v_1,\sigma),\sigma)\, |v-v_1|^2}{(\sigma\circ(v-v_1))^2}\right)p(r,\psi 
(v,v_1,\sigma),t)\, d\sigma\, p_\gamma(r,v_1,t)\, dv_1 \nonumber \\ 
\end{eqnarray} 
where, for the last equality sign, we have taken into consideration the 
identity (\ref{2.6}). For $p(\cdot,\cdot,t)\ge 0$ and $\|p(\cdot,\cdot,t) 
\|_{L^1(\Omega\times V)}\le 1$ it follows from (\ref{2.7}) and condition 
(v),that 
\begin{eqnarray}\label{2.8}
&&\hspace{-.5cm}Q^+(p,p)\, (r,v,t)\vphantom{\int}\nonumber \\ 
&&\hspace{.5cm}\le \int_V\int_{S^{d-1}_+}\left(\frac{2^{d-1}B(v,\vp(v,v_1,e) 
,e)\, |v-v_1|^2}{(e\circ(v-v_1))^2}\right)\, de\, p_\gamma(r,v_1,t)\, dv_1 
\cdot\| p(\cdot,\cdot,t)\|_{L^\infty(\Omega\times V)}\nonumber \\ 
&&\hspace{.5cm}\le b\, \|h_\gamma\|\cdot\| p(\cdot,\cdot,t)\|_{L^\infty(\Omega 
\times V)}\, . \vphantom{\int}
\end{eqnarray} 
\br{1} Denote by $\alpha\equiv\alpha(v,v_1,\sigma)\in(0,\pi]$ the angle 
between $\sigma$ and $v-v_1$. From (\ref{2.7}) we obtain the following 
{\it Carleman type representation} of the collision operator 
\begin{eqnarray*} 
&&\hspace{-.5cm}Q(p,p)\, (r,v,t)=Q^+(p,p)\, (r,v,t)-Q^-(p,p)\, (r,v,t) 
\vphantom{\int} \\ 
&&\hspace{.5cm}=\int_V\int_{S^{d-1}_+}\left(\frac{2^{d-1}\t B_J(v,\vp(v, 
v_1,\sigma),\sigma)}{\cos^2\alpha(v,v_1,\sigma)}\right)p(r,\psi(v,v_1, 
\sigma),t)\, d\sigma\, p_\gamma(r,v_1,t)\, dv_1 \\ 
&&\hspace{1.0cm}-\int_V\int_{S^{d-1}}\t B_J(v,v_1,\sigma)\, p(r,v,t)\, d 
\sigma\, p_\gamma(r,v_1,t)\, dv_1 
\end{eqnarray*} 
where $Q^+$ and $Q^-$ correspond to the two integrals. This representation 
of the collision operator can alternatively be obtained by first adjusting 
Lemma 7 of \cite{GBV09} to the particular form of the items in (iv),(v),(vi). 
In a second step one has to substitute the inner integral with respect to 
the canonical map (using the velocity symbols of \cite{GBV09}) $S^{d-1}_+(v- 
v_\ast')\mapsto\{v':(v'-v)\circ(v_\ast'-v)=0\}$. 
\er 
\br{2} 
For the physical and mathematical background of spatial smearing we 
refer to \cite{CPW98}, \cite{Lo18}, and the references therein. In 
particular, the estimates (\ref{2.8}) and 
\begin{eqnarray*}
\left\|Q(q_1+q_2,q_1-q_2)\right\|_{L^1(\Omega\times V)}\le 2\|h_\gamma\| 
\|B\|\left\|q_1+q_2\right\|_{L^1(\Omega\times V)}\left\|q_1-q_2\right\|_{ 
L^1(\Omega\times V)} 
\end{eqnarray*}  
are important for the analysis in the paper. They rely on spacial 
smearing. In Lemma \ref{Lemma5.1} and Proposition \ref{Proposition5.2} 
the latter estimate is used to prepare the application of Banach's fixed 
point theorem which provides local existence and uniqueness. In Proposition 
\ref{Proposition5.3} the latter estimate yields local Lipschitz continuity 
of the map $q\mapsto Q(q,q)$ in $L^1(\Omega\times V)$, which implies, via 
Lemmata \ref{Lemma5.4} and \ref{Lemma5.5}, global existence and uniqueness 
in Theorems \ref{Theorem5.6}, \ref{Theorem5.8}, and Corollary 
\ref{Corollary5.7}. 

Moreover, the estimate (\ref{2.8}) results in upper bounds in Theorems 
\ref{Theorem5.6}, \ref{Theorem5.8}, and Corollary \ref{Corollary5.7}. 
\er
Let us turn to the last of the global conditions, namely condition (vii) 
below. For this, recall the structure of $\partial\Omega$ from the beginning 
of this section. For $(y,v)\in\overline{\Omega}\times V$ introduce 
\begin{eqnarray*}
T_\Omega\equiv T_\Omega(y,v):=\inf\{s>0:y-sv\not\in\Omega\}  
\end{eqnarray*} 
and 
\begin{eqnarray}\label{2.9}
y^-\equiv y^-(y,v):=y-T_\Omega(y,v)v\, . 
\end{eqnarray} 
For all $(y,v)\in\bigcup_{i=1}^{n_\partial}\Gamma_i\times V=\partial^{ 
(1)}\Omega\times V$ with $v\circ n(y)\le 0$ as well as $y^-(y,R_y(v))\in 
\partial^{(1)}\Omega$, let 
\begin{eqnarray*}
\sigma(y,v):=\left(y^-(y,R_y(v)),R_y(v)\right)\, . 
\end{eqnarray*} 
We note that there exist $m_\partial\equiv m_\partial(\sigma)\in\mathbb{N}$ 
and mutually disjoint sets $G_i\subseteq\{(r,v):r\in\partial^{(1)}\Omega, 
\ v\in V,\ v\circ n(r)\le 0\}$, $i\in\{1,\ldots ,m_\partial\}$, satisfying 
\begin{eqnarray*} 
\bigcup_{i=1}^{m_\partial}\overline{G_i}=\overline{\{(r,v):r\in 
\partial^{(1)}\Omega,\ v\in V,\ v\circ n(r)\le 0\}} 
\end{eqnarray*} 
such that the following holds. {\it For each $i\in\{1,\ldots ,m_\partial 
\}$ there exists $j\in\{1,\ldots ,m_\partial\}$ such that $\sigma$ maps 
$G_i$ bijectively and continuously to $G_j$.} In this sense, there is 
an inverse of $\sigma$ denoted by $\sigma^{-1}$.

In this way we understand $\sigma$ as a map defined for $(y,v)\in 
\partial^{(1)}\Omega\times V$ with $v\circ n(y)\le 0$ and $y^-(y,R_y(v)) 
\in\partial^{(1)}\Omega$. For its $k$-fold composition $\sigma^k(y,v)= 
(y^{(k)}(y,v),v^{(k)}(y,v))\equiv (y^{(k)},v^{(k)})$, $k\in\mathbb{N}$, 
we suppose also that by iteration $y^{(j)}\in\partial^{(1)}\Omega$, $j\in 
\{1,\ldots ,k\}$, and we indicate this by the phrase ``\,for a.e. $y\in 
\partial^{(1)}\Omega\, $". Set $(y^{(0)},v^{(0)}):=(y,v)$. Likewise 
introduce $\sigma^{-k}(y,v)=(y^{(-k)}(y,v),v^{(-k)}(y,v))\equiv (y^{(-k)},
v^{(-k)})$ for $k\in\mathbb{N}$ by $\sigma^{-k}:=(\sigma^{-1})^k$ and use 
the phrase ``\,for a.e. $y\in\partial^{(1)}\Omega\, $" accordingly. The 
following condition on the shape of $\Omega$ will play an important role. 
\begin{itemize} 
\item[(vii)] There exist $k^{(0)}\in\mathbb{N}$ and $\sigma_{min}>0$ 
such that for a.e. $(y,v)\in\partial^{(1)}\Omega\times V$ with $v\circ 
n(y)\le 0$ and all $j\in\mathbb{Z}_+$ we have 
\begin{eqnarray*}
\sigma_{min}\le\left|y^{(j)}-y^{(j+1)}\right|+\left|y^{(j+1)}-y^{(j+2)} 
\right|+\ldots +\left|y^{(j+k^{(0)}-1)}-y^{(j+k^{(0)})}\right|\, .
\end{eqnarray*} 
\end{itemize} 

\bigskip\bigskip

\newcommand{\polygon}[2]{%
  let \n{len} = {2*#2*tan(360/(2*#1))} in
 ++(0,-#2) ++(\n{len}/2,0) \foreach \x in {1,...,#1} { -- ++(\x*360/#1:\n{len})}} 
 
\hspace{0.5cm} 
\fbox{

\begin{tikzpicture}[scale=3.2]
\useasboundingbox (0.07,0.035) rectangle (1.93,0.965);
      \tikzmath{ 
        % initial value, max value, next bounce value, speed
        \x = 0.4; \maxx=2; \bx=\maxx; \dx = 6;   %1;
        \y = 0; \maxy=1; \by=\maxy; \dy =  1;  %1.41421356237;
        for \i in {0,...,1}{
          % save values
          \sx = \x; \sy = \y;
          % time to next bounce
          \tx = (\bx-\x)/\dx;
          \ty = (\by-\y)/\dy;
          if \tx < \ty then { % if bounce on x before y
            \x = \bx;
            \bx = \maxx-\x;
            \dx = -\dx;
            \y = \y + \tx*\dy;
          } else { % if bounce on y before x
            \y = \by;
            \by = \maxy-\y;
            \dy = -\dy;
            \x = \x + \ty*\dx;
          };
          {\draw[color=red,thick] (\sx,\sy) -- (\x,\y);};
        }; % end for
      };
      \tikzmath{ 
        % initial value, max value, next bounce value, speed
        \x = 0; \maxx=2; \bx=\maxx; \dx = 1;   %1;
        \y = 0.9; \maxy=1; \by=\maxy; \dy =  5; %1.41421356237;
        for \i in {0,...,1}{
          % save values
          \sx = \x; \sy = \y;
          % time to next bounce
          \tx = (\bx-\x)/\dx;
          \ty = (\by-\y)/\dy;
          if \tx < \ty then { % if bounce on x before y
            \x = \bx;
            \bx = \maxx-\x;
            \dx = -\dx;
            \y = \y + \tx*\dy;
          } else { % if bounce on y before x
            \y = \by;
            \by = \maxy-\y;
            \dy = -\dy;
            \x = \x + \ty*\dx;
          };
          {\draw[color=blue,thick] (\sx,\sy) -- (\x,\y); };
        }; % end for
      };
        \tikzmath{ 
        % initial value, max value, next bounce value, speed
        \x = 0; \maxx=2; \bx=\maxx; \dx = 6;   %1;
        \y = 0.8; \maxy=1; \by=\maxy; \dy =  1; %1.41421356237;
        for \i in {0,...,1}{
          % save values
          \sx = \x; \sy = \y;
          % time to next bounce
          \tx = (\bx-\x)/\dx;
          \ty = (\by-\y)/\dy;
          if \tx < \ty then { % if bounce on x before y
            \x = \bx;
            \bx = \maxx-\x;
            \dx = -\dx;
            \y = \y + \tx*\dy;
          } else { % if bounce on y before x
            \y = \by;
            \by = \maxy-\y;
            \dy = -\dy;
            \x = \x + \ty*\dx;
          };
          {\draw[color=green,thick] (\sx,\sy) -- (\x,\y);};
        }; % end for 
      };
      \tikzmath{ 
        % initial value, max value, next bounce value, speed
        \x = 1.9; \maxx=2; \bx=\maxx; \dx = 1;   %1;
        \y = 1; \maxy=1; \by=\maxy; \dy = 6; %1.41421356237;
        for \i in {0,...,2}{
          % save values
          \sx = \x; \sy = \y;
          % time to next bounce
          \tx = (\bx-\x)/\dx;
          \ty = (\by-\y)/\dy;
          if \tx < \ty then { % if bounce on x before y
            \x = \bx;
            \bx = \maxx-\x;
            \dx = -\dx;
            \y = \y + \tx*\dy;
          } else { % if bounce on y before x
            \y = \by;
            \by = \maxy-\y;
            \dy = -\dy;
            \x = \x + \ty*\dx;
          };
          {\draw[color=black,thick] (\sx,\sy) -- (\x,\y);};
        }; % end for 
      };
\node[color=blue] at (-0.15,0.95) {$y^{(j)}$};
\node[color=green] at (-0.15,0.8) {\footnotesize$y^{(j)}$};    
\node[color=red] at (-0.15,0.6) {\footnotesize$y^{(j)}$};
\node[color=black] at (1.9,1.1) {\footnotesize$y^{(j)}$};
  
\node[color=blue] at (0.1,1.1) {\footnotesize$y^{(j+1)}$};
\node[color=green] at (1.2,1.1) {\footnotesize$y^{(j+1)}$};    
\node[color=red] at (2.2,0.28) {\footnotesize$y^{(j+1)}$}; 
\node[color=black] at (2.2,0.4) {\footnotesize$y^{(j+1)}$};

\node[color=blue] at (0.18,-0.15) {\footnotesize$y^{(j+2)}$};
\node[color=green] at (2.2,0.88) {\footnotesize$y^{(j+2)}$};    
\node[color=red] at (0.5,-0.15) {\footnotesize$y^{(j+2)}$}; 
\node[color=black] at (1.95,-0.15) {\footnotesize$y^{(j+2)}$}; 

\node[color=black] at (3.1,1.1) {\footnotesize$y^{(j)}$}; 
\node[color=black] at (4.15,0.8) {\footnotesize$y^{(j+1)}$};
\node[color=black] at (4.05,0.05) {\footnotesize$y^{(j+2)}$};

\node at (1.4,0.5) {$\Omega$};
\node at (3.6,0.5) {$\Omega$};

\draw[color=black,scale=0.6] (5.7,0.8) circle (1);
\draw[color=blue,scale=0.6] (5.7,0.8) \polygon{6.1}{0.866};
\draw[color=red,scale=0.6] (5.7,0.8) \polygon{8.1}{0.924};
\draw[color=green,scale=0.6] (5.7,0.8) \polygon{12}{0.966};
\draw[color=blue,scale=0.6] (5.7,0.8) \polygon{16}{0.981};
\draw[color=black,thick,scale=0.6] (5.7,0.8) \polygon{4.5}{0.765};
\end{tikzpicture} 
  }

\bigskip\bigskip\bigskip

\noindent
{\small{\hspace{1.2cm}{\bf Fig.~6} \hspace{8.1cm}{\bf Fig.~7}
\\ \\
Fig.~6: Let the shorter edges of the rectangle $\Omega$ be of length one. 
Then we have (vii) for $k_0=2$ and $\sigma_{min}=1$.   
Fig.~7: Condition (vii) fails if $\Omega$ is a disk. All distances $|y^{(j)} 
-y^{(j+1)}|$ originating from the same starting point $(y,v)\in\partial^{(1)} 
\Omega\times V$ with $v\circ n(y)\le 0$ are equal but can be arbitrarily small. }} 

\bigskip 

\noindent
Condition (vii) is primarily used in the proof of Theorem 
\ref{Theorem5.8} and the subsequent Lemma \ref{Lemma3.1}. This lemma is 
fundamental for the whole paper. Since condition (vii) appears again 
in quite another but more complex context in Subsection \ref{sec:4:1}, we 
postpone the discussion of its verifiability to Remark \ref{10} below.

\section{Basic Properties of the Knudsen Type Semigroup}\label{sec:3} 
\setcounter{equation}{0}

Suppose the global conditions (i)-(iii) and (vii). It will always be 
clear from the context what we understand by the Borel $\sigma$-algebra 
$\mathcal B$ over a subset $E$ of some Euclidean space. The terms Lebesgue 
measure or surface measure or canonical measure on $(E,\mathcal{B})$ refer 
to the measure corresponding to the absolute value of the respective volume 
form. 

In the paper we will work with spaces $L^1(E)$. More precisely, 
let $L^1(E)$ denote the space of all (equivalence classes of) measurable 
functions on $E$ which are absolutely integrable with respect to the 
Lebesgue measure or surface measure or canonical measure on $E$. 

Moreover, it will also always be clear from the context of the actual 
section, subsection, statement, etc. whether we are concerned with a 
space of complex valued functions, or if it is sufficient to deal with 
a space of real valued functions. By the {\it normalization condition} 
of (ii) and the {\it boundary conditions} in (iii), $S(t)$ maps $L^1 
(\Omega\times V)$ linearly to $L^1(\Omega\times V)$ with operator norm 
one.

\subsection{The Knudsen Type Semigroup Preserves a.e. Boundedness}
\label{sec:3:1} 

This subsection is entirely devoted to Lemma \ref{Lemma3.1} below. In 
its proof we will apply \cite{CPW98}, Theorem 2.1. It says that the 
{\it Knudsen} semigroup $S_0(t)$, $t\ge 0$, in $L^1(\Omega\times V)$ 
given by (\ref{2.5}) for $\omega=0$, initial condition $S_0(0)p_0= 
p_0\in L^1(\Omega\times V)$, and boundary conditions
\begin{eqnarray*}
(S_0(t)p_0)(r,v)=J(r,t)(S_0(\cdot)\, p_0)M(r,v)\, ,\quad t>0, 
\end{eqnarray*} 
for all $(r,v)\in\partial\Omega\times V$ with $v\circ n(r)\le 0$, 
admits a unique non-negative stationary element $\overline{g}_0\in 
L^1(\Omega\times V)$ with $\|\overline{g}_0\|_{L^1(\Omega\times V) 
}=1$. In other words, it holds that $S_0(t)\overline{g}_0=\overline 
{g}_0$ for all $t\ge 0$. Moreover, for any $\eta>0$ there exists 
$T_0(\eta)>0$ such that, for any $t\ge T_0(\eta)$ and for any 
probability density $f$ on $\Omega\times V$, we have $\left\|S_0(t) 
f-\overline{g}_0\right\|_{L^1(\Omega\times V)}\le\eta$. 

The following remarks explain in which sense the setup of 
\cite{CPW98}, Theorem 2.1 and its proof, is compatible with our 
framework. 
\br{3} Recall the hypotheses on $\partial\Omega$ from the beginning 
of Section \ref{sec:2}. In particular, recall the notion of   
\begin{eqnarray*} 
\partial^{(1)}\Omega=\bigcup_{i=1}^{n_{\partial}}\Gamma_i\, . 
\end{eqnarray*} 
Theorem 2.1 of \cite{CPW98} is formulated for a domain $\Omega$ with 
``sufficiently smooth" boundary, which means that it satisfies Assumption  
A of its proof. The proof is simplified to the case of a convex physical 
space $\Omega$. 

Rewriting this proof under our hypotheses on the boundary $\partial 
\Omega$, we have to replace the boundary $\partial\Omega$ there, with 
$\partial^{(1)}\Omega$. In addition we have to adjust the assumptions 
A and B within the proof of Theorem 2.1 of \cite{CPW98} to our specific 
situation. However everything from (A.11) in \cite{CPW98} on, can be 
taken over with just minor changes, see Remark \ref{4} below. 
%\medskip

Let $r_\Omega:=\sup\{|x-x_1|:x,x_1\in\partial\Omega,\, \{\alpha x+ 
(1-\alpha) x_1\, ,\ \alpha\in (0,1)\}\subset\Omega\}$. Transition 
densities as elements of $L^1(\Omega \times V)$ do not exist for 
the Knudsen semigroup $S_0(t)$ for small $t$. Keeping in mind the 
diffusive boundary conditions on $S_0$, it is well known that for 
$V=\{v\in \mathbb {R}^d:0<v_{min}<|v|<v_{max}<\infty\}$, transition 
densities as elements of $L^1(\Omega \times V)$ exist from time $(d+1) 
r_\Omega /v_0$ on where $d$ denotes the dimension of $\Omega$. For 
composite boundary conditions consisting of a specular reflective 
and a diffusive part, we have from time $(d+1)r_\Omega /v_0$ on
\begin{eqnarray*}
&&\hspace{-.5cm}S(t)f(x,v)=\int f(x',v') P_{t}(dx',dv',x,v) \\ 
&&\hspace{.5cm}=\int f(x',v') P^a_{t}(x',v',x,v)\, dx'\, dv'+\int 
f(x',v')P^s_{t}(dx',dv',x,v)
\end{eqnarray*}
for some nonnegative measurable function $P^a_{t}(\cdot ,\cdot ,x,v)$ 
on $\Omega \times V$ and some measure $P^s_{t}(\, \cdot\, ,x,v)$ on 
$\Omega\times V$ singular to the Lebesgue measure on $\Omega\times V$, 
where $P^s_{t} (\Omega \times V ,x,v) \stack {t\to\infty} {\lra}0$, 
$x\in\Omega$, $v\in V$. 
\medskip

In particular, we have to verify an adjusted version of (A.3) in 
\cite{CPW98} for a not necessarily convex $\Omega$ satisfying the 
conditions of Section 2. Expressed in terms of the present paper, 
(A.3) in \cite{CPW98} is related to the following. There is a time 
$t_0>0$ such that the Knudsen semigroup $S_0(t)$, $t\ge 0$, has the 
representation 
\begin{eqnarray*}
S_0(t_0)f=\int_{\Omega\times V}P_{t_0}(x',v';\cdot,\cdot)f(x',v')\, 
dx'\, dv'\, ,\quad f\in L^1(\Omega\times V), 
\end{eqnarray*} 
for some non-negative $P_{t_0}(x',v';\cdot,\cdot)\in L^1(\Omega\times 
V)$, $(x',v')\in\Omega\times V$, and the following holds. There exists 
$\gamma>0$ and $\beta>0$ such that 
\begin{itemize}
\item[(A.3)] $\qquad\displaystyle\inf_{x',v'}\inf_{(x,v)\in\mathcal{M} 
(\beta)}P_{t_0}(x',v';x,v)\ge\gamma$
\end{itemize}
where with $y^-=y-T_\Omega(y,v)v$ as in (\ref{2.9}), 
\begin{eqnarray*}
\mathcal{M}(\beta):=\left\{(y,v)\in\Omega\times V:y^-\in\partial^{(1)} 
\Omega,\ |v\circ n(y^-)|M(y^-,v)\ge\beta\right\}\, . 
\end{eqnarray*} 

For our purposes, Assumption A of \cite{CPW98} has to be reformulated as 
follows. 
\medskip

\nid 
(A') {\it There exist $\ve>0$ and $n_0\in\mathbb{N}$ such that for all 
$y,y'\in\partial^{(1)}\Omega$, there exist $y_1,\ldots,y_n\in\partial^{(1)} 
\Omega$ for some $n\le n_0$ in the following way. The straight connections 
between the points $y$ and $y_1$, between $y_j$ and $y_{j-1}$ for $j\in\{2, 
\ldots,n\}$, and between the points $y'$ and $y_n$, excluding $y,y',y_1, 
\ldots ,y_n$, lie entirely in $\Omega$. Furthermore, it holds that} 
\begin{itemize}
\item[(A.4)] $\min(|y_1-y|,|y_2-y_1|,\ldots ,|y_n-y_{n-1}|,|y'-y_n|)\ge\ve$ 
{\it and,} 
\item[(A.5)] {\it denoting $e(z,z_1):=(z_1-z)/|z_1-z|$, $z,z_1\in\partial^{ 
(1)}\Omega$, we have}
\begin{eqnarray*} 
|e(y,y_1)\circ n(y_1)|\ge\ve\, ,\quad |e(y,y_1)\circ n(y)|\ge\ve 
\end{eqnarray*} 
{\it and }
\begin{eqnarray*} 
|e(y_2,y_1)\circ n(y_2)|\ge\ve\, ,\ \ldots\, ,\ |e(y_n,y_{n-1})\circ 
n(y_n)|\ge\ve\, ,\quad |e(y_n,y')\circ n(y')|\ge\ve.
\end{eqnarray*} 
\end{itemize} 

It turns out that in order to take over the ideas of the 
first part of the proof of Theorem 2.1 in [8] it is 
sufficient to just consider $P^a_{t}(\cdot ,\cdot ,x,v)$ 
and to ignore $P^s_{t}(\, \cdot\, ,x,v)$. In particular, 
$(A')$ guarantees the following. Two  arbitrary points 
$y',y\in\partial^{(1)}\Omega$ can be connected trough a 
sequence of at most $n_0$ points $y_1,\ldots,y_n \in 
\partial^{(1)} \Omega$ such that the straight connections 
between the points $y$ and $y_1$, between $y_{j-1}$ and 
$y_{j}$ for $j\in \{2,\ldots ,n\}$, and between $y_n$ and 
$y'$, excluding $y,y',y_1, \ldots ,y_n$, lie entirely in 
$\Omega$. By (A.4) of condition $(A')$ the individual 
lengths of these connections are at least $\varepsilon >0$. 
Following the path from $y'$ to $y$ by straight connections 
in the order $y',y_n, y_{n-1}\ldots ,y_1,y,y_1,y,\ldots 
,y_1,y$ over a period of time of length $t_0$, where $4\vee 
(n_0+1)r_\Omega /v_{min}\le t_0 < 4\vee (n_0+1)r_\Omega 
/v_{min}+2$, there is just a bounded number of boundary 
hits possible. Assumption (A.5) of condition $(A')$ and the 
positive boundedness from below on the redistribution density 
$M$ in the diffusive part of the boundary conditions, cf. 
(ii),  imply now that there is a positive lower bound on 
$P_{t_0}(\cdot ,x,v)$ in the sense of the following 
adjustment of (A.3) in \cite{CPW98}, 
\begin{eqnarray*} 
\inf_{\vp\ge 0,\ \|\vp\|_{L^1(\Omega\times V)}=1}\ \inf_{(x,v) 
\in\mathcal{M}(\beta)}\int\vp(x',v') P_{t_0}(dx',dv',x,v)\ge 
\gamma\, .  
\end{eqnarray*} 

Assumption B in \cite{CPW98}, proof of Theorem 2.1, follows from 
our hypothesis (ii). 
\er 
\br{4} 
There are just a few minor change necessary in the proof of \cite{CPW98}, 
Theorem 2.1, to show the following. 
{\it There is a non-negative element $\overline{g}\in L^1(\Omega\times 
V)$ with $\|\overline{g}\|_{L^1(\Omega\times V)}=1$ which is stationary 
under the Knudsen type semigroup $S(t)$, $t\ge 0$, i.e. 
\begin{eqnarray}\label{3.1}
S(t)\overline{g}=\overline{g}\, ,\quad\mbox{\rm for all}\ t\ge 0.  
\end{eqnarray} 
Furthermore, for any $\eta>0$ there exists $T(\eta)>0$ such that, for 
any $t\ge T(\eta)$ and for any probability density $f$, we have} 
\begin{eqnarray}\label{3.2} 
\left\|S(t)f-\overline{g}\right\|_{L^1(\Omega\times V)}\le\eta\, . 
\end{eqnarray} 
The first change is the subsequent one. Using our notation, the two 
sentences after (A.8) in \cite{CPW98} have to be updated as follows. 
{\it The velocities $u=(y-y_1)/t^\ast$ and $-u$ allow us to go both 
directions from the boundary points $y_1$ and respectively $y$ with 
probability densities uniformly bounded from below. For this, we focus 
on the diffusive part of the boundary conditions. Indeed, according 
to (A.5) and (A.6), we have for $u$ and $y_1$}
\begin{eqnarray*} 
&&\hspace{-.5cm}(1-\omega)\, |u\circ n(y_1)|M(y_1,u) \\ 
&&\hspace{.5cm}= (1-\omega)|u|\left|\frac{u}{|u|}\circ n(y_1) 
\, \right|M(y_1,u)\ge (1-\omega)v_{min}\, \ve M_{\rm min}\, . 
\end{eqnarray*} 
The next change concerns the part between (A.11) and (A.18) in 
\cite{CPW98}. Since for the semigroup $S(t)$, $t\ge 0$, there are no 
transition densities, we replace $P_{t_0}(x',v';x,v)$ by 
\begin{eqnarray*} 
P_{t_0}^{(\alpha)}(x',v';x,v):=\frac{\int_\Omega\int_V\alpha^{-2d}\vp 
((x'-y)/\alpha,(v'-w)/\alpha)P_{t_0}(dy,dw;x,v)}{\int_\Omega\int_V 
\alpha^{-2d}\vp((x'-y)/\alpha,(v'-w)/\alpha)\, dy\, dw}
\end{eqnarray*} 
where $\vp$ is the $2d$-dimensional standard mollifier function and 
$\alpha\in (0,1)$. Note that $\int_\Omega\int_V P_{t_0}^{(\alpha)} 
(x',v';x,v)\, dx\, dv=1$, set $S^{(\alpha)}f(x,v):=\int_\Omega\int_V 
P_{t_0}^{(\alpha)}f(x',v')(x',v';x,v)\, dx\, dv$, and follow the proof 
in \cite{CPW98} from (A.11) until (A.18) with $P_{t_0}^{(\alpha)} 
(x',v';x,v)$ instead of $P_{t_0}(x',v';x,v)$. 

In this second part of our adjustment of the proof of Theorem 2.1 in 
[8], i.e. taking over the ideas and calculations from (A.11) through 
(A.18) in [8] for $P_{t_0}^{(\alpha)}$ instead of $P_{t_0}$, a certain 
joint representation of $S^{(\alpha)}f$ and $S^{(\alpha)}g$ is 
explicitly calculated. That is a Markov kernel whose marginals are 
$S^{(\alpha)}f$ and $S^{(\alpha)}g$. It is shown that the joint 
representation is exponentially decaying with respect to the number 
of time steps. 

Together with the adjusted result of the first part of the proof of 
Theorem 2.1 in [8] we arrive at the following exponential estimate, 
$\|(S^{(\alpha)})^nf-(S^{(\alpha)})^ng\|_{L^1}\le (1-\ve)^n$. Since 
$(S^{(\alpha)})^nf\stack {\alpha\to 0}{\lra}S^nf$ and $(S^{(\alpha)} 
)^n\stack {\alpha\to 0}{\lra}$ $S^ng$ in $L^1$, we obtain (A.18) of 
\cite{CPW98}, i.e.~$\|S^nf-S^ng\|_{L^1}\le (1-\ve)^n$. From (A.18)
on the proof in \cite{CPW98} can be taken over, word for word. 
\er 

Taking into consideration the modification of the proof of Theorem 2.1 
in \cite{CPW98} according to Remark \ref{3}, by means of Remark \ref{4} 
we have verified the existence of a non-negative function $\overline{g} 
\in L^1(\Omega\times V)$ with $\|\overline {g}\|_{L^1(\Omega\times V)} 
=1$ and (\ref{3.1}), (\ref{3.2}). 
\begin{lemma}\label{Lemma3.1} 
Let $\omega\in (0,1)$ and $p_0\in L^\infty(\Omega\times V)$. There are 
finite real numbers $p_{0,{\rm min}}$ and $p_{0,{\rm max}}$ such that 
\begin{eqnarray}\label{3.3}
p_{0,{\rm min}}\le S(t)\, p_0\le p_{0,{\rm max}}\quad\mbox{\rm a.e. on } 
\Omega\times V
\end{eqnarray} 
for all $t\ge 0$. In particular, if $p_0\ge 0$ a.e. and $\|1/p_0\|_{L^\infty 
(\Omega\times V)}<\infty$ then there exists $p_{0,{\rm min}}>0$. In addition, 
$\|\bar{g}\|_{L^\infty(\Omega\times V)}<\infty$ and $\|1/\bar{g}\|_{L^\infty 
(\Omega\times V)}<\infty$. 
\end{lemma}
Proof. {\it Step 1 } Without loss of generality, we may suppose 
that $\overline{g}:\Omega\times V\to [0,\infty]$ is defined everywhere 
on $\Omega\times V$ such that for every $(r,v)\in\Omega\times V$ there 
are $a\equiv a(r,v)<0$ and $b\equiv b(r,v)>0$ with $r+av\in\partial 
\Omega$ and $r+bv\in\partial\Omega$, $\{r+cv:c\in (a,b)\}\subset\Omega$, 
and $\overline{g}(r+cv,v)=\overline{g}(r,v)$, $c\in (a,b)$. In addition, 
we even may suppose that $\overline{g}(r+a v,v)=\overline{g}(r+b v,v)= 
\overline{g}(r,v)$. For such a boundary point $y=r+bv$ we recall the 
notation $y^-(y,v)=r+av$, see (\ref{2.9}).  
\medskip

Let $L_w^1$ be the space of all equivalence classes of measurable functions 
$f$ defined on $(y,v)\in\partial\Omega\times V$ with $v\circ n(y)\le 0$ such 
that 
\begin{eqnarray*}
\|f\|_{L_w^1}:=\int_{y\in\partial^{(1)}\Omega}\int_{v\circ n(y)\le 0}|v\circ 
n(y)||f(y,v)|\, dv\, dy<\infty\, .
\end{eqnarray*}
Furthermore, the map $\mathcal{S}f(y,v):=f\left(y^-(y,R_y(v)),R_y(v)\right)= 
f(\sigma(y,v))$, $(y,v)\in\partial\Omega\times V$ with $v\circ n(y)\le 0$, 
is by 
\begin{eqnarray}\label{3.4}
&&\hspace{-.5cm}\|\mathcal{S}f\|_{L_w^1}=\int_{y\in\partial^{(1)}\Omega}\int_{v\circ 
n(y)\le 0}|v\circ n(y)|\left|f\left(y^-(y,R_y(v)),R_y(v)\right)\right|\, dv 
\, dy\nonumber \\ 
&&\hspace{.5cm}=\int_{y\in\partial^{(1)}\Omega}\int_{v_{min}}^{v_{max}}\int_{ 
S_+^{d-1}(n(y))}\alpha e\circ n(y)\left|f\left(y^-(y,R_y(-\alpha e)),R_y(- 
\alpha e)\right)\right|\, de\alpha^{d-1}\, d\alpha\, dy\nonumber \\ 
&&\hspace{.5cm}=\int_{y\in\partial^{(1)}\Omega}\int_{S_+^{d-1}}\int_{v_{min} 
}^{v_{max}}\alpha^d e\circ n(y)\left|f\left(y^-(y,R_y(-\alpha e)),R_y(-\alpha 
e)\right)\right|\, d\alpha\, de\, dy\nonumber \\ 
&&\hspace{.5cm}=\int_{y\in\partial^{(1)}\Omega}\int_{S_+^{d-1}}\int_{v_{min} 
}^{v_{max}}\alpha^d e\circ n(y)\left|f\left(y^-(y,\alpha e),\alpha e\right) 
\right|\, d\alpha\, de\, dy\nonumber \\ 
&&\hspace{.5cm}=\int_{y\in\partial^{(1)}\Omega}\int_{v\circ n(y)\le 0}|v\circ 
n(y)|\left|f\left(y^-(y,-v),-v\right)\right|\, dv\, dy\nonumber \\ 
&&\hspace{.5cm}=\int_{y^-\in\partial^{(1)}\Omega}\int_{v\circ n(y^-)\ge 0}v
\circ n(y^-)|f(y^-,-v)|\cdot\left(\frac{|v\circ n(y)|}{v\circ n(y^-)}\cdot\frac 
{dy}{dy^-}\right)\, dv\, dy^-\nonumber \\ 
&&\hspace{.5cm}=\int_{r\in\partial^{(1)}\Omega}\int_{w\circ n(r)\le 0}|w\circ 
n(r)||f(r,w)|\, dw\, dr=\|f\|_{L_w^1} 
\end{eqnarray}
a linear operator $\mathcal{S}:L_w^1\mapsto L_w^1$ with operator norm one. For this 
calculation we have used $e\circ n(y)=R_y(-e)\circ n(y)$ as well as $de=d 
R_y(-e)$ in order to obtain the fourth line from the third. 
\medskip

\nid 
{\it Step 2 } Next we aim to demonstrate that $\overline{g}\in L_w^1$. For 
this we denote $(\partial\Omega)_y:=\{r\in\partial^{(1)}\Omega:\{r+\alpha 
(y-r):\alpha\in (0,1)\}\subset\Omega\}$, $y\in\overline{\Omega}$, and 
$(\Omega)_r:=\{y\in\Omega:\{r+\alpha(y-r):\alpha\in (0,1)\}\subset\Omega\}$, 
$r\in\partial\Omega$. Furthermore, we introduce 
\begin{eqnarray*} 
\rho_d(r)=\int_{y\in (\Omega)_r}|y-r|^{1-d}\cdot\frac{n(r)\circ (r-y)} 
{|y-r|}\, dy\, ,\quad r\in\partial^{(1)}\Omega,  
\end{eqnarray*}
and observe that, by the piecewise smoothness of $\partial\Omega$, there 
is a constant $c_d>0$ only depending on $\Omega$ such that $c_d\le\rho_d 
(r)<\infty$ for all $r\in\partial\Omega$. Set $C_{v,d}:=(v_{max}^d-v_{min 
}^d)/d$ as well as $C_d:=(1-\omega)M_{\rm min}\cdot C_{v,d}$, and denote 
by $l_S$ the surface measure on $(S^{d-1},\mathcal{B}(S^{d-1}))$. 

We will write $v=\alpha e$ where $\alpha\in (v_{min}, v_{max})$ and $e 
\in S^{d-1}$. Here we mention that $y^-(y,\alpha e)\in\partial\Omega$ is 
independent of $\alpha\in(v_{min},v_{max})$ and therefore may appear 
as $y^-(y,\cdot\, e)$. We obtain 
\begin{eqnarray}\label{3.5}
&&\hspace{-.5cm}1=\int_\Omega\int_V\overline{g}(y,v)\, dv\, dy=\int_\Omega 
\int_V\overline{g}(y^-(y,v),v)\, dv\, dy\nonumber \\ 
&&\hspace{.5cm}\ge (1-\omega)M_{\rm min}\int_\Omega\int_VJ(y^-(y,v),\cdot) 
(\overline{g})\, dv\, dy\nonumber \\
&&\hspace{.5cm}=(1-\omega)M_{\rm min}\int_\Omega\int_{v_{min}}^{v_{max}} 
\alpha^{d-1}\int_{S^{d-1}}J(y^-(y,\alpha e),\cdot)(\overline{g})\, dl_S(e) 
\, d\alpha\, dy\nonumber \\ 
&&\hspace{.5cm}=C_d\int_\Omega\int_{S^{d-1}}J(y^-(y,\cdot\, e), \cdot) 
(\overline{g})\, dl_S(e)\, dy\nonumber \\ 
&&\hspace{.5cm}=C_d\int_{y\in\Omega}\int_{r\in(\partial\Omega)_y}|y-r|^{1-d} 
\cdot\frac{n(r)\circ (r-y)}{|y-r|}J(r,\cdot)(\overline{g})\, dr\, dy\nonumber 
 \\ 
&&\hspace{.5cm}=C_d\int_{\partial^{(1)}\Omega}J(r,\cdot)(\overline{g})\rho_d 
(r)\, dr\ge c_dC_d\int_{\partial\Omega}J(r,\cdot)(\overline{g})\, dr\, . 
\end{eqnarray}
Furthermore, we note that  
\begin{eqnarray}\label{3.6}
&&\hspace{-.5cm}J(y,\cdot)(\overline{g})=\int_{v\circ n(y)\ge 0}v\circ n 
(y)\overline{g}(y,v)\, dv\nonumber \\ 
&&\hspace{.5cm}=\omega\int_{v\circ n(y)\le 0}|v\circ n(y)|\overline{g}(y, 
R_y(v))\, dv+(1-\omega)\int_{v\circ n(y)\le 0}|v\circ n(y)|M(y,v)J(y,\cdot) 
(\overline{g})\, dv\nonumber \\ 
&&\hspace{.5cm}=\int_{v\circ n(y)\le 0}|v\circ n(y)|\overline{g}(y,v)\, dv 
\, ,\quad y\in\partial^{(1)}\Omega.   
\end{eqnarray} 
It follows now from (\ref{3.5}) and (\ref{3.6}) that
\begin{eqnarray}\label{3.7} 
\|\overline{g}\|_{L_w^1}=\int_{r\in\partial^{(1)}\Omega}\int_{v\circ n(r) 
\le 0}|v\circ n(r)|\overline{g}(r,v)\, dv\, dr=\int_{\partial\Omega}J(r, 
\cdot)(\overline{g})\, dr\le\frac{1}{c_dC_d}\, . 
\end{eqnarray} 
In other words, we have $\overline{g}\in L_w^1$. 
\medskip

\nid
{\it Step 3 } In this step we apply the results of Steps 1 and 2. According 
to (iii) we have the boundary conditions on $\overline{g}$
\begin{eqnarray}\label{3.8}
\overline{g}(y,v)=\omega\, \overline{g}(y,R_y(v))+(1-\omega)\, J(y,\cdot) 
(\overline{g})\cdot M(y,v)\, ,\quad y\in\partial^{(1)}\Omega,\ v\circ 
n(y)\le 0. 
\end{eqnarray} 
These boundary conditions on $\overline{g}$ can be rewritten as  
\begin{eqnarray*} 
&&\hspace{-.5cm}\overline{g}(y,v)=(1-\omega)M(y,v)J(y,\cdot)(\overline 
{g})+\omega\, \overline{g}(y,R_y(v))\vphantom{\dot{f}} \\ 
&&\hspace{.5cm}=(1-\omega)M(y,v)J(y,\cdot)(\overline{g})+\omega\, 
\overline{g}(y^-(y,R_y(v)),R_y(v))\vphantom{\dot{f}} \\ 
&&\hspace{.5cm}=(1-\omega)M(y,v)J(y,\cdot)(\overline{g})+\omega\, (\mathcal{S} 
\overline{g})(y,v)\vphantom{\dot{f}}\, ,\quad y\in\partial^{(1)}\Omega, 
\ v\circ n(y)\le 0. 
\end{eqnarray*}
Together with the just shown $\overline{g}\in L_w^1$ and (\ref{3.4}) 
the latter says that, among other things, that $MJ(\overline{g})\in L_w^1$. 
Therefore 
\begin{eqnarray}\label{3.9}
&&\hspace{-.5cm}\overline{g}=(1-\omega)\sum_{k=0}^\infty\omega^k\mathcal{S}^k(MJ 
(\overline{g}))\nonumber \\ 
&&\hspace{.5cm}=(1-\omega)MJ(\overline{g})+(1-\omega)\sum_{k=0}^\infty 
\omega^{k+1}\mathcal{S}^{k+1}(MJ(\overline{g}))
\end{eqnarray} 
a.e. on $\{(y,v)\in\partial^{(1)}\Omega\times V:v\circ n(y)\le 0\}$, where 
the infinite sums converge in $L_w^1$. 
\medskip

\nid
{\it Step 4 } The next two steps are devoted to upper bounds of $J(y, 
\cdot)(\overline{g})$. In the present one we still construct upper bounds 
depending on $y\in\partial\Omega$. In Step 5 below we will derive a bound
that is independent of $y$ and use it to show $\|\overline{g}\|_{L^\infty 
(\Omega\times V)}<\infty$.  

Taking into consideration that according to (ii), there exist constants 
$M_{\rm min},M_{\rm max}\in (0,\infty)$ such that $M_{\rm min}\le M(y,v) 
\le M_{\rm max}$, it turns out that $\mathcal{S}^kM$ is uniformly bounded for all 
$k\in\mathbb{Z}_+$. Moreover extending $J(\overline{g})$ to all $(y,v) 
\in\partial^{(1)}\Omega\times V$ with $v\circ n(y)\le 0$ by $J(\overline 
{g})(y,v):=J(y,\cdot)(\overline{g})$, from (\ref{3.4}) and $MJ(\overline 
{g})\in L^1_w$ we may conclude that $\mathcal{S}^kJ(\overline{g})\in L^1_w$ for all 
$k\in\mathbb{Z}_+$. It follows now from (\ref{3.9}) and (\ref{3.6}) 
that 
\begin{eqnarray}\label{3.10}
&&\hspace{-.5cm}J(y,\cdot)(\overline{g})=(1-\omega)\int_{v\circ n(y)\le 0} 
|v\circ n(y)|\sum_{k=0}^\infty\omega^k \mathcal{S}^{k+1}(MJ(\overline{g}))(y,v)\, dv 
\, ,\quad y\in\partial^{(1)}\Omega. \qquad
\end{eqnarray} 
Introducing 
\begin{eqnarray}\label{3.11}
\delta(k,y,\omega):=(1-\omega)\int_{w\circ n(y)\le 0}|w\circ n(y)|\mathcal{S}^kM(y,w) 
\, dw\, ,\quad k\in\mathbb{N}, 
\end{eqnarray} 
we obtain from (\ref{3.10}) that  
\begin{eqnarray*} 
&&\hspace{-.5cm}(1-\omega)\int_{w\circ n(y)\le 0}|w\circ n(y)|\mathcal{S}(MJ 
(\overline{g}))(y,w)\, dw \\ 
&&\hspace{.0cm}=(1-\omega)\frac{\delta(1,y,\omega)}{1-\delta(1,y, 
\omega)}\cdot\int_{v\circ n(y)\le 0}|v\circ n(y)|\sum_{k=1}^\infty 
\omega^{k}\mathcal{S}^{k+1}(MJ(\overline{g}))(y,v)\, dv\, ,\quad y\in\partial^{ 
(1)}\Omega. 
\end{eqnarray*} 
Inserting this in (\ref{3.10}) gives 
\begin{eqnarray}\label{3.12}
&&\hspace{-.5cm}J(y,\cdot)(\overline{g})=\frac{1-\omega}{1-\delta 
(1,y,\omega)}\int_{v\circ n(y)\le 0}|v\circ n(y)|\sum_{k=1}^\infty 
\omega^k\mathcal{S}^{k+1}(MJ(\overline{g}))(y,v)\, dv\, . \qquad
\end{eqnarray} 
Letting $k_0$ be the integer introduced in condition (vii) and 
iterating the calculations from (\ref{3.10}) to (\ref{3.12}) $k_0$ 
times we verify 
\begin{eqnarray*}
\theta(y)\equiv\theta(k_0,y,\omega):=\sum_{k=1}^{k_0}\omega^{k-1} 
\delta(k,y,\omega)\le 1 
\end{eqnarray*} 
as well as 
\begin{eqnarray}\label{3.13}
&&\hspace{-.5cm}J(y,\cdot)(\overline{g})=\frac{1-\omega}{1-\theta 
(y)}\int_{v\circ n(y)\le 0}|v\circ n(y)|\sum_{k=k_0}^\infty\omega^k 
\mathcal{S}^{k+1}(MJ(\overline{g}))(y,v)\, dv\, ,\quad y\in\partial^{(1)} 
\Omega.\qquad
\end{eqnarray} 
Similarly we also obtain $\theta(k_0,y,\omega)+\omega^{k_0}\delta 
(k_0+1,y,\omega)=\theta(k_0+1,y,\omega)\le 1$. By $M\ge M_{\rm min} 
>0$ we get $\mathcal{S}^{k_0+1}M\ge M_{\rm min}>0$ and therefore from 
(\ref{3.11})
\begin{eqnarray*}
(1-\omega)M_{\rm min}\int_{w\circ n(y)\le 0}|w\circ n(y)|\, dw\le 
\delta(k_0+1,y,\omega)\, ,\quad y\in\partial^{(1)}\Omega. 
\end{eqnarray*} 
Thus 
\begin{eqnarray}\label{3.14}
\theta(y)\equiv\theta(k_0,y,\omega)\le 1-\omega^{k_0}(1-\omega)M_{\rm 
min}\int_{w\circ n(y)\le 0}|w\circ n(y)|\, dw\, ,\quad y\in\partial^{ 
(1)}\Omega,
\end{eqnarray} 
where we observe that the right-hand side does not depend on $y$ and 
is smaller than one. In other words, (\ref{3.14}) implies the existence 
of $\kappa<1$ such that
\begin{eqnarray*}
\theta(y)\equiv\theta(k_0,y,\omega)\le\kappa<1\, ,\quad y\in\partial^{ 
(1)}\Omega. 
\end{eqnarray*} 

Furthermore according to (ii), there exists $0<M_{\rm max} <\infty$ 
such that $M(y,v)\le M_{\rm max}$ and hence $\mathcal{S}^{k+1}M(y,v)=M(\sigma^{k+1} 
(y,v))\le M_{\rm max}$ for $(y,v)\in\partial^{(1)}\Omega\times V$ with 
$v\circ n(y)\le 0$. From (\ref{3.13}) it follows now that 
\begin{eqnarray}\label{3.15}
&&\hspace{-.5cm}J(y,\cdot)(\overline{g})\le\frac{1-\omega}{1-\kappa} 
\sum_{k=k_0}^\infty\omega^k\int_{v\circ n(y)\le 0}|v\circ n(y)|\mathcal{S}^{k+ 
1}M(y,v)\cdot \mathcal{S}^{k+1}J(\overline{g})(y,v)\, dv\nonumber \\ 
&&\hspace{.5cm}\le\frac{1-\omega}{1-\kappa}\, M_{\rm max}\sum_{k=k_0 
}^\infty\omega^k\int_{v\circ n(y)\le 0}|v\circ n(y)|\, \mathcal{S}^{k+1}J 
(\overline{g})(y,v)\, dv\, \quad y\in\partial^{(1)}\Omega. 
\end{eqnarray} 

\nid
{\it Step 5 } In this step we show $\|\overline{g}\|_{L^\infty(\Omega 
\times V)}<\infty$. Recall the notation from the end of Section 
\ref{sec:2}. Let $(y,v)\in\partial^{(1)}\Omega\times V$ with $v\circ 
n(y)\le 0$. For $e:=-v/|v|\in S_+^{d-1}(n(y))$, let $r:=y^-(y,e)$, 
i.e. $e=(y-r)/|y-r|$. In the remainder of this step we suppose $r\in 
\partial^{(1)}\Omega$ and iteratively also $y^{(j)}\in\partial^{(1)} 
\Omega$, $j\in\mathbb{N}$, and we indicate this in the text by the 
phrase ``\,for a.e. $y\in\partial^{(1)}\Omega\, $". 

Furthermore, denote $e_k:=-v^{(k)}/|v^{(k)}|$, $k\in\mathbb{N}$. For 
a.e. $y\in\partial^{(1)}\Omega\, $ introduce the distance 
\begin{eqnarray*}
|y-y^{(k)}|_k:=\sum_{j=1}^k|y^{(j)}-y^{(j-1)}|\, ,\quad k\in \mathbb{N}. 
\end{eqnarray*}
Note that 
\begin{eqnarray*}
\frac{dy^{(k)}}{dr}=\left.\frac{dy^{(k)}}{de_1}\right/\frac{dr}{de_1}= 
\left.\frac{dy^{(k)}}{de_1}\right/\frac{dr}{de}=\frac{|y-y^{(k)}|_k^{d 
-1}|e\circ n(r)|}{|y-r|^{d-1}|e_{k}\circ n(y^{(k)})|}\, ,\quad k\in 
\mathbb{N},  
\end{eqnarray*}
where 
\begin{eqnarray*}
\frac{dr}{de_1}=\frac{dr}{de}\cdot\frac{de}{de_1}=\frac{dr}{de} 
\end{eqnarray*}
by symmetry of $e$ and $-e_1$ about $n(y)$, and 
\begin{eqnarray*}
\frac{dy^{(k)}}{de_1}=\frac{|y-y^{(k)}|_k^{d-1}}{|e_{k}\circ n(y^{ 
(k)})|} 
\end{eqnarray*}
is usually motivated by means of a ray from $y$ in direction of $e_1$ 
as follows. The ray passes through $\Omega$ until it hits $y^{(1)}\in 
\partial^{(1)}\Omega$. Then $\Omega$ is reflected about the straight 
line ($d=2$) or plane ($d=3$) orthogonal to $n(y^{(1)})$ and containing 
$y^{(1)}$. In this way the ray passes through $k-1$ more consecutively 
reflected copies of $\Omega$.  

With these preparations in mind we obtain for a.e. $y\in\partial^{(1)} 
\Omega$ and $k\in\mathbb{N}$ 
\begin{eqnarray*}
&&\hspace{-.5cm}\int_{v\circ n(y)\le 0}|v\circ n(y)|\mathcal{S}^kJ(\overline{g}) 
(y,v)\, dv=\int_{v\circ n(y)\le 0}|v\circ n(y)|J(\overline{g})(y^{(k)} 
,v^{(k)})\, dv\nonumber \\ 
&&\hspace{.5cm}=\int_{S_+^{d-1}}\int_{v_{min}}^{v_{max}}\alpha^d e\circ 
n(y)J(y^{(k)},\cdot)(\overline{g})\, d\alpha\, de\nonumber \\ 
&&\hspace{.5cm}=C_{v,d+1}\, \int_{S_+^{d-1}}e\circ n(y)J(y^{(k)},\cdot) 
(\overline{g})\, de\nonumber \\ 
&&\hspace{.5cm}=C_{v,d+1}\, \int_{r\in\partial\Omega}\frac{e\circ n(y) 
\cdot |e\circ n(r)|}{|y-r|^{d-1}}J(y^{(k)},\cdot)(\overline{g})\, dr 
\nonumber \\ 
&&\hspace{.5cm}=C_{v,d+1}\, \int_{r\in\partial\Omega}\frac{e\circ n(y) 
\cdot |e\circ n(r)|}{|y-r|^{d-1}}\int_{w\circ n(y^{(k)})\le 0}|w\circ n 
(y^{(k)})|\overline{g}(y^{(k)},w)\, dw\, dr\nonumber \\ 
&&\hspace{.5cm}=C_{v,d+1}\, \int_{y^{(k)}\in\partial\Omega}\frac{e\circ 
n(y)\cdot |e\circ n(y^{(k)})|}{|y-y^{(k)}|_k^{d-1}}\int_{w\circ n(y^{(k 
)})\le 0}|w\circ n(y^{(k)})|\overline{g}(y^{(k)},w)\, dw\, dy^{(k)} 
\end{eqnarray*}
where we have applied (\ref{3.6}) in the second last line. By condition 
(vii) we have for all $k\ge k_0$
\begin{eqnarray*}
\frac{e\circ n(y)\cdot |e_k\circ n(y^{(k)})|}{|y-y^{(k)}|_k^{d-1}} 
\le\sigma_{min}^{1-d}  
\end{eqnarray*}
for a.e. $(y,v)\in\partial^{(1)}\Omega\times V$ with $v\circ n(y)\le 0$. 
Thus, for a.e. $y\in\partial^{(1)}\Omega$, 
\begin{eqnarray}\label{3.16}
&&\hspace{-.5cm}\int_{v\circ n(y)\le 0}|v\circ n(y)|\mathcal{S}^kJ(\overline{g}) 
(y,v)\, dv\nonumber \\ 
&&\hspace{.5cm}\le C_{v,d+1}\sigma_{min}^{1-d}\, \int_{y^{(k)}\in 
\partial\Omega}\int_{w\circ n(y^{(k)})\le 0}|w\circ n(y^{(k)})| 
\overline{g}(y^{(k)},w)\, dw\, dy^{(k)}\nonumber \\ 
&&\hspace{.5cm}=C_{v,d+1}\sigma_{min}^{1-d}\, \|\overline{g}\|_{L^1_w} 
\, ,\quad k\ge k_0.\vphantom{\int_{v\circ n(y)\le 0}}
\end{eqnarray}
It follows now from (\ref{3.7}), (\ref{3.15}), and (\ref{3.16}) 
that, for a.e. $y\in\partial^{(1)}\Omega$,
\begin{eqnarray*}
J(y,\cdot)(\overline{g})\le C_{v,d+1}\sigma_{min}^{1-d}M_{\rm max}\frac 
{\omega^{k_0}}{1-\kappa}\cdot\|\overline{g}\|_{L^1_w}=:C_J<\infty\, . 
\end{eqnarray*} 
Noting that this implies $\mathcal{S}^k(MJ(\overline{g}))(y,v)\le M_{\rm max}C_J$ 
for a.e. $(y,v)\in\partial^{(1)}\Omega\times V$ with $v\circ n(y)\le 0$ 
and $k\in \mathbb{Z}_+$, we may now conclude from (\ref{3.9}) and the 
first paragraph of Step 1 that 
\begin{eqnarray*}
\|\overline{g}\|_{L^\infty(\Omega\times V)}\le M_{\rm max}C_J<\infty\, . 
\end{eqnarray*} 
{\it Step 6 } Let us demonstrate $\|1/\overline{g}\|_{L^\infty 
(\Omega\times V)}<\infty$. For $1>\ve>0$ let $C_\Omega(y,\ve)$ denote 
the open cone in $\mathbb{R}^d$ with vertex $y\in\partial^{(1)}\Omega$ 
given by $\{x\in\mathbb{R}^d:((y-x)/|x-y|)\circ n(y)>\ve\}$. Furthermore, 
for $y\in\partial^{(1)}\Omega$ introduce 
\begin{eqnarray*}
(\partial\Omega)_{y,\ve}:=\left\{r\in(\partial\Omega)_y\cap C_\Omega 
(y,\ve):((r-y)/|r-y|)\circ n(r)>\ve\right\}\, . 
\end{eqnarray*} 
Using the notation $e:=(y-r)/|y-r|$, for every $1>\ve>0$ it holds that
\begin{eqnarray*} 
\inf_{y\in\partial^{(1)}\Omega,\, r\in (\partial\Omega)_{y,\ve}}\left\{ 
\frac{e\circ n(y)\cdot |e\circ n(r)|}{|y-r|^{d-1}}\right\}=:c_\ve >0\, .
\end{eqnarray*} 
As already explained in the beginning of Step 1 of this proof we have 
$\overline{g}(y^-(y,v),v)=\overline{g}(y,v)$ for $v\circ n(y)\ge 0$ and 
$y\in\partial^{(1)}\Omega$. Using the boundary conditions (\ref{3.8}) 
we obtain now 
\begin{eqnarray}\label{3.17}
&&\hspace{-.5cm}J(y,\cdot)(\overline{g})=\int_{v\circ n(y)\ge 0}v\circ 
n(y)\cdot\overline{g}(y,v)\, dv\nonumber \\ 
&&\hspace{.5cm}=\int_{v\circ n(y)\ge 0}v\circ n(y)\cdot \overline{g}(y^- 
(y,v),v)\, dv\nonumber \\ 
&&\hspace{.5cm}\ge\, (1-\omega)M_{\rm min}\int_{v\circ n(y)\ge 0}v\circ 
n(y)\cdot J(y^-(y,v),\cdot)(\overline{g})\, dv\nonumber \\ 
&&\hspace{.5cm}=(1-\omega)M_{\rm min}\int_{S_+^{d-1}(n(y))}\int_{v_{min} 
}^{v_{max}}\alpha^d e\circ n(y)\cdot J(y^-(y,\alpha e),\cdot)(\overline 
{g})\, d\alpha\, de\nonumber \\ 
&&\hspace{.5cm}=\, (1-\omega)M_{\rm min}C_{v,d}\int_{r\in (\partial\Omega 
)_y}\frac{|(r-y)\circ n(y)|\cdot (r-y)\circ n(r)}{|y-r|^{d+1}}J(r,\cdot) 
(\overline{g})\, dr\nonumber \\ 
&&\hspace{.5cm}\ge\, c_\ve(1-\omega)M_{\rm min}C_{v,d}\int_{r\in(\partial 
\Omega)_{y,\ve}}J(r,\cdot)(\overline{g})\, dr\, ,\quad y\in\partial^{(1)} 
\Omega.  
\end{eqnarray} 
Assume that there exist $y\in\partial\Omega$ and a sequence $y_k\in 
\partial^{(1)}\Omega$, $k\in \mathbb{N}$, with $y_k\stack {k\to\infty} 
{\lra}y$ as well as $J(y_k,\cdot)(\overline{g})\stack{k\to\infty}{\lra} 
0$. From (\ref{3.17}) it follows that for every $\ve>0$ and $\delta>0$ 
there is a $k\in\mathbb{N}$ such that 
\begin{eqnarray*} 
\delta\ge\int_{r\in(\partial\Omega)_{y_k,\ve}}J(r,\cdot)(\overline{g}) 
\, dr\, .  
\end{eqnarray*} 
Thus, $J(r,\cdot)(\overline{g})=0$ for a.e. $r\in(\partial\Omega)_y$. 
Plugging this in the left-hand side of (\ref{3.17}), and iterating 
the last conclusion, we may even state that $J(r,\cdot)(\overline{g}) 
=0$ for a.e. $r\in\partial^{(1)}\Omega$. Recalling (\ref{3.6}) and 
introducing $r^-(r,v)$ as $y^-(y,v)$ in (\ref{2.9}), it turns out that 
\begin{eqnarray*} 
&&\hspace{-.5cm}0=\int_{r\in\partial^{(1)}\Omega}J(r,\cdot)(\overline 
{g})\, dr=\int_{r\in\partial^{(1)}\Omega}\int_{v\circ n(r)\le 0}|v\circ 
n(r)|\overline{g}(r,v)\, dv\, dr \\ 
&&\hspace{.5cm}=\int_{v\in V}\int_{\{r\in\partial^{(1)}\Omega:v\circ n(r) 
\le 0\}}|v\circ n(r)|\overline{g}(r,v)\, dr\, dv \\ 
&&\hspace{.5cm}=\int_{v\in V}\int_{\{r\in\partial^{(1)}\Omega:v\circ n(r) 
\le 0\}}\frac{|v\circ n(r)|}{|r-r^-(r,-v)|}\int_0^{|r-r^-(r,-v)|}\overline 
{g}(r+\alpha v/|v|,v)\, d\alpha\, dr\, dv \\ 
&&\hspace{.5cm}=\int_V\int_\Omega\frac{1}{|r^-(r,v)-r^-(r,-v)|}\, \overline 
{g}(r,v)\, dr\, dv\, . 
\end{eqnarray*} 
In other words, the above assumption would lead to $\overline{g}=0$ a.e. 
on $\Omega\times V$. This proves the existence of a lower bound $c_J>0$
on $J(r,\cdot)(\overline{g})$, uniformly for a.e. $r\in\partial^{(1)} 
\Omega$. 

This yields $\mathcal{S}^k(MJ(\overline{g}))(y,v)\ge M_{\rm min}c_J>0$ for a.e. 
$(y,v)\in\partial^{(1)}\Omega\times V$ with $v\circ n(y)\le 0$ and $k\in 
\mathbb{Z}_+$. It follows now from (\ref{3.9}) and the first paragraph 
of Step 1 that 
\begin{eqnarray*}
\|\overline{g}\|_{L^\infty(\Omega\times V)}\ge M_{\rm min}c_J>0\, , 
\end{eqnarray*} 
i.e. $\|1/\overline{g}\|_{L^\infty(\Omega\times V)}<\infty$. 
\medskip 

\nid 
{\it Step 7 } Let $p_0\in L^\infty(\Omega\times V)$. By the result 
of Step 6 there exists $a>0$ with  $-a\, \overline{g}\le p_0\le a\, 
\overline{g}$ a.e. on $(\Omega\times V)$. Furthermore, in case of 
$p_0\ge 0$ a.e. and $\|1/p_0\|_{L^\infty(\Omega\times V)}<\infty$, there 
is $b>0$ such that $b\, \overline{g}\le p_0$ a.e. on $(\Omega\times 
V)$, see Step 5. This implies 
\begin{eqnarray*} 
-a\, \|\overline{g}\|_{L^\infty(\Omega\times V)}\le -a\, \overline 
{g}\le S(t)p_0\le a\, \overline{g}\le a\, \|\overline{g}\|_{L^\infty 
(\Omega\times V)}\quad\mbox{\rm a.e. on }(\Omega\times V)\ \mbox{\rm 
for all }t\ge 0 
\end{eqnarray*} 
and, if $p_0\ge 0$ a.e. and $\|1/p_0\|_{L^\infty (\Omega\times V)} 
<\infty$, also 
\begin{eqnarray*} 
0<b\cdot\mathop{\mathrm{ess~inf}}\limits_{(y,w)\in\Omega\times V}\, 
\overline{g}(y,w)\le b\, \overline{g}\le S(t)p_0\quad\mbox{\rm a.e. 
on }(\Omega\times V)\ \mbox{\rm for all }t\ge 0. 
\end{eqnarray*} 
The lemma follows. 
\qed 

\subsection{Construction of the Knudsen Type Semigroup}\label{sec:3:2}

Let us now introduce the concept of rays and paths. Let $(r,v)\in\Omega 
\times V$ or $(r,v)\in\partial^{(1)}\Omega\times V$ with $v\circ n(r)\ge 
0$. The map $[\tau_0,\tau_1)\ni\tau\mapsto r-\tau v$ with $0\le\tau_0< 
\tau_1\le \infty$ is called a {\it ray}. We mention that the time 
$T_\Omega\equiv T_\Omega(r,v)=\inf\{s>0:r-sv\not\in\Omega\}$ can be 
interpreted as the {\it first exit time from $\Omega$ of the ray} $[0, 
\infty)\ni\tau\mapsto r-\tau v$. 
\medskip

Suppose we are given $t\ge 0$ and $(r,v)\equiv(r_0,v_0)\in\Omega\times 
V$. If $T_{\Omega}(r,v)<t$ and $r_1:=r-T_{\Omega}(r,v)v\in\partial^{(1)} 
\Omega$ then from the point $r_1$ we simultaneously follow all rays $[0, 
\infty)\ni\tau\mapsto r_1-\tau v_1$ for which $v_1\circ n(r_1)\ge 0$ until 
$T_{\Omega}(r,v)+\tau=t\le T_{\Omega}(r,v)+T_{\Omega}(r_1,v_1)$ or these 
rays exit from $\Omega$ for the first time at some $r_2\in\partial\Omega$.

We start over, simultaneously from all $r_2\in\partial^{(1)}\Omega$, along 
all rays $[0,\infty)\ni\tau\mapsto r_2-\tau v_2$ for which $v_2\circ n(r_2) 
\ge 0$, and continue in this manner until $\sum_{l=0}^{m-1}T_{\Omega}(r_l, 
v_l)<t$ but $\sum_{l=0}^mT_{\Omega}(r_l,v_l)\ge t$ for some $m\in \mathbb{Z 
}_+$. Here we use the convention $\sum_0^{-1}=0$. We note that $m$ depends 
on $r_0,v_0,v_1,\ldots ,v_{m-1}$. 

Any ray $[0,\infty)\ni\tau\mapsto r_m-\tau v_m$ takes at time $\tau=t 
-\sum_{l=0}^{m-1}T_{\Omega}(r_l,v_l)$ the value of some  point $r_e\in 
\overline{\Omega}$. In this way, we have constructed sequences or rays 
for which we now consider only the restrictions to $[0,T_{\Omega}(r_0, 
v_0))$, $[0,T_{\Omega}(r_1,v_1))$, $\ldots\, $, $[0,T_{\Omega}(r_{m-1}, 
v_{m -1}))$, $[0,t-\sum_{l=0}^{m-1}T_{\Omega} (r_l,v_l)]$. 

To prepare the next definition we re-parametrize the rays to maps 
over consecutive intervals $\tau\in [0,T_{\Omega}(r_0,v_0))$, $\tau 
\in [T_{\Omega}(r_0,v_0),T_{\Omega}(r_0,v_0)+T_{\Omega}(r_1,v_1))$, 
$\ldots\, $, $\tau\in [\sum_{l=0}^{m-1}T_{\Omega}(r_l,v_l),t]$. It is 
important to understand that $\tau$ stands for reversed time from $t$ 
to zero, i.e. at time $t$ we have $\tau=0$ and at time zero we have 
$\tau=t$. 
\begin{definition}\label{Definition3.2}
{\rm Let $t>0$ and $(r,v)\in\Omega\times V$. \\ 
(a) A {\it path $\pi$ with time range} $[0,t]$ {\it pinned in at 
time $t$ in} $(r,v)$ is a finite collection of re-parametrized rays 
of the form 
\begin{eqnarray}\label{3.18} 
\left[\textstyle{\sum_{l=0}^{k-1}}T_{\Omega}(r_l,v_l),\textstyle 
{\sum_{l=0}^k}T_{\Omega}(r_l,v_l)\right)\ni\tau\mapsto r_k-\left(\tau 
-\textstyle{\sum_{l=0}^{k-1}}T_{\Omega}(r_l,v_l)\right)v_k\, , 
\end{eqnarray} 
$k=0,\ldots ,m-1$, such that $\textstyle{\sum_{l=0}^{m-1}}T_{\Omega} 
(r_l,v_l)<t$ as well as $\textstyle{\sum_{l=0}^m}T_{\Omega}(r_l,v_l) 
\ge t$ and 
\begin{eqnarray}\label{3.19} 
\left[\textstyle{\sum_{l=0}^{m-1}}T_{\Omega}(r_l,v_l),t\right]\ni 
\tau\mapsto r_m-\left(\tau-\textstyle{\sum_{l=0}^{m-1}}T_{\Omega}(r_l, 
v_l)\right)v_m\, ,  
\end{eqnarray} 
$m\in \mathbb{Z}_+$, where we use the convention $\sum_0^{-1}=0$. Here, 
we suppose
\begin{eqnarray*} 
v_k\circ n(r_k)\ge 0\, ,\quad r_k=r_{k-1}-T_{\Omega}(r_{k-1},v_{k-1}) 
v_{k-1}\in\partial^{(1)}\Omega\, ,\quad k=1,\ldots ,m, 
\end{eqnarray*} 
and $r_0:=r$ as well as $v_0:=v$. Furthermore, $v_e:=v_m$ and we suppose 
\begin{eqnarray*}
r_e:=r_m-\left(t-\textstyle{\sum_{l=0}^{m-1}}T_{\Omega}(r_l,v_l)\right) 
v_m\in\Omega\, . 
\end{eqnarray*} 
(b) For fixed $r,v,t$ as above, let {\boldmath${\pi}$}$(r,v,t)$ be 
the set of all paths $\pi$ with time range $[0,t]$ pinned at time $t$ 
in $(r,v)$. 
}
\end{definition} 

\bigskip\bigskip\bigskip\bigskip\medskip 

\hspace{3.5cm}
\fbox{ 

\begin{tikzpicture}[scale=3.2]

\useasboundingbox (0.06,0.023) rectangle (1.94,0.975);

         \draw[color=black,thick,->] (0.75,0.4) -- (-0.01,0.8);
         \draw[color=red,thick,->] (-0.01,0.8) -- (1,1);
         \draw[color=black,thick,->] (1,1) -- (1.3,-0.01);
         \draw[color=red,thick,->] (1.3,-0.01) -- (2,0.6);
         \draw[color=black,thick,->] (2,0.6) -- (1.7,1);
         \draw[color=red,thick,->] (1.7,1) -- (-0.01,0.3);
         \draw[color=black,thick,->] (-0.01,0.3) -- (0.5,0.2);
      
      \node[color=black] at (1.17,0.4) {\footnotesize$(r,v)\equiv (r_0,v_0)$};
      \node[color=black] at (0.75,0.18) {\footnotesize$(r_e,v_e)$};
      \node[color=red] at (-0.5,0.8) {\footnotesize$\displaystyle{(r_1,v_1) 
\ \mbox{\rm at }\atop \tau=T_\Omega(r_0,v_0)}$ };
      \node[color=black] at (1,1.2) {\footnotesize$\displaystyle{(r_2,v_2) 
\ \mbox{\rm at }\atop \tau=T_\Omega(r_0,v_0)+T_\Omega(r_1,v_1)}$ };
      \node[color=black] at (2.85,0.17) {\footnotesize$(r_e,v_e)$ at time 0, 
i.e.~$\tau=t$};
      \node[color=black] at (2.80,0.4) {\footnotesize$(r,v)$ at time $t$, 
i.e.~$\tau=0$};
      \node[color=black] at (-0.6,0.25) {\footnotesize$\displaystyle{(r_m,v_m) 
\ \mbox{\rm at }\atop \tau={\textstyle\sum_{l=0}^{m-1}}T_\Omega(r_l,v_l)}$};
      
\end{tikzpicture} 
 }

\bigskip 

\noindent
{\small{\hspace{4.2cm}{\bf Fig.~8}
\\ \\ 
Fig.~8: The figure displays a path $\pi$ with time range $[0,t]$ pinned in at 
time $t$ in $(r,v)$. The arrows correspond to the rays. Their lengths do not 
represent the absolute values of the associated velocities.}}

\bigskip

\br{5} Let $p_0$ be a probability density everywhere defined on 
$\Omega\times V$. A more intuitive explanation of the term {\it path 
$\pi$ with time range $[0,t]$ pinned at time $t$ in} $(r,v)$ is to 
follow $S(t)p_0$ along a path $\pi$ backward in time. This is 
from time $t$ to time 0, with the understanding that starting with 
$S(t)p_0(r,v)$, after $\tau\in [0,t]$ units backward in time we 
have arrived at some $S(t-\tau)p_0(r',v')$. In particular, we keep 
in mind the boundary conditions of (iii), together with (i). In 
this sense, a path $\pi$ consists of $m+1$ re-parametrized rays along 
each of which $S(t-\cdot)p_0$ is constant. These are 
\begin{itemize}
\item if $T_{\Omega}(r_0,v_0)\equiv T_{\Omega}(r,v)<t$, i.e. $m\ge 1$, 
the ray $r_0-\tau v_0$ from $r=r_0\in\Omega$ to $r_1\in\partial^{(1)} 
\Omega$ with velocity $v_0=v$, i.e. over the time range $\tau\in [0,
T_{\Omega}(r_0,v_0))$, 
\item if ${\sum_{l=0}^k}T_{\Omega}(r_l,v_l)<t$, the re-parametrized ray 
$r_k-\left(\tau-{\sum_{l=0}^{k-1}}T_{\Omega}(r_l,v_l)\right)v_k$ from 
$r_k\in\partial^{(1)}\Omega$ to $r_{k+1}\in\partial^{(1)}\Omega$ with 
velocity $v_k$, and therefore over the time range of (\ref{3.18}), 
$k=1,\ldots ,m-1$, if $m\ge 1$, 
\item and the re-parametrized ray $r_m-\left(\tau-\sum_{k=0}^{m-1}T_{ 
\Omega}(r_k,v_k)\right)v_m$ from $r_m\in\partial^{(1)}\Omega$ to $r_e 
\in\Omega$ with velocity $v_m=v_e$, over the time range of (\ref{3.19}), 
if $m\ge 1$. In case of $m=0$, the ray $r_0-\tau v_0$ from $r=r_0\in 
\Omega$ to $r_e\in\Omega$ with velocity $v_0=v=v_e$, over the time range 
$\tau\in [0,t]$. 
\end{itemize} 
\er

Let $p_0$ be a probability density everywhere defined on $\Omega 
\times V$ and let $(r,v,t)\in\Omega\times V\times [0,\infty)$. 
Introduce 
\begin{eqnarray*} 
t_k\equiv t_{\Omega,k}(r_0,v_0,\ldots ,v_{k-1}):=\sum_{l=1}^kT_{\Omega} 
(r_{l-1},v_{l-1})\, ,\quad 1\le k\le m\quad\mbox{\rm for } m\in\mathbb{N}.  
\end{eqnarray*} 
Recalling Definition \ref{Definition3.2} and Remark \ref{5} we may 
state that 
\begin{eqnarray}\label{3.20} 
&&\hspace{-.5cm}S(t)p_0(r,v)=\chi_{\{0\}}(m)p_0(r_e,v_0)+\sum_{k=1 
}^\infty\left[\omega\, \chi_{\{k\}}(m)S(t-t_1)p_0\left(r_1,R_{r_1} 
(v_0)\vphantom{l^1}\right)\vphantom{\int_{v_1\circ n(r_1)\ge 0}} 
\right.\nonumber \\
&&\hspace{.5cm}+\ (1-\omega)M(r_1,v_0)\int_{v_1\circ n(r_1)\ge 0}(v_1 
\circ n(r_1))\ldots\times\nonumber \\ 
&&\hspace{1.0cm}\times\left(\omega\, \chi_{\{k\}}(m)S(t-t_k)p_0\left( 
r_k,R_{r_k}(v_{k-1})\vphantom{l^1}\right)+(1-\omega)M(r_k,v_{k-1}) 
\times\vphantom{\int_{v_1\circ n(r_1)\ge 0}}\right.\nonumber \\ 
&&\hspace{1.0cm}\left.\left.\times\int_{v_k\circ n(r_k)\ge 0}(v_k 
\circ n(r_k))p_0(r_e,v_k)\vphantom{\int_{v_1\circ n 
(r_1)\ge 0}}\, \chi_{\{k\}}(m)\, dv_k\right)\, \ldots\, dv_1\right] 
\vphantom{\int_{v_1\circ n(r_1)\ge 0}}
\end{eqnarray} 
where we note that for fixed $(r,v,t)\equiv(r_0,v_0,t)\in\Omega\times 
V\times[0,\infty)$, the number $m$ and the point $r_e$ are functions 
of $v_1,\ldots ,v_m$. Furthermore $\chi_k(m)$ can be regarded as an 
abbreviation of $\chi_{\{m=k\}}(v_1,\ldots ,v_k)$. The sum in 
(\ref{3.20}) converges since any partial sum consists of non-negative 
items and is bounded by $S(t)p_0(r,v)$. 

We recall that {\boldmath${\pi}$}$(r,v,t)$ is the set of paths used in 
(\ref{3.20}). We observe that, in accordance with Definition 
\ref{Definition3.2} and Remark \ref{5}, none of the paths in 
{\boldmath${\pi}$}$(r,v,t)$ ever visits an edge or vertex of $\partial 
\Omega$. By definition, all paths contained in {\boldmath${\pi}$}$(r,v, 
t)$ take at time zero, which is $\tau=t$, some value $(r_e,v_e)\in\Omega 
\times V$. 

Furthermore, there is a corresponding stochastic Markov process $X_t$, $t\ge 
0$, which is given by the explicit construction of the paths in Definition 
\ref{Definition3.2}. The random mechanism is successively introduced to 
the paths whenever they hit the boundary $\partial^{(1)}\Omega$. They change 
direction according to the combination of non-random and random reflections 
given by the boundary conditions of hypothesis (iii), see (\ref{3.20}). 
Denoting by $Q_{r,v,t}$ the distribution over the paths $\pi\in$ 
{\boldmath${\pi}$}$(r,v,t)$ and by $(r_e(\pi),v_e(\pi))\in\Omega\times V$ the 
initial point of a single path $\pi\in${\boldmath ${\pi}$}$(r,v,t)$ at time 
zero, the function 
\begin{eqnarray}\label{3.21}
\Omega\times V\ni (r,v)\mapsto S(t)p_0(r,v)=\int_{\pi\in\mbox{\small 
\boldmath${\pi}$}(r,v,t)}p_0(r_e(\pi),v_e(\pi))\, Q_{r,v,t}(d\pi)
\end{eqnarray} 
represents the probability density with respect to the Lebesgue measure 
on $\Omega\times V$ of the stochastic process at time $t\ge 0$, provided 
it was started at time zero with the probability density $p_0$. 
\medskip

Below we will frequently refer to Remark \ref{5} and the following 
Remark \ref{6}, for example in the proofs of Lemma \ref{Lemma3.3} 
and Theorem \ref{Theorem4.5}. 
\medskip 

\noindent 
\br{6} The subsequent explicit representation of the semigroup $S(t)$, 
$t\ge 0$, is a consequence of (\ref{3.20}). Let $p_0$ be a probability 
density everywhere defined on $\Omega\times V$ and let $(r,v,t)\in\Omega 
\times V\times [0,\infty)$. We have 
\begin{eqnarray}\label{3.22} 
&&\hspace{-.5cm}S(t)p_0(r,v)=\chi_{\{0\}}(m)p_0(r_e,v_0)+\left(\omega 
\chi_{\{1\}}(m) p_0(r_e,R_{r_{1}}(v_0))\vphantom{\int_{v_1\circ n(r_1) 
\ge 0}}\right.\nonumber \\ 
&&\hspace{.5cm}\left.+(1-\omega)M(r_1,v_0)\int_{v_1\circ n(r_1)\ge 0} 
v_1\circ n(r_1)\cdot\chi_{\{1\}}(m)p_0(r_e,v_1)\, dv_1\right)\nonumber 
 \\ 
&&\hspace{.5cm}+\left(\omega\left[\omega\chi_{\{2\}}(m) p_0(r_e,R_{r_{2} 
}(R_{r_{1}}(v_0)))\vphantom{\int_{v_1\circ n(r_1)\ge 0}}\right.\right.
\nonumber \\ 
&&\hspace{.5cm}\left.\left.+(1-\omega)M(r_2,R_{r_{1}}(v_0))\int_{v_2 
\circ n(r_2)\ge 0}v_2\circ n(r_2)\cdot\chi_{\{2\}}(m)p_0(r_e,v_2)\, dv_2 
\right]\right.\nonumber \\
&&\hspace{.5cm}+(1-\omega)M(r_1,v_0)\int_{v_1\circ n(r_1)\ge 0}v_1\circ 
n(r_1)\left[\omega\chi_{\{2\}}(m)p_0(r_e,R_{r_{2}}(v_1))\vphantom{\int_{ 
v_1\circ n(r_1)\ge 0}}\right.\nonumber \\ 
&&\hspace{.5cm}\left.\left.+(1-\omega)M(r_2,v_1)\int_{v_2\circ n(r_2) 
\ge 0}v_2\circ n(r_2)\cdot\chi_{\{2\}}(m)p_0(r_e,v_2)\, dv_2\right]\, d 
v_1\right)\nonumber \\
&&\hspace{.5cm}+\chi_{\{3\}}(m)\mbox{\rm -term}+\ldots\, ,\quad (r,v,t) 
\in\Omega\times V\times [0,\infty),\vphantom{\int} 
\end{eqnarray} 
where we recall that, for fixed $(r,v,t)\in\Omega\times V\times [0,\infty 
)$, the number $m$ is a function of $v_1,v_2,\ldots, v_e$, cf. Definition 
\ref{Definition3.2} and Remark \ref{5}. 
\er
\begin{lemma}\label{Lemma3.3}
The semigroup $S(t)$, $t\ge 0$, is strongly continuous in $L^1(\Omega\times 
V)$. 
\end{lemma}
Proof. {\it Step 1 } Let $p_0$ be a probability density everywhere defined 
on $\Omega\times V$ and let $(r,v,t)\in\Omega\times V\times [0,\infty)$. Let 
us use the notation of Remarks \ref{5} and  \ref{6} and let us continue 
from (\ref{3.20}). For $\ve>0$ set $\Omega_\varepsilon:=\{x\in\Omega:|y-x| 
>\varepsilon$ for all $y\in\partial\Omega\}$. Now, choose $\ve>0$ such that 
$\Omega_{2\ve}\neq\emptyset$. 

In fact, if $t<\varepsilon/v_{max}$ then for all $r\in\Omega_\varepsilon$ 
we have $T_{\Omega}(r_0,v_0)=T_{\Omega}(r,v)>t$. Consequently $m=0$ and 
$(r,v)=(r_e+tv_e,v_e)=(r_e+tv,v)$. Furthermore, (\ref{3.20}) implies for 
sufficiently small $\varepsilon>0$
\begin{eqnarray}\label{3.23}
S(t)p_0(r,v)=p_0(r_e,v_e)=p_0(r-tv,v)\, ,\quad r\in\Omega_\varepsilon, 
\ v\in V,\ t<\varepsilon/v_{max}. 
\end{eqnarray} 
We also observe that $\t p_0=p_0$ a.e. on $\Omega\times V$ implies $S(t) 
\t p_0=S(t)p_0$ for all $t\ge 0$. 
\medskip 

\noindent
{\it Step 2 } Let $p_0\in L^\infty(\Omega\times V)$ and let $\vp\in C_b 
(\Omega\times V)$ be uniformly continuous and non-negative on $\Omega 
\times V$. Introduce $a(t):=\sup_{v\in V}\sup_{r\in\Omega_\varepsilon}| 
\vp(r-tv,v)-\vp(r,v)|$. Choosing $t=\ve/(2v_{max})$ and recalling 
(\ref{3.23}) we verify 
\begin{eqnarray}\label{3.24}
&&\hspace{-.5cm}\left|\int_\Omega\int_V\left(S(t)p_0-p_0\right)\vp\, dv 
\, dr\right|\nonumber \\ 
&&\hspace{.5cm}\le\left|\int_{\Omega_\varepsilon}\int_V\left(p_0(r-tv,v) 
-p_0(r,v)\vphantom{l^1}\right)\vp(r,v)\, dv\, dr\right|+\int_{\Omega 
\setminus\Omega_\varepsilon}\int_V\left(S(t)p_0+p_0\right)\vp\, dv\, dr 
\nonumber \\ 
&&\hspace{.5cm}\le\left|\int_{\Omega_\varepsilon}\int_V\left(p_0(r-tv,v) 
\vp(r-tv,v)-p_0(r,v)\vp(r,v)\vphantom{l^1}\right)\, dv\, dr\right|
\nonumber \\ 
&&\hspace{1.0cm}+a(t)\, \int_{\Omega_\varepsilon}\int_Vp_0(r-tv,v)\, dv 
\, dr+\int_{\Omega\setminus\Omega_\varepsilon}\int_V\left(S(t)p_0+p_0 
\right)\vp\, dv\, dr\nonumber \\ 
&&\hspace{.5cm}\le\|\vp\|\int_{\Omega\setminus\Omega_{2\varepsilon}}\int_V 
p_0\, dv\, dr+a(t)+\|\vp\|\int_{\Omega\setminus\Omega_\varepsilon}\int_V 
S(t)p_0\, dv\, dr+\|\vp\|\int_{\Omega\setminus\Omega_\varepsilon}\int_V 
p_0\, dv\, dr\nonumber \\ 
&&\hspace{.5cm}\le 2\|\vp\|\int_{\Omega\setminus\Omega_{2\varepsilon}} 
\int_Vp_0\, dv\, dr+a(t)+\|\vp\|\left(1-\int_{\Omega_\varepsilon}\int_Vp_0 
(r-tv,v)\, dv\, dr\right)\nonumber \\ 
&&\hspace{.5cm}\le 3\|\vp\|\int_{\Omega\setminus\Omega_{2\varepsilon}} 
\int_Vp_0\, dv\, dr+a(t)\stack{t\to 0}{\lra}0\, .
\end{eqnarray} 
According to Lemma \ref{Lemma3.1} and $p_0\in L^\infty(\Omega\times V)$ there 
is some $b\in (0,\infty)$ with $\|S(t)p_0-p_0\|_{L^\infty(\Omega\times V)}<b$ 
for all $t\ge 0$. By suitable approximation of test functions in  $L^1(\Omega 
\times V)$ we verify $\int_\Omega\int_V\left(S(t)p_0-p_0\right)\vp\, dv\, dr 
\stack{t\to 0}{\lra}0$ for all $\vp\in L^1(\Omega\times V)$. In particular,   
\begin{eqnarray}\label{3.25}
\int_\Omega\int_V\left(S(t)p_0-p_0\right)\vp\, dv\, dr\stack{t\to 0}{\lra}0 
\, ,\quad\vp\in L^\infty(\Omega\times V).
\end{eqnarray} 
{\it Step 3 } Now let $p_0\in L^1(\Omega\times V)$ and $\vp\in L^\infty 
(\Omega\times V)$. Let us recall from the beginning of Section \ref{sec:2} 
that for every $t\ge 0$, $S(t)$ maps $L^1(\Omega\times V)$ into $L^1(\Omega 
\times V)$ with operator norm one. By (\ref{3.25}) we have for $\vp\in 
L^\infty(\Omega\times V)$, $N\in \mathbb{N}$, and $p_{0,N}:=p_0-(p_0\wedge N)$
\begin{eqnarray*}
&&\hspace{-.5cm}\left|\limsup_{t\to 0}\int_\Omega\int_V\left(S(t)p_0-p_0\right) 
\vp\, dv\, dr\right| \\ 
&&\hspace{.5cm}=\left|\lim_{t\to 0}\int_\Omega\int_V\left(S(t)(p_0\wedge N)-(p_0 
\wedge N)\right)\vp\, dv\, dr\right.\\
&&\hspace{1.0cm}+\left.\limsup_{t\to 0}\int_\Omega\int_V\left(S(t)p_{0,N} 
-p_{0,N}\right)\vp\, dv\, dr\right| \\
&&\hspace{.5cm}=\left|\limsup_{t\to 0}\int_\Omega\int_V\left(S(t)p_{0,N}-p_{ 
0,N}\right)\vp\, dv\, dr\right|\le 2\|p_{0,N}\|_{L^1(\Omega\times V)}\, \|\vp 
\|_{L^\infty(\Omega\times V)}\, .
\end{eqnarray*} 
Since the right-hand side can be made arbitrarily small we have (\ref{3.25}) 
for all $p_0\in L^1(\Omega\times V)$. In other words, we have $S(t)p_0\stack 
{t\to 0}{\lra}p_0$ weakly in $L^1(\Omega\times V)$. It follows now from 
\cite{Pa83}, Theorem 1.4 of Chapter 2, or \cite{EN06}, Theorem 1.6 of Chapter 
1, that $S(t)p_0\stack{t\to 0}{\lra}p_0$ strongly in $L^1(\Omega\times V)$. 
\qed 
\medskip 

\section{Spectral Properties of the Knudsen Type Semigroup and Group} 
\label{sec:4}
\setcounter{equation}{0}

Let us suppose that the global conditions (i)-(iii) and (vii) are satisfied. 
For $f\in L^1(\partial\Omega\times V)$ recall the definition 
\begin{eqnarray*}
J(r,\cdot)(f)=\int_{w\circ n(r)\ge 0}w\circ n(r)\, f(r,w)\, dw\, ,
\quad r\in\partial^{(1)}\Omega,  
\end{eqnarray*} 
where, as in Section \ref{sec:3}, the notation $J(r,\cdot)(f)$ indicates 
that $f$ does not depend on $t$. Recall also the boundary conditions 
\begin{eqnarray}\label{4.1}
f(r,v)=\omega\, f(r,R_r(v))+(1-\omega)\, M(r,v)J(r,\cdot)(f)
\end{eqnarray} 
for $(r,v)\in\partial^{(1)}\Omega\times V$ with $v\circ n(r)\le 0$. Consider 
the space $L_b^1\equiv L_b^1(\partial\Omega\times V)$ of all equivalence 
classes $f$ of measurable functions on $\partial\Omega\times V$ satisfying 
a.e. the boundary conditions (\ref{4.1}) such that 
\begin{eqnarray*}
\|f\|_{L_b^1}:=\int_{r\in\partial^{(1)}\Omega}\int_V|w\circ n(r)||f(r,w)| 
\, dw\, dr<\infty\, .
\end{eqnarray*}
Furthermore, let $L_u^1$ be the space of all equivalence classes of 
measurable functions defined on $\{(r,w)\in\partial^{(1)}\Omega\times V: 
w\circ n(r)\ge 0\}$ such that 
\begin{eqnarray*}
\|f\|_{L_u^1}:=\int_{r\in\partial^{(1)}\Omega}\int_{w\circ n(r)\ge 0}w 
\circ n(r)|f(r,w)|\, dw\, dr<\infty\, .
\end{eqnarray*}
Obviously the restriction of $f\in L_b^1$ to $\{(r,w)\in\partial^{(1)} 
\Omega\times V:w\circ n(r)\ge 0\}$ belongs to $L_u^1$. In addition every 
$f\in L_u^1$ can uniquely be extended to an element of $L_b^1$ by the 
boundary conditions (\ref{4.1}). Let us introduce a map $U$ defined on 
$L_u^1$ by 
\begin{eqnarray}\label{4.2}
Uf(r,v):=\omega\, f(r^-,R_{r^-}(v))+(1-\omega)\, M(r^-,v)J(r^-,\cdot)(f) 
\end{eqnarray} 
where $(r,v)\in\partial\Omega\times V$ such that $(r^-,v)\equiv(r^-(r,v) 
,v)\in\partial^{(1)}\Omega\times V$ and $v\circ n(r^-)\le 0$. Let $\mu\in 
\mathbb{C}$. For $f\in L_u^1$, $g\in L^1(\Omega\times V)$ and $(r,w)\in 
\partial^{(1)}\Omega\times V$ with $(r^-,w)\in\partial^{(1)}\Omega\times 
V$ and $w\circ n(r^-)\le 0$ introduce 
\begin{eqnarray*}
a(\mu)f(r,w):=e^{-T_\Omega(r,w)\mu}Uf(r,w)  
\end{eqnarray*} 
as well as 
\begin{eqnarray*}
b(\mu)g(r,w):=\int_0^{T_\Omega(r,w)}e^{-\beta\mu}g(r-\beta w,w)\, d\beta\, . 
\end{eqnarray*} 

\subsection{Spectral Properties on the boundary $\partial\Omega$} 
\label{sec:4:1}

Let $(A,D(A))$ denote the infinitesimal generator of the strongly continuous 
semigroup $S(t)$, $t\ge 0$, in $L^1(\Omega\times V)$, cf. also Lemma 
\ref{Lemma3.3}. Furthermore, let diam$(\Omega)=\sup\{|r_1-r_2|:r_1,r_2\in 
\Omega\}$ denote the diameter of $\Omega$. 
\begin{lemma}\label{Lemma4.1} 
(a) The operator $U$ is a linear operator $L_u^1\mapsto L_u^1$ with operator 
norm one. If for $f\in L_u^1$ we have $f\ge 0$ or $f\le 0$ on $\{(r,w)\in 
\partial^{(1)}\Omega\times V:w\circ n(r)\ge 0\}$, then $\|Uf\|_{L_u^1}=\|f 
\|_{L_u^1}$. \\ 
(b) Let $\mu\in\mathbb{C}$ and $g\in L^1(\Omega\times V)$. There is a unique 
$f\in D(A)$ such that $\mu f-Af=g$ if and only if there is a unique $\t f\in 
L_u^1$ such that
\begin{eqnarray}\label{4.3}
\t f=a(\mu)\t f+b(\mu)g\, .   
\end{eqnarray} 
\end{lemma} 
Proof. {\it Step 1 } Let us verify part (a). Observe that by the boundary 
conditions (\ref{4.1}) for $f\in L_b^1$ and $(r,w)\in\partial^{(1)}\Omega 
\times V$ with $w\circ n(r)\ge 0$ we have 
\begin{eqnarray}\label{4.4}
Uf(r,w)=f(r^-(r,w),w)\, ,\quad (r,w)\in\partial^{(1)}\Omega\times V\quad 
\mbox{\rm with}\quad w\circ n(r)\ge 0\, . 
\end{eqnarray} 
Furthermore, recall that the restriction of $f\in L_b^1$ to $\{(r,w)\in 
\partial^{(1)}\Omega\times V:w\circ n(r)\ge 0\}$ belongs to $L_u^1$ and that 
$f\in L_u^1$ can uniquely be extended to an element of $L_b^1$. For $f\in 
L_b^1$ it holds that   
\begin{eqnarray*}
&&\hspace{-.5cm}\int_{r\in\partial^{(1)}\Omega}\int_{w\circ n(r)\ge 0}w\circ n(r) 
\left|f\left(r^-(r,w),w\right)\right|\, dw\, dr \\  
&&\hspace{.5cm}=\int_{r^-\in\partial^{(1)}\Omega}\int_{w\circ n(r^-)\le 0}|w\circ 
n(r^-)||f(r^-,w)|\cdot\left(\frac{w\circ n(r)}{|w\circ n(r^-)|}\cdot\frac 
{dr}{dr^-}\right)\, dw\, dr^- \\ 
&&\hspace{.5cm}=\int_{y\in\partial^{(1)}\Omega}\int_{v\circ n(y)\ge 0}v\circ n(y) 
|f(y,-v)|\, dv\, dy \\ 
&&\hspace{.5cm}\le\int_{y\in\partial^{(1)}\Omega}\int_{v\circ n(y)\ge 0}v\circ n(y) 
|f(y,v)|\, dv\, dy  
\end{eqnarray*}
where the last $\le$ sign holds because of (\ref{4.1}) and (ii). We 
mention also that the equality sign holds whenever $f$ is non-negative 
or non-positive. 
\medskip

\noindent 
{\it Step 2 } We demonstrate the first part of (b). Let $g\in L^1(\Omega 
\times V)$ and $\mu\in\mathbb{C}$. Let us assume that there is a unique 
$f\in D(A)$ with $\mu f-Af=g$. Recalling condition (iii) this equation 
takes the form 
\begin{eqnarray}\label{4.5}
\mu f(y,w)+w\circ\nabla_yf(y,w)=g(y,w)\, ,\quad (y,w)\in \Omega\times V.  
\end{eqnarray} 
Here, $w\circ\nabla_yf$ is a directional derivative with respect to the 
norm in $L^1(\Omega\times V)$, that is $\|\frac1\tau(f(y+\tau w,w)-f(y,w)) 
-w\circ\nabla_yf(y,w)\|_{L^1(\Omega\times V)}\stack{\tau\to 0}{\lra}0$. It 
follows that 
\begin{eqnarray}\label{4.6} 
f(r-\gamma w,w)=e^{\gamma\mu}\left(f(r,w)-\int_0^\gamma e^{-\beta\mu}g(r- 
\beta w,w)\, d\beta\right)\, ,\quad\gamma\in [0,T_\Omega(r,w)],  
\end{eqnarray} 
a.e. on $\{(r,w)\in\partial^{(1)}\Omega\times V:w\circ n(r)\ge 0\}$. 
Conversely, applying the differentiation $w\circ\nabla_y$ in $L^1(\Omega 
\times V)$ to (\ref{4.6}) we obtain (\ref{4.5}). In addition, as a 
consequence of (\ref{4.6}), for a.e. $\{(r,w)\in\partial^{(1)}\Omega 
\times V:w\circ n(r)\ge 0\}$ and $r^-\equiv r^-(r,w)=r-T_\Omega(r,w)w$ 
\begin{eqnarray*} 
f(r,w)=e^{-T_\Omega(r,w)\mu}f(r^-,w)+\int_0^{T_\Omega(r,w)}e^{-\beta\mu} 
g(r-\beta w,w)\, d\beta\, .    
\end{eqnarray*} 
As a consequence of (iii) and Lemma \ref{Lemma3.3}, the restriction 
$\bar{f}$ of $f\in D(A)$ to $\partial\Omega\times V$ belongs to $L_b^1$. 
Recalling (\ref{4.4}), the last equation coincides therefore a.e. on 
$\{(r,w)\in\partial^{(1)}\Omega\times V:w\circ n(r)\ge 0\}$ with 
\begin{eqnarray*} 
f(r,w)=e^{-T_\Omega(r,w)\mu}Uf(r,w)+\int_0^{T_\Omega(r,w)}e^{-\beta\mu}g 
(r-\beta w,w)\, d\beta\, .    
\end{eqnarray*} 
Since the restriction $\t f$ of $f\in D(A)$ to $\{(r,w)\in\partial^{(1)} 
\Omega\times V:w\circ n(r)\ge 0\}$ is also the restriction of $\bar{f}\in 
L_b^1$ to $\{(r,w)\in\partial^{(1)}\Omega\times V:w\circ n(r)\ge 0\}$ we 
have $\t f\in L_u^1$. Thus the last equality coincides with (\ref{4.3}). 
In addition, uniqueness of the solution $f$ to $\mu f-Af=g$ implies 
uniqueness of $\t f$ in (\ref{4.3}) as follows. Replace in (\ref{4.6}) 
$f(r,w)$ with $\t f(r,w)$ and use the already mentioned equivalence of 
(\ref{4.6}) with (\ref{4.5}) and hence also with $\mu f-Af=g$. 
\medskip 

\nid 
{\it Step 3 } We show the second part of (b). Let us assume that there is a 
unique $\t f\in L_u^1$ such that we have (\ref{4.3}). As mentioned above 
$\t f\in L_u^1$ can uniquely be extended to some $\bar{f}\in L_b^1$. Keeping 
in mind relation (\ref{4.4}) from (\ref{4.3}) it follows that 
\begin{eqnarray*} 
\t f(r,w)=\bar{f}(r,w)=e^{-T_\Omega(r,w)\mu}\bar{f}(r^-,w)+\int_0^{T_\Omega 
(r,w)}e^{-\beta\mu}g(r-\beta w,w)\, d\beta     
\end{eqnarray*} 
for a.e. $(r,w)\in\partial^{(1)}\Omega\times V$ with $w\circ n(r)\ge 0$. From 
here we define a function $\hat{f}$ by 
\begin{eqnarray}\label{4.7}
\hat{f}(r-\gamma w,w):=e^{\gamma\mu}\left(\t f(r,w)-\int_0^\gamma e^{ 
-\beta\mu}g(r-\beta w,w)\, d\beta\right)\, ,\quad\gamma\in [0,T_\Omega 
(r,w)],  
\end{eqnarray} 
for a.e. $(r,w)\in\partial^{(1)}\Omega\times V$ with $w\circ n(r)\ge 0$. We 
observe
\begin{eqnarray}\label{4.8}
&&\hspace{-.5cm}\left\|\hat{f}\right\|_{L^1(\Omega\times V)}=\int_V\left 
\|\hat{f}(\cdot ,w)\right\|_{L^1(\Omega)}\, dw\nonumber \\ 
&&\hspace{.5cm}=\int_V\int_{\{r\in\partial^{(1)}\Omega:\, w\circ n(r)\ge 0\}} 
\frac{w}{|w|}\circ n(r)\int_0^{T_\Omega(r,w)\cdot |w|}\left|\hat{f}(r- 
\beta w/|w|,w)\right|\, d\beta\, dr\, dw\nonumber \\ 
&&\hspace{.5cm}=\int_{r\in\partial^{(1)}\Omega}\int_{w\circ n(r)\ge 0}w\circ 
n(r)\int_0^{T_\Omega(r,w)\cdot |w|}\frac1{|w|}\left|\hat{f}(r-\beta w/|w| 
,w)\right|\, d\beta\, dw\, dr\nonumber \\ 
&&\hspace{.5cm}=\int_{r\in\partial^{(1)}\Omega}\int_{w\circ n(r)\ge 0}w\circ 
n(r)\int_0^{T_\Omega(r,w)}\left|\hat{f}(r-\beta w,w)\right|\, d\beta\, dw 
\, dr\, , 
\end{eqnarray} 
no matter whether $\|\hat{f}\|_{L^1(\Omega\times V)}$ is finite or not. 
Now we substitute $\hat{f}(r-\beta w,w)$ by (\ref{4.7}) and estimate
\begin{eqnarray*}
&&\hspace{-.5cm}\left\|\hat{f}\right\|_{L^1(\Omega\times V)}\le\int_{r 
\in\partial^{(1)}\Omega}\int_{w\circ n(r)\ge 0}w\circ n(r)\int_0^{ 
T_\Omega(r,w)}\left|e^{\beta\mu}\right|\left|\t{f}(r,w)\right|\, d\beta 
\, dw\, dr \\ 
&&\hspace{1.0cm}+\int_{r\in\partial^{(1)}\Omega}\int_{w\circ n(r)\ge 0} 
w\circ n(r)\int_0^{T_\Omega(r,w)}\int_0^\beta e^{(\beta-\alpha)\mu} 
|g(r-\alpha w,w)|\, d\alpha\, d\beta\, dw\, dr \\ 
&&\hspace{.5cm}\le\int_0^{{\rm diam}(\Omega)/v_{min}}\left|e^{\beta\mu} 
\right|\, d\beta\cdot\|\t f\|_{L_u^1}+e^{{\rm diam}(\Omega)\, |\mathfrak 
{Re}(\mu)|/v_{min}}\cdot({\rm diam}(\Omega)/v_{min}\times \\ 
&&\hspace{1.0cm}\times\int_{r\in\partial^{(1)}\Omega}\int_{w\circ n(r) 
\ge 0}w\circ n(r)\int_0^{T_\Omega(r,w)} |g(r-\alpha w,w)|\, d\alpha\, 
dw\, dr  \\ 
&&\hspace{.5cm}\le\int_0^{{\rm diam}(\Omega)/v_{min}}\left|e^{\beta\mu} 
\right|\, d\beta\cdot\|\t f\|_{L_u^1}  \\ 
&&\hspace{1.0cm}+e^{{\rm diam}(\Omega)\, |\mathfrak{Re}(\mu)|/v_{min}} 
\cdot({\rm diam}(\Omega)/v_{min}\cdot\|g\|_{L^1(\Omega\times V)}<\infty 
\vphantom{\int}
\end{eqnarray*} 
where for the last line we have applied (\ref{4.8}) once again. This says 
$\hat{f}\in L^1(\Omega\times V)$. 

Our task is now to show that $\hat{f}$ belongs to $D(A)$ and that we have 
$\mu\hat{f}-A\hat{f}=g$. Choosing $\gamma=0$ in (\ref{4.7}) we deduce 
that $\t f$ coincides a.e.~with the restriction of $\hat f$ to $\{(r,w)\in 
\partial\Omega\times V:w\circ n(r)\ge 0\}$. In other words, we have 
(\ref{4.6}) with $\hat{f}$ instead of $f$. Its equivalence to (\ref{4.5}) 
has already been mentioned in Step 2 of this proof, i.e. we have $\mu\hat 
{f}-A\hat{f}=g$. Noting that (\ref{4.7}) is (\ref{4.6}) in the setup of 
the present Step 3, the  equivalence of (\ref{4.5}) and (\ref{4.6}) shows 
also that uniqueness of $\t f\in L_u^1$ in (\ref{4.3}) implies uniqueness 
of $f$ in $\mu f-Af=g$.
\qed 
\medskip

In order to solve (\ref{4.3}) for $f\in L_u^1$ let us examine the operator 
$id-a(\mu)$. Since $(A,D(A))$ is the generator of a strongly continuous 
semigroup in $L^1(\Omega\times V)$, cf. Lemma \ref{Lemma3.3}, by the 
Hille-Yosida theorem and our Lemma \ref{Lemma4.1} (b) it is sufficient 
to focus on $\mu\in\mathbb C$ with $\mathfrak{Re}\mu\le 0$. For $(r,w)\in 
\partial^{(1)}\Omega\times V$ with $w\circ n(r)\ge 0$ introduce 
\begin{eqnarray}\label{4.9} 
A(\mu)f(r,w):=\omega e^{-T_\Omega(r,w)\mu}f\left(r^-,R_{r^-}(w)\right)
\end{eqnarray} 
and 
\begin{eqnarray*}
B(\mu)f(r,w):=(1-\omega)e^{-T_\Omega(r,w)\mu} M(r^-,w)J(r^-,\cdot)(f) 
\end{eqnarray*} 
which, because of Lemma \ref{Lemma4.1} (a) and (ii), are well-defined 
on $f\in L_u^1$. Relation (\ref{4.2}) yields the decomposition 
\begin{eqnarray}\label{4.10}
a(\mu)f(r,w)=A(\mu)f(r,w)+B(\mu)f(r,w)\, . 
\end{eqnarray} 
By Lemma \ref{Lemma4.1} (a), $a(\mu)$, $A(\mu)$, and $B(\mu)$ are 
bounded linear operators on $L_u^1$. Next we are concerned with the 
bijectivity of $id-A(\mu)$ on $L_u^1$. 
\medskip

For this we need some preparations. For $(r,w)\in\partial^{(1)}\Omega 
\times V$ with $w\circ n(r)\ge 0$ let us introduce $\mathbb{T}^0(r,w) 
:=(r,w)$, 
\begin{eqnarray*}
\mathbb{T}(r,w):=\left(r^-(r,w),R_{r^-}(w)\right)\, , 
\end{eqnarray*} 
and the abbreviations $l_k\equiv l_k(r,w):=T_\Omega\left(\mathbb{T}^k 
(r,w)\right)$ as well as $\left(r_k,w_k\right)\equiv\left(r_k(r,w),w_k 
(r,w)\right)$ $:=\mathbb{T}^k(r,w)$, $k\in\mathbb{Z}_+$. Noting that 
\begin{eqnarray*}
\mathbb{T}^{-1}(r,w)=\left(r^-(r,-R_r(w)),R_r(w)\right) 
\end{eqnarray*} 
for all $(r,w)\in\partial^{(1)}\Omega\times V$ with $w\circ n(r)\ge 0$ 
we also abbreviate $l_k\equiv l_k(r,w):=T_\Omega\left(\mathbb{T}^k(r,w) 
\right)$ and $\left(r_k,w_k\right)\equiv\left(r_k(r,w),w_k(r,w)\right) 
:=\mathbb{T}^k(r,w)$, $-k\in\mathbb{N}$. We emphasize that the notation 
just introduced is adjusted to the remainder of this section. It is not 
compatible to Definition \ref{Definition3.2} and Remarks \ref{5} as 
well as \ref{6}, but will not cause any ambiguity. However note the  
compatibility to condition (vii) by the equivalence of $\sigma^{-1}
(r,-w)=(y,v)$ and $\mathbb{T}(r,w)=(y,-v)$ for $(r,w)$ as above. 
\medskip

\nid
\br{7} Let us take advantage of a well known fact from the ergodic 
theory of $d$-dimensional mathematical billiards, $d=2,3$. For a.e. 
$(r,w)\in\partial^{(1)}\Omega\times V$ with $w\circ n(r)\ge 0$ there 
is a $\tau(r,w)\le\, $diam$(\Omega)/v_{min}$ such that 
\begin{eqnarray}\label{4.11}
\lim_{n\to\infty}\frac1n\sum_{k=0}^{n-1}l_k=\lim_{n\to\infty}\frac1n 
\sum_{k=1}^nl_{-k}=\tau(r,w)\, . 
\end{eqnarray} 
Relation (\ref{4.11}) holds even a.e. on $(r,v)\in\{\partial^{(1)} 
\Omega\times S^{d-1}:v\circ n(r)\ge 0\}$, where $S^{d-1}$ denotes the 
unit sphere. 

One may recover these results by using the excellent source 
\cite{CM03} in the following way. First recall the Ergodic Theorem of 
Birkhoff-Khinchin, Theorem II.1.1 together with the Corollary II.1.4. 
Then recall from \cite{CM03}, (IV.2.3) and the formula above (IV.2.3), 
that the billiard map $T$ is measure preserving with respect to some 
probability measure $\nu$ on $\{(r,v)\in\partial^{(1)}\Omega\times 
S^{d-1}:v\circ n(r)\ge 0\}$. In fact, the Radon-Nikodym derivative of 
$\nu$ with respect to the canonical measure $\lambda$ on $\{\partial^{(1) 
}\Omega\times S^{d-1}:v\circ n(r)\ge 0\}$ is 
\begin{eqnarray}\label{4.12}
\frac{d\nu}{d\lambda}=c_\nu\cdot n(r)\circ v 
\end{eqnarray} 
where $c_\nu>0$ is a normalizing constant. Now it remains to mention 
that $\mathbb{T}$ leaves the modulus of the velocity $|w|$ invariant 
an that, assuming for a moment that $S^{d-1}\subset V$, the map 
$\mathbb{T}$ is for $|w|=1$ identical with the billiard map $T$ of 
\cite{CM03}, up to orientation. 
\er

Introduce 
\begin{eqnarray*}
\mathbf{d}:=\{(r,v)\in\partial^{(1)}\Omega \times S^{d-1}:v\circ n(r) 
\ge 0,\ \mbox{\rm and we have (\ref{4.11})}\}\, .
\end{eqnarray*} 

\noindent
\br{8} For the well-definiteness of the term (\ref{4.16}) below let us  
verify $\tau(r,v)>0$ for a.e. $(r,v)\in\mathbf{d}$, which also will imply 
$\tau(r,w)>0$ for a.e. $(r,w)\in\partial^{(1)}\Omega\times V$ with $w 
\circ n(r)\ge 0$. For this we will use 
\begin{eqnarray}\label{4.13}
\int_\mathbf{d}\lim_{n\to\infty}\frac1n\#\left\{0\le j<n:T^j(r,v)\in A 
\right\}\, \nu(d(r,v))=\nu(A)
\end{eqnarray} 
for all Borel subsets $A$ of $\mathbf{d}$, see \cite{CM03} below Theorem 
II.1.1. Keep in mind $|v|=1$ and consider $l_{-k}\equiv l_{-k}(r,v)$ as a 
function $l_{-k}(r,v)\equiv l(T^k(r,v))$, $k\in\mathbb{N}$, $(r,v)\in 
\mathbf{d}$. 

Indeed, let us assume that for $B:=\{(r,v)\in\mathbf{d}:\tau(r,v)=0\}$ we 
have $\nu(B)>0$. With $\ve>0$, $f_\ve (r,v):=\chi_{[\ve,\infty)}(l(T(r,v)) 
)$, $(r,v)\in\mathbf{d}$, and $A_\ve:=\left\{(r,v)\in\mathbf{d}:l(T(r,v)) 
\ge\ve\right\}$ the following chain of inclusions holds. 
\begin{eqnarray}\label{4.14}
&&\hspace{-.5cm}B=\{(r,v)\in\mathbf{d}:\tau(r,v)=0\}\subseteq\left\{(r,v) 
\in\mathbf{d}:\lim_{n\to\infty}\frac1n\sum_{k=1}^nf_\ve\left(T^k(r,v)
\right)=0\right\}\nonumber \\ 
&&\hspace{.5cm}=\left\{(r,v)\in\mathbf{d}:\lim_{n\to\infty}\frac1n\#\left 
\{0\le j<n:T^j(r,v)\in A_\ve\right\}=0\right\}=:B_\ve\, .
\end{eqnarray} 
Now replace $A$ in (\ref{4.13}) by $A_\ve$ and let $\ve\to 0$. Then the 
{\it right-hand side of (\ref{4.13}) tends to one} by the definition of 
$A_\ve$. However, the integral on the left-hand side can be decomposed into 
the two parts $\int_{B_\ve}$ and $\int_{\mathbf{d}\setminus B_\ve}$. We 
observe $\lim_{\ve\to 0}\int_{B_\ve}=0$ by (\ref{4.14}) and $\limsup_{\ve 
\to 0}\int_{\mathbf{d}\setminus B_\ve}\le 1-\nu(B)$ by construction and 
$\nu(B_\ve)\ge\nu(B)$. The assumption $\nu(B)>0$ implies now that limsup of 
the {\it left-hand side of (\ref{4.13}) is bounded by $1-\nu(B)<1$}. Thus 
the above assumption cannot hold, i.e. $\tau(r,v)>0$ for a.e. $(r,v)\in 
\mathbf{d}$ and also $\tau(r,w)>0$ for a.e. $(r,w)\in\partial^{(1)}\Omega 
\times V$ with $w\circ n(r)\ge 0$. 
\er

In other words, in Remark \ref{8} it is demonstrated that with 
\begin{eqnarray*}
\mathbf{D}:=\left\{(r,w)\in\partial\Omega\times V:w=|w|\cdot v,\ (r,v)\in 
\mathbf{d},\ \tau(r,w)>0\right\} 
\end{eqnarray*} 
the canonical measure of $\{(r,w)\in\partial^{(1)}\Omega\times V:w\circ 
n(r)\ge 0\}\setminus\mathbf{D}$ is zero. In particular we stress that if 
$w=|w|\cdot v\in V$, $(r,v)\in\mathbf{d}$, and $\tau(r,v)>0$ then we have 
$|w|\cdot\tau(r,w)=\tau(r,v)$. 

\br{9} Let us introduce a probability measure $\rm N$ on $\{(r,w)\in 
\partial^{(1)}\Omega\times V:w\circ n(r)\ge 0\}$ by its Radon-Nikodym 
density with respect to the canonical measure. For this, let $\lambda$ 
denote the canonical measure on $\{(r,v)\in\partial^{(1)}\Omega\times 
S^{d-1}:v\circ n(r)\ge 0\}$ and let $\Lambda$ be the canonical measure 
on $\{(r,w)\in\partial^{(1)}\Omega\times V:w\circ n(r)\ge 0\}$. Set 
\begin{eqnarray}\label{4.15}
\frac{d{\rm N}}{d\Lambda}(r,w):=c_{\rm N}\cdot|w|\, \frac{d\nu}{d\lambda} 
(r,w/|w|)=c_{\rm N}\cdot w\circ n(r)
\end{eqnarray} 
where $c_{\rm N}>0$ is a normalizing constant, cf. also (\ref{4.12}). Let 
us recall the following facts. The billiard map $T$ of $\cite{CM03}$ 
preserves the measure $\nu$ of $\cite{CM03}$, Section IV.2. Assuming for 
a moment that $S^{d-1}\subset V$, the map $\mathbb{T}$ coincides for fixed 
modulus of the velocity $|w|=1$ with $T$, up to orientation. In addition, 
$\mathbb{T}$ leaves the modulus of the velocity $|w|$ invariant. As a 
consequence we conclude that the map $\mathbb{T}$ preserves the measure 
$\rm N$. Furthermore, the definition (\ref{4.15}) shows the identity of 
$L^1(\mathbf{D},{\rm N})$ with $L_u^1$, up to zero-sets; we indicate this 
by $L^1(\mathbf{D},{\rm N})\simeq L_u^1$. 
\er

Recall the definition of $\tau(r,w)\cdot |w|$ in (\ref{4.11}) and 
introduce 
\begin{eqnarray*}
c_\tau:=\mathop{\mathrm{ess~inf}}\limits_{(r,w)\in\mathbf{D}}\tau(r,w) 
\cdot |w| 
\end{eqnarray*} 
as well as 
\begin{eqnarray}\label{4.16}
C_\tau\equiv C_{\tau,k_0}:=\mathop{\mathrm{ess~sup}}\limits_{(r,w)\in 
\mathbf{D},\, k\ge k_0}\frac{k\cdot\tau(r,w)-\sum_{i=1}^kl_i(r,w)}{k 
\cdot\tau (r,w)}\, ,\quad k_0\in\mathbb{N}.
\end{eqnarray} 
\br{10} Below we will suppose $C_\tau\equiv C_{\tau,k_0}<1$. To 
motivate this assumption let $\Omega$ be a convex polygon in dimension 
$d=2$ or a convex polyhedron in dimension $d=3$. Suppose that the inner 
angles between two neighboring edges if $d=2,3$ or faces if $d=3$ are at 
least $\pi/2$. 

Consider points $(r,v)\in\partial^{(1)}\Omega\times S^{d-1}$ with $v\circ 
n(r)>0$ and $(r_k,v_k)=\mathbb{T}_k(r,v)$ such that $r_k$ is neither a 
vertex nor belongs to an edge if $d=3$, for any $k\in \mathbb{Z}$. For 
$d=2$ denote by $s_{min}$ the length of the shortest edge of $\partial 
\Omega$. For $d=3$ denote by $s_{min}$ the infimum over all distances 
between two points $y_1$ and $y_2$ belonging to two faces of $\partial 
\Omega$ having no edge or vertex in common. 

As above, let diam$(\Omega)$ denote the diameter of $\Omega$. Then for 
any $j\in\mathbb{Z}$ it holds that 
\begin{eqnarray*} 
s_{min}\le\left\{
\begin{array}{lc}
|r_j-r_{j+1}|+|r_{j+1}-r_{j+2}|&\quad\mbox{\rm if}\ d=2 \\
|r_j-r_{j+1}|+|r_{j+1}-r_{j+2}|+|r_{j+2}-r_{j+3}|&\quad\mbox{\rm if}\ 
d=3
\end{array}\right.\, . 
\end{eqnarray*} 
This is obvious for $d=2$. For $d=3$ it follows from the fact that 
there cannot be more than three consecutive reflections on faces 
having a vertex in common. Now we may choose $k_0=2$ if $d=2$ and 
$k_0=3$ if $d=3$ to obtain 
\begin{eqnarray*}
C_\tau\le\frac{(2d-1)\cdot\mbox{\rm diam}(\Omega)-s_{min}}{(2d-1)\cdot 
\mbox{\rm diam}(\Omega)}<1\, .
\end{eqnarray*} 
It follows also that 
\begin{eqnarray*}
c_\tau\ge\frac13\, s_{min}\, .
\end{eqnarray*} 
\er

Introduce 
\begin{eqnarray*}
m:=\frac{\log\omega\cdot v_{max}}{c_\tau}\quad \mbox{\rm if }c_\tau>0\ 
\mbox{\rm and }m:=-\infty\ \mbox{\rm if }c_\tau=0\, . 
\end{eqnarray*} 
Throughout this section, let us use the convention $\frac{\mu}{0} 
\sum_{i=0}^{-1}=0$. 

\begin{lemma}\label{Lemma4.2} 
Let $\omega\in (0,1)$. (a) There exists an at most countable subset 
$\mathfrak{M}$ of $[m,0]$ such that we have the following. For 
\begin{eqnarray*} 
\mu\in\mathbb{M}:=\{\lambda\in\mathbb{C}:\mathfrak{Re}\lambda\in 
(-\infty,0]\setminus\mathfrak{M}\}
\end{eqnarray*} 
and any given a.e. bounded measurable $\psi$ there exists an a.e. unique 
measurable $\vp$, both defined on $\{(r,w)\in\partial^{(1)}\Omega\times 
V:w\circ n(r)\ge 0\}$, such that 
\begin{eqnarray}\label{4.17}
\vp(r,w)-A(\mu)\vp(r,w)\equiv\vp(r,w)-\omega e^{-T_\Omega(r,w)\mu}\vp 
\left(r^-,R_{r^-}(w)\right)=\psi(r,w) 
\end{eqnarray} 
for a.e. $(r,w)\in\partial^{(1)}\Omega\times V$ with $w\circ n(r)\ge 0$. \\ 
(b) Suppose $(1-C_{\tau,k_0})\, c_\tau>0$ for some $k_0\in\mathbb{N}$. 
Then for 
\begin{eqnarray*}
\mu\in\mathbb{M}_m:=\left\{\mu\in\mathbb{C}:\mathfrak{Re}\mu<\frac{m}{1 
-C_\tau}\right\} 
\end{eqnarray*} 
the operator $id-A(\mu)$ is bijective in $L_u^1$. 
\end{lemma}
Proof. {\it Step 1 } We point out two equivalent representations of 
(\ref{4.17}), namely (\ref{4.18}) and (\ref{4.20}) below. Relation 
(\ref{4.17}) is equivalent to 
\begin{eqnarray}\label{4.18}
&&\hspace{-.5cm}\vp(r,w)-\exp\left\{-n\cdot\left(\frac{\mu}{n}\sum_{k=0 
}^{n-1}l_k-\log\omega\right)\right\}\vp(r_n,w_n)\nonumber \\ 
&&\hspace{.5cm}=\sum_{k=0}^{n-1}\exp\left\{-k\cdot\left(\frac{\mu}{k} 
\sum_{i=0}^{k-1}l_i-\log\omega\right)\right\}\psi(r_k,w_k)\, ,\quad n\in 
\mathbb{N}. 
\end{eqnarray} 
Recalling the definitions above Remark \ref{7} and replacing $(r,w)$ with 
$\left(r^-(r,-R_r(w)),R_r(w)\right)$, and hence $\left(r^-,R_{r^-}(w) 
\right)$ with $(r,w)$, equation $(\ref{4.17})$ turns into 
\begin{eqnarray}\label{4.19}
&&\hspace{-.5cm}\vp(r,w)-\frac{e^{T_\Omega\left(r^-(r,-R_r(w)),R_r(w) 
\right)\mu}}{\omega}\vp\left(r^-(r,-R_r(w)),R_r(w)\right)\nonumber \\ 
&&\hspace{.5cm}=-\frac{e^{T_\Omega\left(r^-(r,-R_r(w)),R_r(w)\right)\mu}} 
{\omega}\psi\left(r^-(r,-R_r(w)),R_r(w)\right)\, . 
\end{eqnarray} 
Similarly to (\ref{4.18}) we obtain from (\ref{4.19})
\begin{eqnarray}\label{4.20}
&&\hspace{-.5cm}\vp(r,w)-\exp\left\{n\cdot\left(\frac{\mu}{n}\sum_{k=1}^n 
l_{-k}-\log\omega\right)\right\}\vp(r_{-n},w_{-n})\nonumber \\ 
&&\hspace{.5cm}=-\sum_{k=1}^n\exp\left\{k\cdot\left(\frac{\mu}{k}\sum_{i= 
1}^kl_{-i}-\log\omega\right)\right\}\psi(r_{-k},w_{-k})\, ,\quad n\in 
\mathbb{N}. 
\end{eqnarray} 
{\it Step 2 } In this step, we construct the set $\mathbb{M}$. Recall the 
hypothesis $\omega\in (0,1)$. Let 
\begin{eqnarray*}
\mathbf{D}^<_\mu:=\{(r,w)\in\mathbf{D}:\mathfrak{Re}\mu\cdot\tau(r,w)< 
\log\omega\}\, , 
\end{eqnarray*} 
and 
\begin{eqnarray*}
\mathbf{D}^>_\mu:=\{(r,w)\in\mathbf{D}:\mathfrak{Re}\mu\cdot\tau(r,w)>
\log\omega\}\, ,\quad\mathfrak{Re}\mu\le 0.
\end{eqnarray*} 
Furthermore introduce 
\begin{eqnarray*}
\mathbf{D}\setminus\left(\mathbf{D}^<_\mu\cup\mathbf{D}^>_\mu\right)=\left 
\{(r,w)\in\mathbf{D}:\mathfrak{Re}\mu=\frac{\log\omega}{\tau(r,w)}\right\} 
=:\mathbf{D}^=_\mu\, ,\quad\mathfrak{Re}\mu\le 0.
\end{eqnarray*} 
We observe that for $\omega\in (0,1)$ and $a\in (-\infty,0)$, the set 
$\mathbf{D}^=_{a+bi}$ is independent of $b\in\mathbb{R}$. Since, for 
fixed $\omega\in (0,1)$, there are at most countably many $a\in 
(-\infty,0)$ such that $\mathbf{D}^=_{a+bi}$ has positive canonical 
measure the set 
\begin{eqnarray}\label{4.21}
\mathbb{M}:=\left\{\mu\in\mathbb{C}:\mathfrak{Re}\mu\le 0,\ \mathbf 
{D}^=_\mu\ \mbox{\rm has zero canonical measure}\right\} 
\end{eqnarray} 
coincides with $\{\mu\in\mathbb{C}:\mathfrak{Re}\mu\le 0\}$ except 
for a union of at most countably many vertical lines of the form 
$\{\mu\in\mathbb{C}:\mathfrak{Re}\mu\in\mathfrak{M}\}$ where 
$\mathfrak{M}$ is an at most countable set of non-positive real 
numbers. We mention that for $\mu\in\mathbb{M}$ we have a.e. 
on $\{(r,w)\in\partial^{(1)}\Omega\times V:w\circ n(r)\ge 0\}$ the 
alternative that either $(r,w)\in\mathbf{D}^>_\mu$ or $(r,w)\in 
\mathbf{D}^<_\mu$. Furthermore if $c_\tau>0$, we have $\mathfrak 
{Re}\mu\cdot c_\tau <\log\omega\cdot v_{max}$ if and only if 
$\mathfrak{Re}\mu<m$ and therefore
\begin{eqnarray}\label{4.22}
\mathbf{D}^>_\mu=\emptyset\quad\mbox{\rm and}\quad\mathbf{D}^=_\mu= 
\emptyset\, ,\quad\mbox{\rm if }c_\tau>0\ \mbox{\rm and }\mathfrak{Re} 
\mu<m.  
\end{eqnarray} 
In other words, we have $\mathfrak{M}\subset [m,0]$ and $\mathbb{M}= 
\{\mu\in\mathbb{C}:\mathfrak{Re}\mu\in (-\infty,0]\setminus\mathfrak 
{M}\}$, no matter if $c_\tau>0$ or $c_\tau=0$. 
\medskip

\nid 
{\it Step 3 } Next we will be concerned with the existence/non-existence 
of the limits as $n\to\infty$ in (\ref{4.18}) and (\ref{4.20}). In 
this step we restrict the analysis to $(r,w)\in\mathbf{D}^<_\mu\cup 
\mathbf{D}^>_\mu$. It holds that 
\begin{eqnarray}\label{4.23}
&&\hspace{-.5cm}\vp^<_1(r,w):=-\sum_{k=1}^\infty\exp\left\{k\cdot 
\left(\frac{\mu}{k}\sum_{i=1}^{k}l_{-i}(r,w)-\log\omega\right)\right\} 
\nonumber \\ 
&&\hspace{0.5cm}=-e^{\mu\cdot l_{-1}(r,w)-\log\omega}\cdot\sum_{k=0 
}^\infty\exp\left\{k\cdot\left(\frac{\mu}{k}\sum_{i=1}^kl_{-i}(r_{-1}, 
w_{-1})-\log\omega\right)\right\}\nonumber \\ 
&&\hspace{0.5cm}=e^{\mu\cdot l_{-1}(r,w)-\log\omega}\cdot\left(-1+ 
\vp^<_1(r_{-1},w_{-1})\right)\, ,\vphantom{\left\{\sum_{k=0}^{m-1} 
\right\}}  
\end{eqnarray} 
no matter whether the sum $\vp^<_1(r,w)$ converges or not. We obtain 
immediately that the sum $\vp^<_1(r,w)$ converges if and only if 
$\vp^<_1(r_k,w_k)$ converges for every $k\in\mathbb{Z}$. 

Iterating (\ref{4.23}) we find 
\begin{eqnarray}\label{4.24}
&&\hspace{-0.5cm}\vp^<_1(r,w)=-\sum_{k=1}^n\exp\left\{k\cdot\left( 
\frac{\mu}{k}\sum_{i=1}^kl_{-i}(r,w)-\log\omega\right)\right\} 
\nonumber \\ 
&&\hspace{1.0cm}+\exp\left\{\sum_{i=1}^n(\mu\cdot l_{-i}(r,w)-\log 
\omega)\right\}\vp^<_1(r_{-n},w_{-n})\, ,\quad n\in\mathbb{N}.
\end{eqnarray} 
Recalling that $(r,w)\in\mathbf{D}^<_\mu\cup\mathbf{D}^>_\mu$, we 
observe 
\begin{eqnarray}\label{4.25}
\sum_{i=1}^\infty(\mathfrak{Re}\mu\cdot l_{-i}(r,w)-\log\omega)= 
\mbox{\rm sign }(\mathfrak{Re}\mu\cdot\tau(r,w)-\log\omega)\cdot 
\infty\, . 
\end{eqnarray} 

In the remainder of this step we show the equivalence of the 
following. 
\begin{itemize}
\item[(1)] The sum $\vp^<_1(r,w)$ given in the first line of 
(\ref{4.23}) converges. 
\item[(2)] All $\vp^<_1(r_k,w_k)$, $k\in\mathbb{Z}$, converge. 
\item[(3)] We have $\mathfrak{Re}\mu\cdot\tau(r,w)<\log\omega$, 
i.e. $(r,w)\in\mathbf{D}^<_\mu$.
\end{itemize}
The equivalence of (1) and (2) has been noted above. Let us proceed 
from (1) to (3) and back to (1). If the sum $\vp^<_1(r,w)$ converges 
then the first term on the right-hand side of (\ref{4.24}) converges 
to $\vp^<_1(r,w)$ as $n\to\infty$. Since $\vp^<_1(r_n,w_n)$ cannot 
tend to zero as $n\to\infty$ by (\ref{4.23}), from (\ref{4.24}) and 
(\ref{4.25}) we may now conclude the following.
\begin{eqnarray}\label{4.26}
\mbox{\rm If }\ \vp^<_1(r,w)\ \mbox{\rm converges then }\lim_{n\to 
\infty}\sum_{i=1}^n(\mathfrak{Re}\mu\cdot l_{-i}(r,w)-\log\omega) 
=-\infty\, .  
\end{eqnarray} 
By (\ref{4.25}) we have (3). In this case, there is a $k_1\in\mathbb 
{N}$ such that for $j>k_1$ 
\begin{eqnarray*}
\mathfrak{Re}\mu\cdot\sum_{i=1}^jl_{-i}(r,w)-j\log\omega<\frac{j}{2} 
(\mathfrak{Re}\mu\cdot\tau(r,w)-\log\omega)\, . 
\end{eqnarray*} 
Letting $k,l>k_1$ and using the first line of (\ref{4.23}) we obtain 
\begin{eqnarray}\label{4.27}
&&\hspace{-.5cm}\left|\sum_{j=1}^{k}\exp\left\{j\cdot\left(\frac{\mu} 
{j}\sum_{i=1}^jl_{-i}(r,w)-\log\omega\right)\right\}-\sum_{j=1}^{l} 
\exp\left\{j\cdot\left(\frac{\mu}{j}\sum_{i=1}^jl_{-i}(r,w)-\log\omega 
\right)\right\}\right|\nonumber \\ 
&&\hspace{0.5cm}\le\sum_{j=k\wedge l+1}^\infty\exp\left\{j\cdot\left( 
\frac{\mathfrak{Re}\mu}{j}\sum_{i=1}^jl_{-i}(r,w)-\log\omega\right) 
\right\}\nonumber \\ 
&&\hspace{0.5cm}\le\sum_{j=k\wedge l+1}^\infty\exp\left\{\frac{j}{2} 
(\mathfrak{Re}\mu\cdot\tau(r,w)-\log\omega)\right\}\stack{k,l\to\infty} 
{\lra}0\, , 
\end{eqnarray} 
i.e. the sum $\vp^<_1(r,w)$ converges provided that $\mathfrak{Re}\mu 
\cdot\tau(r,w)<\log\omega$. We may now state the equivalence of the 
properties (1)-(3). 

Introducing 
\begin{eqnarray*}
\vp^>_1(r,w):=\sum_{k=0}^\infty\exp\left\{-k\cdot\left(\frac{\mu}{k} 
\sum_{i=0}^{k-1}l_i(r,w)-\log\omega\right)\right\} 
\end{eqnarray*} 
the equivalence of the following properties (4)-(6) can be proved in 
a similar way. 
\begin{itemize}
\item[(4)] The sum $\vp^>_1(r,w)$ converges. 
\item[(5)] All $\vp^>_1(r_k,w_k)$, $k\in\mathbb{Z}$, converge. 
\item[(6)] We have $\mathfrak{Re}\mu\cdot\tau(r,w)>\log\omega$, i.e. 
$(r,w)\in\mathbf{D}^>_\mu$. 
\end{itemize}

\nid 
{\it Step 4 } Let us prove part (a). Without mentioning this again in 
the present step, we take advantage of the equivalences (1)-(3) and 
(4)-(6) of Step 3. From (\ref{4.18}) and (\ref{4.20}) we obtain 
necessarily  
\begin{eqnarray}\label{4.28}
&&\hspace{-.5cm}\vp(r,w)=\sum_{k=0}^\infty\exp\left\{-k\cdot\left(\frac 
{\mu}{k}\sum_{i=0}^{k-1}l_i-\log\omega\right)\right\}\psi(r_k,w_k)\, , 
\end{eqnarray} 
for a.e. $(r,w)\in\mathbf{D}^>_\mu$, as well as 
\begin{eqnarray}\label{4.29}
&&\hspace{-.5cm}\vp(r,w)=-\sum_{k=1}^\infty\exp\left\{k\cdot\left( 
\frac{\mu}{k}\sum_{i=1}^kl_{-i}-\log\omega\right)\right\}\psi(r_{-k}, 
w_{-k})\, , 
\end{eqnarray} 
for a.e. $(r,w)\in\mathbf{D}^<_\mu$ whenever $\mu\in\mathbb{M}$. Checking 
back with (\ref{4.17}) we verify that (\ref{4.28}), (\ref{4.29}) is 
the unique solution to (\ref{4.17}). 
\medskip

\noindent
{\it Step 5 } In order to verify part (b), let us also recall the 
definition of $(r_{-k},w_{-k})$, $k\in\mathbb{N}$, above Remark \ref{7} 
and the definition of $\tau(r,w)$ in (\ref{4.11}). Review thoroughly 
Remark \ref{9}. The maps $\mathbb{T}$ and $\mathbb {T}^{-1}$ preserve 
the measure $\rm N$ given by (\ref{4.15}). 

Note that for $\mathfrak{Re}\mu<m$ it holds that $\mathbf{D}=\mathbf{D}^<$ 
by (\ref{4.22}). By the construction of $m$ and $\mathbb{M}_m$ before 
and in the formulation of Lemma \ref{Lemma4.2}, $\, \mathfrak{Re}\mu\cdot 
(1-C_\tau)c_\tau/v_{max}-\log\omega$ is for $\mu\in\mathbb{M}_m$ negative 
on $\mathbf{D}$. 
\medskip

By means of the above preparations it follows as in (\ref{4.27}) that for 
$\psi\in L_u^1\simeq L^1(\mathbf{D},{\rm N})$ and $\mu\in\mathbb{M}_m$ we 
have 
\begin{eqnarray*}
&&\hspace{-.5cm}\left|\sum_{k=n_0}^n\exp\left\{k\cdot\left(\frac{\mu}{k} 
\sum_{i=1}^kl_{-i}-\log\omega\right)\right\}\psi(r_{-k},w_{-k})\right| \\ 
&&\hspace{.5cm}\le\sum_{k=n_0}^n\exp\left\{k\cdot\left(\frac{\mathfrak{Re} 
\mu}{k}\sum_{i=1}^kl_{-i}-\log\omega\right)\right\}|\psi(r_{-k},w_{-k})|  
 \\ 
&&\hspace{.5cm}=\sum_{k=n_0}^n\exp\left\{k\cdot\mathfrak{Re}\mu\left( 
\frac{1}{k}\sum_{i=1}^kl_{-i}-\tau(r,w)\right)\right\}\cdot e^{k\cdot\left( 
\mathfrak{Re}\mu\cdot\tau(r,w)-\log\omega\right)}|\psi(r_{-k},w_{-k})| \\ 
&&\hspace{.5cm}\le\sum_{k=n_0}^ne^{-k\cdot\mathfrak{Re}\mu\cdot C_\tau\tau 
(r,w)}\cdot e^{k\cdot\left(\mathfrak{Re}\mu\cdot\tau(r,w)-\log\omega\right)} 
|\psi(r_{-k},w_{-k})| \\ 
&&\hspace{.5cm}\le\sum_{k=n_0}^n e^{k\cdot\left(\mathfrak{Re}\mu\cdot(1- 
C_\tau)c_\tau/v_{max}-\log\omega\right)}|\psi(r_{-k},w_{-k})|\, ,\quad k_0 
\le n_0\le n,\ n\in\mathbb{N}. 
\end{eqnarray*} 
We obtain 
\begin{eqnarray*}
&&\hspace{-.5cm}\int_{\mathbf{D}}\left|\sum_{k=n_0}^n\exp\left\{k\cdot 
\left(\frac{\mu}{k}\sum_{i=1}^kl_{-i}-\log\omega\right)\right\}\psi(r_{-k} 
,w_{-k})\right|\, d{\rm N} \\ 
&&\hspace{.5cm}\le\int_{\mathbf{D}}\sum_{k=n_0}^n e^{k\cdot\left(\mathfrak 
{Re}\mu\cdot (1-C_\tau)c_\tau/v_{max}-\log\omega\right)}|\psi(r_{-k},w_{-k} 
)|\, d{\rm N}\vphantom{\left|\sum_{i=1}^k\right|} \\ 
&&\hspace{.5cm}=\sum_{k=n_0}^n e^{k\cdot\left(\mathfrak{Re}\mu\cdot (1- 
C_\tau)c_\tau/v_{max}-\log\omega\right)}\int_{\mathbf{D}}|\psi(r,w)|\, d 
{\rm N}\vphantom{\left|\sum_{i=1}^k\right|} \\ 
&&\hspace{.5cm}<\frac{e^{n_0\cdot\left(\mathfrak{Re}\mu\cdot (1-C_\tau) 
c_\tau/v_{max}-\log\omega\right)}}{1-e^{\mathfrak{Re}\mu\cdot (1-C_\tau) 
c_\tau/v_{max}-\log\omega}}\int_{\mathbf{D}}|\psi(r,w)|\, d{\rm N}\vphantom 
{\left|\sum_{i=1}^k\right|} \\ 
&&\hspace{.0cm}\stack{}{\lra}0\quad \mbox{\rm as }n_0\to\infty\, ,\vphantom 
{\left|\sum_{i=1}^k\right|}
\end{eqnarray*} 
the third line because of the fact that $\mathbb{T}$ preserves the measure 
${\rm N}$. This together with the a.e. result (\ref{4.29}) shows that the 
term on the right-hand side of (\ref{4.20}) converges to (\ref{4.29}) in 
$L_u^1\simeq L^1(\mathbf{D},{\rm N})$. In particular, for $\vp\in L_u^1$ and 
$\mu\in\mathbb{M}_m$ we get 
\begin{eqnarray*}
\exp\left\{n\cdot\left(\frac{\mu}{n}\sum_{k=1}^nl_{-k}-\log\omega\right) 
\right\}\vp(r_{-n},w_{-n})\stack{n\to\infty}{\lra}0\quad\mbox{\rm in}\ L_u^1.  
\end{eqnarray*} 

Considering (\ref{4.17}) and the equivalent equation (\ref{4.20}) for 
$\vp,\psi\in L_u^1$ we obtain necessarily (\ref{4.29}) where the infinite 
sum converges in $L_u^1$. On the other hand, we may check representation 
(\ref{4.29}) of $\vp$ back with (\ref{4.17}), now as an equation in 
$L_u^1$. Recalling that $A(\mu)$ is a bounded linear operator in $L_u^1$ 
we establish bijectivity of $id-A(\mu)$ in $L_u^1$ for $\mu\in\mathbb{M}_m$. 
\qed
\medskip

Let 
\begin{eqnarray*}
&&\hspace{-.55cm}W_k(r,w):=\left\{
\begin{array}{ll}
\exp\left\{-k\cdot\left(\frac{\mu}{k}\sum_{i=0}^{k-1}l_i(r,w)-\log\omega\right) 
\right\} & \mbox{\rm if} \ (r,w)\in\mathbf{D}^>_\mu \\ 
-\exp\left\{(k+1)\cdot\left(\frac{\mu}{k+1}\sum_{i=1}^{k+1}l_{-i}(r,w)-\log\omega 
\right)\right\} & \mbox{\rm if} \ (r,w)\in\mathbf{D}^<_\mu 
\end{array}
\right. ,\quad k\in\mathbb{Z}_+,
\end{eqnarray*} 
and, for definiteness, $W_k(r,w):=0$ if $(r,w)\in\{(r,w)\in\partial^{(1)}\Omega 
\times V:w\circ n(r)\ge 0\}\setminus(\mathbf{D}^<_\mu\cup\mathbf{D}^>_\mu)$, 
$k\in\mathbb{Z}$. In addition, let  
\begin{eqnarray*}
&&\hspace{-.5cm}\left(p_k,u_k\right)\equiv\left(p_k(r,w),u_k(r,w)\right) 
\vphantom{\sum} \\
&&\hspace{.5cm}:=\left\{
\begin{array}{ll}
(r_k(r,w),w_k(r,w)) & \mbox{\rm if} \ (r,w)\in\mathbf{D}^>_\mu \\ 
(r_{-k-1}(r,w),w_{-k-1}(r,w)) & \mbox{\rm if} \ (r,w)\in\mathbf{D}^<_\mu 
\end{array}
\right.\, ,\quad k\in\mathbb{Z}_+.  
\end{eqnarray*} 

\noindent
\br{11} For $\mu\in\mathbb{M}$ and any given bounded measurable 
$\psi$ on $\{(r,w)\in\partial^{(1)}\Omega\times V:w\circ n(r)\ge 0\}$, 
the function 
\begin{eqnarray}\label{4.30}
\vp(r,w):=\sum_{k=0}^\infty W_k(r,w)\cdot\psi\left(p_k^-(r,w),u_k(r,w) 
\right) 
\end{eqnarray} 
is a.e. finite, see (\ref{4.28}) as well as (\ref{4.29}) and the 
equivalences (1)-(3) as well as (4)-(6) of Step 3 of the proof of 
Lemma \ref{Lemma4.2}. It is the unique solution to (\ref{4.17}), 
i.e.
\begin{eqnarray*} 
\vp(r,w)-A(\mu)\vp(r,w)\equiv\vp(r,w)-\omega e^{-T_\Omega(r,w)\mu} 
\vp\left(r^-,R_{r^-}(w)\right)=\psi(r,w) 
\end{eqnarray*} 
a.e. on $\{(r,w)\in\partial^{(1)}\Omega\times V:w\circ n(r)\ge 0\}$. 

Furthermore in the context of Lemma \ref{Lemma4.2} (b), for $\psi 
\in L_u^1$, the element $\left( id-A(\mu)\right)^{-1}\psi\in L_u^1$ 
is also given by (\ref{4.30}) where the infinite sum converges in 
$L_u^1$. Note that under the hypotheses of Lemma \ref{Lemma4.2} (b), 
i.e. in particular or $\mu\in\mathbb{M}_m$, we have only the case 
$(r,w)\in\mathbf{D}^>_\mu$, cf. (\ref{4.22}). 
\er

In order to solve (\ref{4.3}) for $\mu$ belonging to some subset 
of $\mathbb{M}_m$ we will now apply the decomposition (\ref{4.10}) 
where we use the bijectivity of $id-A(\mu)$ in $L_u^1$ established in 
Lemma \ref{Lemma4.2} (b). We observe  
\begin{eqnarray}\label{4.31}
&&\hspace{-.5cm}id-a(\mu)=id-A(\mu)-B(\mu)\nonumber \\ 
&&\hspace{.5cm}=(id-A(\mu))\left(id-(id-A(\mu))^{-1}B(\mu)\right)
\end{eqnarray} 
which means that in order to establish invertibility of $id-a(\mu)$ 
we should prove bijectivity of $id-(id-A(\mu))^{-1}B(\mu)$ first. 
In this case we may even conclude that $id-a(\mu)$ is bijective. 
\begin{proposition}\label{Proposition4.3}
Let $c_\tau$, $C_\tau$, and $\mathbb{M}_m$ be given as in Lemma 
\ref{Lemma4.2}. Recall also the notation above Remark \ref{7}. Suppose 
\begin{eqnarray*}
(1-C_{\tau,k_0})\, c_\tau>0\quad\mbox{\rm for some }k_0\in\mathbb{N}.  
\end{eqnarray*} 
(a) Let $\omega\in (0,1)$. For $\mu\in\mathbb{M}_m$ and $f\in L_u^1$ 
the sum 
\begin{eqnarray}\label{4.32}
&&\hspace{-.5cm}X_\mu\, f:=-\frac{(1-\omega)}{\omega}\sum_{k=1}^\infty 
\, \exp\left\{(k-1)\cdot\left(\frac{\mu}{k-1}\sum_{i=1}^{k-1}l_{-i}-\log 
\omega\right)\right\}\times\nonumber \\  
&&\hspace{1.0cm}\times M\left(r_{1-k},R_{r_{1-k}}(w_{1-k})\right)J(r_{ 
1-k},\cdot)(f)\vphantom{ \sum_{k=1}^\infty}
\end{eqnarray} 
converges in $L_u^1$. It defines a bounded linear operator $X_\mu$ in 
$L_u^1$. \\ 
(b) Suppose 
\begin{eqnarray}\label{4.33}
2^{-\frac{1}{k_0}}<\omega<1\, . 
\end{eqnarray} 
Then there exists $-\infty<m_1\le m$ such that for 
\begin{eqnarray*}
\mu\in\mathbb{M}_{m_1}:=\{\lambda\in\mathbb{C}:\mathfrak{Re}\lambda<m_1\} 
\end{eqnarray*} 
we have $\|X_\mu\|_{B(L_u^1,L_u^1)}<1$. For $\mu\in\mathbb{M}_{m_1}$ and 
$\psi\in L_u^1$, the equation 
\begin{eqnarray*}
\vp-(id-A(\mu))^{-1}B(\mu)\vp=\psi
\end{eqnarray*} 
has the unique solution 
\begin{eqnarray}\label{4.34}
\vp=\sum_{n=0}^\infty X_\mu^n\, \psi 
\end{eqnarray} 
where the infinite sum converges in $L_u^1$. Moreover, the operator 
$id-(id-A(\mu))^{-1}B(\mu)$ is bijective in $L_u^1$. 
\end{proposition}
Proof. {\it Step 1 } We show (a). Let $\mu\in\mathbb{M}_m$ and recall 
the definitions of $A(\mu)$ and $B(\mu)$,
\begin{eqnarray*} 
A(\mu)\psi(r,w)=\omega e^{-T_\Omega(r,w)\mu}\psi(r_1,w_1)
\end{eqnarray*}
and, taking into consideration the identity $(r^-,w)=(r_1,R_{r_1}(w_1))$, 
\begin{eqnarray*}
B(\mu)\psi(r,w)=(1-\omega)e^{-T_\Omega(r,w)\mu}M(r_1,R_{r_1}(w_1))J(r_1, 
\cdot)(\psi)\, ,\quad\psi\in L_u^1. 
\end{eqnarray*} 
As a consequence of (\ref{4.2}) and Lemma \ref{Lemma4.1} (a), the 
map $\vp\mapsto M(r^-,w) J(r^-,\cdot)(\vp)$ is a bounded linear operator 
$L_u^1\mapsto L_u^1$. Thus, also $B(\mu)$ is a bounded linear operator 
$L_u^1\mapsto L_u^1$. Furthermore, by Lemma \ref{Lemma4.2} (b), $(id-A 
(\mu))^{-1}$ and therefore also $id-(id-A(\mu))^{-1}B(\mu)$ are bounded 
linear operators $L_u^1\mapsto L_u^1$. 

Now recall that for $c_\tau>0$ and $\mu\in\mathbb{M}_m$ we have $\mathbf 
{D}^>_\mu=\emptyset$ and $\mathbf{D}^=_\mu=\emptyset$ by (\ref{4.22}). 
Using representation (\ref{4.29}) for $(id-A(\mu))^{-1}$ and the above 
definition of $B(\mu)$ it turns out that 
\begin{eqnarray*}
X_\mu=(id-A(\mu))^{-1}B(\mu)
\end{eqnarray*} 
where the sum (\ref{4.32}) converges in $L_u^1$. 
\medskip 

\noindent 
{\it Step 2 } We verify (b). Similar to Step 5 of the proof of 
Lemma \ref{Lemma4.2} for $\psi\in L_u^1$ we obtain 
\begin{eqnarray*}
&&\hspace{-.5cm}\left|(id-A(\mu))^{-1}B(\mu)\psi\right|=\left| 
\sum_{k=1}^\infty\exp\left\{k\cdot\left(\frac{\mu}{k}\sum_{i=1}^k 
l_{-i}-\log\omega\right)\right\}(B(\mu)\psi)(r_{-k},w_{-k})\right| 
 \\ 
&&\hspace{.5cm}\le\frac{(1-\omega)}{\omega}\sum_{k=1}^{k_0}\exp 
\left\{(k-1)\cdot\left(\frac{\mathfrak{Re}\mu}{k-1}\sum_{i=1}^{k-1} 
l_{-i}-\log\omega\right)\right\}\times \\ 
&&\hspace{1.0cm}\times M\left(r_{1-k},R_{r_{1-k}}(w_{1-k})\right) 
J(r_{1-k},\cdot)(|\psi|)\vphantom{\left(\frac{\mathfrak{Re}\mu}{k-1} 
\sum_{i=1}^{k-1}\right)} \\  
&&\hspace{1.0cm}+\frac{(1-\omega)}{\omega}\sum_{k=k_0+1}^\infty\exp 
\left\{(k-1)\cdot\left(\frac{\mathfrak{Re}\mu}{k-1}\sum_{i=1}^{k-1} 
l_{-i}-\log\omega\right)\right\}\times \\ 
&&\hspace{1.0cm}\times M\left(r_{1-k},R_{r_{1-k}}(w_{1-k})\right) 
J(r_{1-k},\cdot)(|\psi|)\vphantom{\left(\frac{\mathfrak{Re}\mu}{k-1} 
\sum_{i=1}^{k-1}\right)} \\  
&&\hspace{.5cm}\le\frac{(1-\omega)}{\omega}\sum_{k=1}^{k_0}\omega^{ 
-(k-1)}M\left(r_{1-k},R_{r_{1-k}}(w_{1-k})\right)J(r_{1-k},\cdot) 
(|\psi|) \\  
&&\hspace{1.0cm}+\frac{(1-\omega)}{\omega}\sum_{k=k_0+1}^\infty e^{ 
(k-1)\cdot\left(\mathfrak{Re}\mu\cdot(1-C_\tau)c_\tau/v_{max}-\log 
\omega\right)}M\left(r_{1-k},R_{r_{1-k}}(w_{1-k})\right)J(r_{1-k}, 
\cdot)(|\psi|)\, .  
\end{eqnarray*} 
Since $\mathbb{T}^{-1}$ preserves the measure ${\rm N}$ we have 
\begin{eqnarray*}
&&\hspace{-.5cm}\int_\mathbf{D}M\left(r_{1-k},R_{r_{1-k}}(w_{1-k}) 
\right)J(r_{1-k},\cdot)(|\psi|)\, d{\rm N}=\int_\mathbf{D}M\left(r, 
R_{r}(w)\right)J(r,\cdot)(|\psi|)\, d{\rm N} \\  
&&\hspace{.5cm}=\int_\mathbf{D}|\psi(r,w)|\, d{\rm N}\, ,\quad\psi 
\in L_u^1,
\end{eqnarray*} 
where we have applied the definition of ${\rm N}$ in (\ref{4.15}) 
and (ii) in the last line. Therefore 
\begin{eqnarray*}
&&\hspace{-.5cm}\int_\mathbf{D}\left|(id-A(\mu))^{-1}B(\mu)\psi 
\right|\, d{\rm N} \\  
&&\hspace{.5cm}\le\frac{(1-\omega)}{\omega}\sum_{k=1}^{k_0}\omega^{ 
-(k-1)}\int_\mathbf{D}|\psi(r,w)|\, d{\rm N} \\  
&&\hspace{1.0cm}+\frac{(1-\omega)}{\omega}\sum_{k=k_0+1}^\infty e^{ 
(k-1)\cdot\left(\mathfrak{Re}\mu\cdot(1-C_\tau)c_\tau/v_{max}-\log 
\omega\right)}\int_\mathbf{D}|\psi(r,w)|\, d{\rm N} \\  
&&\hspace{.5cm}=\frac{(1-\omega)}{\omega}\left(\frac{\omega^{-k_0} 
-1}{\omega^{-1}-1}+\frac{e^{k_0\cdot\left(\mathfrak{Re}\mu 
\cdot(1-C_\tau)c_\tau/v_{max}-\log\omega\right)}}{1-e^{\mathfrak{Re} 
\mu\cdot(1-C_\tau)c_\tau/v_{max}-\log\omega}}\right)\int_\mathbf{D}| 
\psi(r,w)|\, d{\rm N}\, .
\end{eqnarray*} 
Noting that for $\omega\in (0,1)$
\begin{eqnarray*}
\frac{(1-\omega)}{\omega}\cdot\frac{\omega^{-k_0}-1}{\omega^{-1}-1} 
=\omega^{-k_0}-1<1
\end{eqnarray*} 
is equivalent to (\ref{4.33}) we may now claim that there is $m_1\le 
m$ such that for $\mu\in\mathbb{M}_{m_1}=\{\mu\in\mathbb{C}:\mathfrak 
{Re}\mu<m_1\}$ 
\begin{eqnarray*}
\|X_\mu\|_{B(L_u^1,L_u^1)}=\|(id-A(\mu))^{-1}B(\mu)\|_{B(L_u^1,L_u^1)} 
<1\, .
\end{eqnarray*} 

Since $id-A(\mu)$ is bijective in $L_u^1$ by Lemma \ref{Lemma4.2} 
(b), bijectivity of $id-(id-A(\mu))^{-1}B(\mu)$ in $L_u^1$ means that 
for any $\psi\in L_u^1$ there is a unique $\vp\in L_u^1$ such that  
\begin{eqnarray}\label{4.35}
(id-A(\mu))\psi(r,w)=(id-A(\mu))\vp(r,w)-B(\mu)\vp(r,w) 
\end{eqnarray} 
for a.e. $(r,w)\in\partial^{(1)}\Omega\times V$ with $w\circ n(r)\ge 0$. 
We have necessarily obtained (\ref{4.34}). On the other hand, 
checking back representation (\ref{4.34}) of $\vp$ with (\ref{4.35}) 
we verify bijectivity of $id-(id-A(\mu))^{-1}B(\mu)$ in $L_u^1$. 
\qed

\subsection{Spectral Properties on the state space $\Omega$ and the 
Knudsen Type Group}\label{sec:4:2} 

For the readers convenience, let us summarize the facts from spectral 
theory for linear operators in Banach space we use for the subsequent 
Corollary \ref{Corollary4.4} and Theorem \ref{Theorem4.5}.
\br{12}
Let $(X,\|\cdot\|)$ be a Banach space and let $id$ denote the identity 
in $X$. For a closed linear operator $A:X\supseteq D(A)\mapsto X$, let 
$\rho(A):=\{\lambda \in \mathbb{C}: \lambda\, id-A$ is a bijective map 
$D(A) \mapsto X\}$ denote the {\it resolvent set} and let $\sigma(A):= 
\mathbb{C} \backslash \rho(A)$ be the {\it spectrum} of $A$. Denote 
by $P\sigma(A):=\{\lambda\in\mathbb{C}:\lambda\, id-A$ is not 
injective$\}$ the {\it point spectrum} of $A$. 

Furthermore, introduce the {\it approximate point spectrum} $A \sigma 
(A):=\{\lambda\in\mathbb{C}:\lambda\, id-A$ is not injective or the 
range of $\lambda-A$ is not closed in $X\}$. Consider also the {\it 
residual spectrum} $R\sigma(A):=\{\lambda\in\mathbb{C}:$ the range of 
$\lambda -A$ is not dense in $X\}$. Observe that $\sigma(A)=A\sigma(A) 
\cup R\sigma(A)$. 

The {\it Spectral Inclusion Theorem}, Theorem 3.6 of Chapter IV in 
\cite{EN00}, says that for the generator $(A,D(A))$ of a strongly 
continuous semigroup $(T(t))_{t\ge 0}$ on a Banach space $X$, we have 
the inclusions $\sigma(T(t))\supseteq{e}^{t\sigma(A)}$ for 
$t\ge 0$. In particular, for the point-, approximate point-, and 
residual spectrum 
we have 
\begin{eqnarray*}
P \sigma(T(t)) \supseteq {e}^{t P \sigma(A)}\, ,\quad
A \sigma(T(t)) \supseteq {e}^{t A \sigma(A)}\, ,\quad
R \sigma(T(t)) \supseteq {e}^{t R \sigma(A)}\, , 
\end{eqnarray*}
for all $t\ge 0$. In addition the {\it Spectral Mapping Theorem}, Theorem 
3.6 of Chapter IV in \cite{EN00}, for point and residual spectrum states 
the identities
\begin{eqnarray*}
P \sigma(T(t)) \backslash\{0\}={e}^{t P \sigma(A)}\, ,\quad
R \sigma(T(t)) \backslash\{0\}={e}^{t R \sigma(A)}\, ,\quad t\ge 0. 
\end{eqnarray*}
Important for the technical procedure below are Lemma 1.9 and Proposition 
1.10 of Chapter IV in \cite{EN00}. Lemma 1.9 states that for a number 
$\lambda \in \mathbb{C}$ one has $\lambda \in A \sigma(A)$ if and only if 
there exists a sequence $\left(x_{n}\right)_{n \in \mathbb{N}} \subseteq 
D(A)$ such that $\left\|x_{n} \right\|=1$ and $\lim _{n \rightarrow\infty} 
\left\|A x_{n}-\lambda x_{n}\right\|=0$. According to Proposition 1.10, 
the topological boundary $\partial\sigma(A)$ of the spectrum $\sigma(A)$ 
is contained in the approximate point spectrum $A\sigma(A)$.

In the proof of Theorem 4.5 we apply \cite{Pa83}, Theorem 6.5 of Chapter 1, 
according to which a strongly continuous semigroup $T(t)$, $t\ge 0$, of 
bounded linear operators in $X$ can be embedded in a strongly continuous 
group in $X$ if $0\in\rho(T(t_0))$ for some $t_0>0$.
\er

Let us recall that $(A,D(A))$ denotes the infinitesimal generator of 
the strongly continuous semigroup $S(t)$, $t\ge 0$, in $L^1(\Omega 
\times V)$, cf. also Lemma \ref{Lemma3.3} and Lemma \ref{Lemma4.1} 
(b). Note the difference to the operator $A(\mu)$ given by 
(\ref{4.9}). 
\begin{corollary}\label{Corollary4.4} Let $c_\tau$ and $C_\tau$ be 
given as in Lemma \ref{Lemma4.2}. Suppose 
\begin{eqnarray*}
(1-C_{\tau,k_0})\, c_\tau>0\quad\mbox{\rm for some }k_0\in\mathbb{N}   
\end{eqnarray*} 
and $2^{-\frac{1}{k_0}}<\omega<1$. Let $m_1$ and $\mathbb{M}_{m_1}$ 
be as introduced in Proposition \ref{Proposition4.3}. \\ 
(a) For $\mu\in\mathbb{M}_{m_1}$ the operator $id-A(\mu)-B(\mu)=id- 
a(\mu)$ is a bijection in $L_u^1$. \\ 
(b) For $\mu\in\mathbb{M}_{m_1}$ and $g\in L^1(\Omega\times V)$ the 
equation $\mu f-Af=g$ has a unique solution $f\in D(A)$. \\ 
(c) Denoting by $\sigma(A)$ the spectrum of $A$ and by $\sigma(S(t))$ 
the spectrum of $S(t)$ we have the spectral mapping relation 
\begin{eqnarray*}
\sigma(S(t))\setminus\{0\}=\left\{e^{t\mu}:\mu\in\sigma(A)\right\}\, 
,\quad t>0. 
\end{eqnarray*} 
Furthermore, for $t>0$, the resolvent set $\rho(S(t))$ of $S(t)$ 
contains the set $\left\{\lambda=e^{t\mu}:\mu\in\mathbb{M}_{m_1}\right 
\}$ $\cup\left\{\lambda=e^{t\mu}:\mathfrak{Re}\mu>0\right\}$. 
\end{corollary}
Proof. {\it Step 1 } Part (a) is an immediate consequence of (\ref{4.31}), 
Lemma \ref{Lemma4.2} (b), and Proposition \ref{Proposition4.3} (b). 
Furthermore, from (\ref{4.8}) we obtain 
\begin{eqnarray*}
&&\hspace{-.5cm}\int_{r\in\partial^{(1)}\Omega}\int_{w\circ n(r)\ge 0}w\circ 
n(r)\int_0^{T_\Omega(r,w)}\left|e^{-\beta\mu}\right||g(r-\beta w,w)|\, 
d\beta\, dw\, dr \\
&&\hspace{.5cm}\le e^{{\rm diam}(\Omega)\, |\mathfrak{Re}(\mu)|}\int_V 
\|g(\cdot,w)\|_{L^1(\Omega)}\, dw=e^{{\rm diam}(\Omega)\, |\mathfrak{Re} 
(\mu)|}\|g\|_{L^1(\Omega\times V)}<\infty\, ,  
\end{eqnarray*} 
i.e. 
\begin{eqnarray*}
b(\mu)g\in L_u^1\quad \mbox{\rm if }g\in L^1(\Omega\times V).
\end{eqnarray*} 
For part (b) it is now sufficient to note that the equation $\mu f-Af 
=g$ is equivalent to (\ref{4.3}), see Lemma \ref{Lemma4.1} (b). 
\medskip 

\nid
{\it Step 2 } It remains to verify part (c). In this step we prepare  
the proof of the spectral mapping relation $\sigma(S(t))\setminus\{0\} 
=\left\{\lambda=e^{t\mu}:\mu\in\sigma(A)\right\}$, $t>0$. The actual 
proof will be carried out in Steps 3 and 4 below. Let us keep on using 
the symbols and terminology of \cite{EN00}. In particular, let $\sigma 
(\cdot)$ and $A\sigma(\cdot)$ denote the spectrum and approximate point 
spectrum of an operator. According to \cite{EN00}, Theorems 3.6 and 3.7 
of Chapter IV, we have to demonstrate that for any $\lambda\neq 0$ 
belonging to the approximate point spectrum $A\sigma(S(t))$ of $S(t)$ 
we have $\lambda\in\left\{e^{t\mu}:\mu\in A\sigma(A)\right\}$, $t>0$. 
\medskip

Now, let $t>0$ and $e^{t\mu}\in A\sigma(S(t))$. By \cite{EN00}, Lemma 
1.9 of Chapter IV, there exists a sequence $f_n\in L^1(\Omega\times V)$ 
with $\|f_n\|_{L^1(\Omega\times V)}=1$, $n\in\mathbb{N}$, 
and 
\begin{eqnarray}\label{4.36}
\lim_{n\to\infty}\|S(t)f_n-e^{t\mu}f_n\|_{L^1(\Omega\times V)}=0\, . 
\end{eqnarray} 

The continuous Markov process $X_t$, $t\ge 0$, introduced after 
(\ref{3.20}) can be modified in a way that its trajectories never 
reach $\partial\Omega\setminus\partial^{(1)}\Omega$. The trajectories 
move uniformly with constant velocity from time zero until the first 
time they hit $\partial^{(1)}\Omega$. Then they move again uniformly 
with constant velocity until the second time they hit $\partial^{(1)} 
\Omega$, and so on. We also may suppose that the hitting times do 
not accumulate. 
\medskip 

Until (\ref{4.37}) below, suppose we are given a sequence $b_n\in L^1 
(\Omega\times V)$ with $\|b_n\|_{L^1(\Omega\times V)}=1$, and $b_n\ge 0$, 
$n\in\mathbb{N}$. Let $P_{b_n}$ denote the probability measure over the 
$\sigma$-algebra $\mathcal{F}$ generated by the cylindrical sets of the 
trajectories of $X_t$, $t\ge 0$, for which the image measure of $P_{b_n}$ 
under the map $\{X_t,\ t\ge 0\}\mapsto X_0$ is $b_n$ times Lebesgue 
measure on $\overline{\Omega}\times\overline{V}$. Let $E_{b_n}$ stand for 
the expectation with respect to $P_{b_n}$, $n\in\mathbb{N}$. This notation 
allows the interpretation that $E_{b_n}$ is the expectation relative to 
the process $X_t$, $t\ge 0$, started with the initial probability density 
$b_n$.

By $v\le v_{max}<\infty$, cf. Section 2, we have $|X_{t_2}-X_{t_1}| 
\le v_{max}\, |t_2-t_1|$ for all $0\le t_1<t_2$ and all trajectories of 
$X_t$, $t\ge 0$. This implies $E_{b_n}[|X_{t_2}-X_{t_1}|^2]\le v_{max}^2 
\, \left|t_2-t_1\right|^2$. According to \cite{GS69}, Subsections 9.1, 
9.2, and particularly Theorem 2 of Subsection 9.2, there is a subsequence 
$P_{b_{n_k}}$, $k\in\mathbb{N}$, converging weakly to some probability 
measure $P$ on $\mathcal{F}$. Let $\bnu$ denote the image measure of $P$ 
under the map $\{X_t,\ t\ge 0\}\mapsto X_0$. Theorem 2 of Subsection 9.2 
of \cite{GS69} says in particular that for all $t>0$ and all continuous 
maps $g:C([0,t];\overline{\Omega}\times\overline{V})\mapsto\mathbb{R}$ 
the following holds. 
\begin{eqnarray*}
&&\hspace{-.5cm}\lim_{k\to\infty}\int g(\{X_u,\ u\in [0,t]\})\, dP_{b_{ 
n_k}}=\int g(\{X_u,\ u\in [0,t]\})\, dP\, .
\end{eqnarray*} 
Specified to $g(h):=\psi(h(0))$, $h\in C([0,t];\overline{\Omega}\times 
\overline{V})$, $\psi\in C(\overline{\Omega}\times\overline{V})$, it 
holds that 
\begin{eqnarray}\label{4.37}
\lim_{k\to\infty}\int\psi(r,w)b_{n_k}(r,w)\, dr\, dw=\int\psi(r,w)\, 
\bnu(dr\times dw)\, .
\end{eqnarray} 
This gives rise to denote by $E_{\sbnu}$ the expectation with respect 
to $P$. Here we have the interpretation that $E_{\sbnu}$ is the expectation 
relative to the process $X_t$, $t\ge 0$, started with the initial 
probability measure $\bnu$. Thus, for all continuous $g:C([0,t];\overline 
{\Omega}\times\overline{V}))\mapsto\mathbb{R}$ 
\begin{eqnarray}\label{4.38}
&&\hspace{-.5cm}\lim_{k\to\infty}E_{b_{n_k}}[g(\{X_u,\ u\in [0,t]\})]= 
E_{\sbnu}[g(\{X_u,\ u\in [0,t]\})]\, .
\end{eqnarray} 

Let us now return to the sequence $f_n\in L^1(\Omega\times V)$ with 
$\|f_n\|_{L^1(\Omega\times V)}=1$, $n\in\mathbb{N}$, introduced in the 
beginning of this step. Then (\ref{4.38}) and (\ref{4.37}) hold for 
$b_{n_k}$ replaced with $f_{n_k}^+/\|f_{n_k}^+\|_{L^1(\Omega\times V)}$ 
or $f_{n_k}^-/\|f_{n_k}^-\|_{L^1(\Omega\times V)}$, where $f_{n_k}^+:=
f_{n_k}\vee 0$ as well as $f_{n_k}^-:=(-f_{n_k})\vee 0$, and $\bnu$ 
replaced with the respective non-negative Borel probabiliy measures on 
$\overline{\Omega}\times\overline{V}$, $\bnu^+$ or $\bnu^-$. Choosing 
another subsequence if necessary, we may suppose that the limits 
\begin{eqnarray*}
l^+:=\lim_{k\to\infty}\|f_{n_k}^+\|_{L^1(\Omega\times V)}\quad\mbox{\rm 
and}\quad l^-:=\lim_{k\to\infty}\|f_{n_k}^-\|_{L^1(\Omega\times V)}= 
1-l^+
\end{eqnarray*} 
exist. It follows now from (\ref{4.38}) and (\ref{4.37}) that with 
$\bmu:=l^+\bnu^+-l^-\bnu^-$ we have 
\begin{eqnarray*}
\lim_{k\to\infty}\int\psi(r,w)f_{n_k}(r,w)\, dr\, dw=\int\psi(r,w)\, 
\bmu(dr\times dw)
\end{eqnarray*} 
and
\begin{eqnarray}\label{4.39}
&&\hspace{-.5cm}\lim_{k\to\infty}E_{f_{n_k}}[g(\{X_u,\ u\in [0,t]\})] 
\, dr\, dw=E_{\sbmu}[g(\{X_u,\ u\in [0,t]\})] 
\end{eqnarray} 
for $g$ and $\psi$ as above. Here the expectations $E_{f_{n_k}}$ and 
$E_{\sbmu}$ come with the interpretation as weighted differences of the 
expectations relative to the process $X_t$, $t\ge 0$, started with the 
respective initial probability densities and initial probability measures.

Below we shall consider two cases separately, namely that $\bmu$ is the 
zero measure, or alternatively that $\bmu$ is not the zero measure, i.e. 
a certain  signed measure with total variation not greater than one. 
\medskip 

\nid
{\it Step 3 } Assume in this step that $\bmu$ is the zero measure. 
Recall Steps 1 and 2 of the proof of Lemma \ref{Lemma3.3} and the 
notation introduced there. In particular, for $\vp\in C_b(\Omega 
\times V)$, uniformly continuous on $\Omega\times V$ recall $a_\vp(t) 
:=\sup_{v\in V}\sup_{r\in\Omega_\varepsilon}|\vp(r-tv,v)-\vp(r,v)|$. 
Similar to (\ref{3.24}) we obtain
\begin{eqnarray*}
&&\hspace{-.5cm}\left|\int_\Omega\int_V\left(S(u)f_{n_k}-f_{n_k}\right) 
\vp\, dv\, dr\right|\le 3\|\vp\|\int_{\Omega\setminus\Omega_{2\varepsilon}} 
\int_V |f_{n_k}|\, dv\, dr+a_\vp(u)\, ,\quad k\in\mathbb{N},\ u\in (0,t]. 
\end{eqnarray*} 
Noting that the variation $|\bmu|$ of $\bmu$ is also the zero measure 
it follows that 
\begin{eqnarray*}
&&\hspace{-.5cm}\limsup_{k\in\mathbb{N}}\left|\int_\Omega\int_V\left( 
S(u)f_{n_k}-f_{n_k}\right)\vp\, dv\, dr\right| \\ 
&&\hspace{.5cm}\le\limsup_{k\in\mathbb{N}} 3\|\vp\|\int_{\overline 
{\Omega}\setminus\Omega_{2\varepsilon}}\int_{\overline{V}}|f_{n_k}|\, 
dv\, dr+a_\vp(u)\le 3\|\vp\|\cdot|\bmu|\left(\left({\overline{\Omega}} 
\setminus\Omega_{2\varepsilon}\right)\times\overline{V}\right)+a_\vp(u) 
 \\ 
&&\hspace{.0cm}\stack{u\to 0}{\lra}0\vphantom{\int_{\overline{V}}}\, .
\end{eqnarray*} 
For the $\le$ sign in the second last line, observe that $\left({\overline 
{\Omega}}\setminus\Omega_{2\varepsilon}\right)\times\overline{V}$ is closed 
and apply the Portmanteau theorem. Similarly, we get 
\begin{eqnarray*}
&&\hspace{-.5cm}\limsup_{k\in\mathbb{N}}\left|\int_\Omega\int_V\left( 
S(s+u)f_{n_k}-S(s)f_{n_k}\right)\vp\, dv\, dr\right|\le a_\vp(u)\stack 
{u\to 0}{\lra}0
\end{eqnarray*} 
for all $0<u<t$, $0<s<t$ with $s+u<t$. The last two relations imply that 
for any set $\Phi$ of test functions containing all $\vp\in C_b(\Omega 
\times V)$ having a common modulus of continuity and a common bound on 
the norm in $C_b(\Omega\times V)$, the family 
\begin{eqnarray*}
[0,t]\ni u\mapsto\int_\Omega\int_V\left(S(u)f_{n_k}\right)\vp_k\, dv 
\, dr\, ,\quad \vp_k\in\Phi,\ k\in\mathbb{N}, 
\end{eqnarray*} 
is equicontinuous and equibounded. Keeping in mind that every such set 
$\Phi$ is totally bounded in $C_b(\Omega\times V)$ and that there is an 
increasing sequence of such sets $\Phi$ whose union is dense in $C_b( 
\Omega\times V)$ we may select a set $\Phi$ and a sequence $\vp_k\in 
\Phi$ such that $\int_\Omega\int_V f_{n_k}\vp_k\, dv\, dr>\frac12$. It 
follows now precisely as in \cite{EN00}, proof of Lemma 3.9 of Chapter 
IV, part (a) $\Rightarrow$ (b), that there is an $m\in\mathbb{Z}$ such 
that $\mu+\frac{2\pi m\cdot i}{t}\in A\sigma(A)$ provided that $\lambda 
=e^{t\mu}=e^{t(\mu+\frac{2\pi m\cdot i}{t})}\in A\sigma(S(t))$. As 
explained in the first paragraph of Step 2, this implies the spectral 
mapping relation $\sigma(S(t))\setminus\{0\}=\left\{e^{t\mu}:\mu\in 
\sigma(A)\right\}$, $t>0$.
\medskip

\nid
{\it Step 4 } Assume now that $\bmu$ is not the zero measure. For $m 
\in\mathbb{Z}$ and $\psi\in C\left(\overline{\Omega}\times\overline{V} 
\right)$ equations (\ref{3.21}) and (\ref{4.39}) imply the existence 
of the limit
\begin{eqnarray}\label{4.40}
&&\hspace{-.5cm}\left|\int_{[0,t]\times\overline{\Omega}\times 
\overline{V}}e^{-u\left(\frac{2\pi m\cdot i}{t}\right)}e^{-u\mu} 
S(u)f_{n_k}(r,w)\cdot\psi(r,w)\, dw\, dr\, du\right|\nonumber \\  
&&\hspace{.5cm}=\left| 
E_{f_{n_k}}\left(\int_0^te^{-u\left(\frac{2\pi m\cdot i}{t}\right)} 
e^{-u\mu}\cdot \psi(X_u)\, du\right)\right|\nonumber \\  
&&\hspace{.0cm}\stack{k\to\infty}{\lra}\left|E_{\sbmu}\left(\int_0^t 
e^{-u\left(\frac{2\pi m\cdot i}{t}\right)}e^{-u\mu}\cdot\psi(X_u)\, du 
\right)\right|\vphantom{\int_0^t}\, . 
\end{eqnarray} 
By the Stone-Weierstra\ss\ Theorem, the (complex) $C\left([0,t]\times 
\overline{\Omega}\times\overline{V}\right)$ coincides with the closed 
linear span of 
\begin{eqnarray*} 
\left\{e^{-u\left(\frac{2\pi m\cdot i}{t}\right)}\cdot\psi:m\in\mathbb 
{Z},\ \psi\in C\left(\overline{\Omega}\times\overline{V}\right)\right\} 
\, . 
\end{eqnarray*} 
Since $\bmu$ is not the zero measure in (\ref{4.40}) there exist a 
particular $m\in\mathbb{Z}$ and a particular $\psi\in C\left( 
\overline{\Omega}\times\overline{V}\right)$ with $\|\psi\|=1$ such 
that 
\begin{eqnarray}\label{4.41} 
&&\hspace{-.5cm}\lim_{k\to\infty}\left|\int_{[0,t]\times\overline 
{\Omega}\times\overline{V}}e^{-u\left(\frac{2\pi m\cdot i}{t}\right)} 
e^{-u\mu}S(u)f_{n_k}(r,w)\cdot\psi(r,w)\, dw\, dr\, du\right|>0\, . 
\end{eqnarray} 
As a consequence there is $c>0$ and $k_0\in\mathbb{N}$ such that for 
all $k>k_0$ we have 
\begin{eqnarray}\label{4.42}
&&\hspace{-.5cm}\left\|\int_0^t e^{-u\left(\frac{2\pi m\cdot i}{t} 
\right)}e^{-u\mu}S(u)f_{n_k}\, du\right\|_{L^1(\Omega\times V)}\nonumber 
 \\ 
&&\hspace{.5cm}\ge\left|\int_{[0,t]\times\overline{\Omega}\times\overline 
{V}}e^{-u\left(\frac{2\pi m\cdot i}{t}\right)}e^{-u\mu}S(u)f_{n_k}(r,w) 
\cdot\psi(r,w)\, dw\, dr\, du\right|\nonumber \\  
&&\hspace{.5cm}>c\, .\vphantom{\int_0^t}
\end{eqnarray} 
On the other hand, for $m\in\mathbb{Z}$ chosen in (\ref{4.41}), it 
holds that 
\begin{eqnarray}\label{4.43}
&&\hspace{-.5cm}f_{n_k}-e^{-t\mu}S(t)f_{n_k}=f_{n_k}-e^{-t\left(\mu 
+\frac{2\pi m\cdot i}{t}\right)}S(t)f_{n_k}\vphantom{\int_0^t} 
\nonumber \\ 
&&\hspace{.5cm}=\left(\left(\mu+\frac{2\pi m\cdot i}{t}\right)\, id- 
A\right)\cdot\int_0^te^{-u\left(\frac{2\pi m\cdot i}{t}\right)}e^{-u 
\mu}S(u)f_{n_k}\, du\, ,\quad k\in\mathbb{N}, 
\end{eqnarray} 
cf. \cite{EN00}, Lemma 1.9 of Chapter II. Introduce 
\begin{eqnarray*}
\hat{f}_k:=\int_0^te^{-u\left(\frac{2\pi m\cdot i}{t}\right)}e^{-u\mu} 
S(u)f_{n_k}\, du\, ,\quad k\in\mathbb{N}.  
\end{eqnarray*} 
Relation (\ref{4.42}) says that $\liminf_{k\to\infty}\|\hat{f}_k\|_{ 
L^1(\Omega\times V)}>c>0$. Imposing now (\ref{4.36}) on the left-hand 
side of (\ref{4.43}) we verify $\lim_{k\to\infty}((\mu+\frac{2\pi m 
\cdot i}{t})\hat{f}_k-A\hat{f}_k)=0$ in $L^1(\Omega\times V)$ and, 
moreover with $F_k:=\hat{f}_k/\|\hat {f}_k\|_{L^1(\Omega\times V)}$ 
for sufficiently large $k\in\mathbb{N}$, 
\begin{eqnarray*}
\lim_{k\to\infty}\left\|\left(\mu+\frac{2\pi m\cdot i}{t}\right)F_k-
AF_k\right\|_{L^1(\Omega\times V)}=0\, . 
\end{eqnarray*} 
Also in the case when $\bmu$ is not the zero measure, we have 
demonstrated that $\mu+\frac{2\pi m\cdot i}{t}\in A\sigma(A)$ 
provided that $\lambda=e^{t\mu}=e^{t(\mu+\frac{2\pi m\cdot i}{t})} 
\in A\sigma(S(t))$ where $m\in\mathbb{Z}$ is the number chosen in 
(\ref{4.41}). In other words, together with \cite{EN00}, Theorems 3.6 
and 3.7 of Chapter IV, we have proved the spectral mapping relation 
\begin{eqnarray}\label{4.44}
\sigma(S(t))\setminus\{0\}=\left\{e^{t\mu}:\mu\in\sigma(A)\right\}\, , 
\quad t>0. 
\end{eqnarray} 

\nid
{\it Step 5 } For the remainder of part (c) note that $\rho(S(t))$ 
contains the set $\left\{\lambda=e^{t\mu}:\mu\in\mathbb{M}_{m_1} 
\right\}$ because of part (b) of this corollary and the spectral 
mapping relation (\ref{4.44}). 

Recall also from Lemma \ref{Lemma3.3} and the construction in (iii) 
that $S(t)$, $t\ge 0$, is a strongly continuous semigroup with 
operator norm $\|S(t)\|=1$ in $B(L^1(\Omega\times V),L^1(\Omega 
\times V))$ for all $t\ge 0$. Consequently, $\rho(S(t))$ contains 
the set $\left\{\lambda=e^{t\mu}:\mathfrak{Re}\mu>0\right\}$ for all 
$t>0$ by the Hille-Yosida Theorem and again the spectral mapping 
relation (\ref{4.44}). 
\qed
\medskip

Let us conclude this section with some results concerning the existence 
of $S(t)$ for time $t\in\mathbb{R}$. 
\begin{theorem}\label{Theorem4.5} 
Let $c_\tau$ and $C_\tau$ be given as in Lemma \ref{Lemma4.2}. Suppose 
\begin{eqnarray*}
(1-C_{\tau,k_0})\, c_\tau>0\quad\mbox{\rm for some }k_0\in\mathbb{N}   
\end{eqnarray*} 
and (\ref{4.33}), i.e. 
\begin{eqnarray*}
2^{-\frac{1}{k_0}}<\omega<1\, . 
\end{eqnarray*} 
Then $S(t)$, $t\ge 0$, extends to a strongly continuous group in 
$L^1(\Omega\times V)$ which we will denote by $S(t)$, $t\in\mathbb{R}$. 
\end{theorem}
Proof. As an immediate consequence of Lemma \ref{Lemma3.3} and the 
classical result \cite{Pa83}, Theorem 6.5 of Chapter 1, we have to show 
that $\lambda=0$ does not belong to the spectrum of $S(t_0)$ for some 
$t_0>0$.

\nid
{\it Step 1 } Assuming the contrary it follows from Corollary 
\ref{Corollary4.4} (c) that $\lambda=0$ is an isolated point in the 
spectrum of $S(t)$ for all $t>0$. By \cite{EN00}, Proposition 1.10 in 
Chapter IV, $\lambda=0$ belongs to the approximate spectrum of $S(t)$, 
$t>0$, in the context of \cite{EN00}, Chapter IV. This means by 
\cite{EN00}, Lemma 1.9 of Chapter IV, that for all $t>0$ there is a 
sequence $f_n\in L^1(\Omega\times V)$ with 
\begin{eqnarray}\label{4.45}
\|f_n\|_{L^1(\Omega\times V)}=1\, ,\ n\in\mathbb{N},\ \mbox{\rm and } 
\|S(t)f_n\|_{L^1(\Omega\times V)}\stack{n\to\infty}{\lra}0. 
\end{eqnarray}
{\it Step 2 } It is now our turn to demonstrate that (\ref{4.45})
cannot hold. For this let $0<t<k_0(1-C_{\tau,k_0})\, c_\tau/v_{ 
max}$ which implies by the definitions of $c_\tau$ and $C_\tau$ 
\begin{eqnarray*}
0<t<k_0(1-C_{\tau,k_0})\, c_\tau/v_{max}\le k_0(1-C_{\tau,k_0})\tau 
(r,w)\le\sum_{i=1}^{k_0}l_i(r,w)\, ,\quad (r,w)\in\mathbf{D}\, .
\end{eqnarray*}
Figuratively, this means that particles just being deterministically 
reflected at the boundary during the time interval $(0,t)$, hit the 
boundary $\partial\Omega$ maximally $k_0$ times. Mathematically, we 
have the explicit representation of $S(t)$ in (\ref{3.22}), see 
Remark \ref{6}. Replacing now in (\ref{3.22}) $p_0$ by an 
arbitrary $f_0\in L^1(\Omega\times V)$ and letting $t$ be fixed as 
above, representation (\ref{3.22}) holds for a.e. $(r,v)\in\Omega 
\times V$. 

Recall the terminology of Definition \ref{Definition3.2}. It follows 
from (\ref{3.22}) and the above choice of $t$ that  
\begin{eqnarray}\label{4.46}
&&\hspace{-.5cm}S(t)f_0(r,v)=\chi_{\{0\}}(m)f_0(r_e,v_0)+\omega 
\chi_{\{1\}}(m) f_0(r_e,R_{r_{1}}(v_0))\vphantom{\dot{f}}\nonumber 
 \\ 
&&\hspace{1.0cm}+\ldots+\omega^{k_0}\chi_{\{k_0\}}(m)f_0(r_e,R_{ 
r_{k_0}}(\ldots R_{r_{1}}(v_0)\ldots ))\vphantom{\dot{f}} 
\nonumber \\ 
&&\hspace{1.0cm}+\mathbf{R}(f_0;r,v,t)\quad\mbox{\rm a.e. on } 
(r,v)\in\Omega\times V\vphantom{\dot{f}}
\end{eqnarray}
where the term $\mathbf{R}(f_0;r,v,t)$ contains all the remaining 
items of the right-hand side of (\ref{3.22}) not appearing in the 
first two lines of (\ref{4.46}). Introduce also
\begin{eqnarray*}
&&\hspace{-.5cm}\mathbf{Q}(f_0;r,v,t):=\chi_{\{0\}}(m)f_0(r_e,v_0)+ 
\omega\chi_{\{1\}}(m)f_0(r_e,R_{r_{1}}(v_0))\nonumber \\ 
&&\hspace{1.0cm}+\ldots+\omega^{k_0}\chi_{\{k_0\}}(m)f_0(r_e,R_{r_{ 
k_0}}(\ldots R_{r_{1}}(v_0)\ldots )) 
\end{eqnarray*}
where we recall that all terms which contribute to $\mathbf{R}(f_0; 
\cdot,\cdot,t)$ as well as $\mathbf{Q}(f_0;\cdot,\cdot,t)$, except 
for $\omega$, are functions of $(r,v)\in\Omega\times V$. We have 
\begin{eqnarray}\label{4.47}
&&\hspace{-.5cm}\|\mathbf{R}(f_0;\cdot,\cdot,t)\|_{L^1(\Omega\times 
V)}\le\|\mathbf{R}(|f_0|;\cdot,\cdot,t)\|_{L^1(\Omega\times V)}
\nonumber \\ 
&&\hspace{.5cm}=\|S(t)|f_0|\|_{L^1(\Omega\times V)}-\|\mathbf{Q}(|f_0| 
;\cdot,\cdot,t)\|_{L^1(\Omega\times V)}
\end{eqnarray}
since, replacing in (\ref{4.46}) $f_0$ with $|f_0|$, all terms there 
are non-negative. In the context of Definition \ref{Definition3.2} 
we keep in mind the following. For $(r,v),(\bar{r},\bar{v})\in\Omega 
\times V$ with $(r,v)\neq (\bar{r},\bar{v})$ let us follow the two 
paths $\pi$ with time range $[0,t]$ pinned at time $t$ in $(r,v)$ and 
$(\bar{r},\bar{v})$ generated by deterministic reflections at the 
boundary $\partial\Omega$. Supposing that these paths do not terminate 
in an edge or vertex of $\partial\Omega$, we observe that they do not 
intersect in both, the space as well as the velocity variable, at any 
time between $t$ and $0$. Consequently, there is an one-to-one map 
a.e. on $\Omega\times V$ between the starting point $(r,v)\equiv (r_0, 
v_0)$ at time $t$ and the end point 
\begin{eqnarray}\label{4.48}
\left\{
\begin{array}{cl}
(r_e,v_0) & \mbox{\rm if no reflection between $t$ and $0$} \\  
(r_e,R_{r_{k}}(\ldots R_{r_{1}}(v_0)\ldots ))) &  \mbox{\rm if $k$ 
reflections between $t$ and $0$} 
\end{array} 
\right.\, ,\quad k\in\{1,\ldots,k_0\}, 
\end{eqnarray} 
at time 0 of those deterministic paths. For use in (\ref{4.50}) below, 
we observe also that the volume element $ dr\, dv$ is preserved under 
this map. Since the map between $(r,v)$ and (\ref{4.48}) is one-to-one, 
by the definition of $\mathbf{Q}$ we have 
\begin{eqnarray*}
\|\mathbf{Q}(|f_0|;\cdot,\cdot,t)\|_{L^1(\Omega\times V)}=\|\mathbf{Q} 
(f_0;\cdot,\cdot,t)\|_{L^1(\Omega\times V)}  
\end{eqnarray*}
which together with (\ref{4.47}) gives 
\begin{eqnarray}\label{4.49}
&&\hspace{-.5cm}\|\mathbf{R}(f_0;\cdot,\cdot,t)\|_{L^1(\Omega\times V)} 
+\|\mathbf{Q}(f_0;\cdot,\cdot,t)\|_{L^1(\Omega\times V)}\nonumber \\ 
&&\hspace{0.5cm}\le\|S(t)|f_0|\|_{L^1(\Omega\times V)}=\||f_0|\|_{L^1( 
\Omega\times V)}=\|f_0\|_{L^1(\Omega\times V)}\, . 
\end{eqnarray}
It follows now from (\ref{4.46}) and (\ref{4.49}) and the just mentioned 
preservation of the volume element $dr\, dv$ that 
\begin{eqnarray}\label{4.50}
&&\hspace{-.5cm}\|S(t)f_0\|_{L^1(\Omega\times V)}\ge 2\|\mathbf{Q}(f_0;
\cdot,\cdot,t)\|_{L^1(\Omega\times V)}-\|f_0\|_{L^1(\Omega\times V)} 
\nonumber \\ 
&&\hspace{0.5cm}\ge 2\omega^{k_0}\|\chi_{\{0\}}(m)f_0(r_e,v_0)+\chi_{\{1 
\}}(m)f_0(r_e,R_{r_{1}}(v_0))\nonumber \\ 
&&\hspace{1.0cm}+\ldots+\chi_{\{k_0\}}(m)f_0(r_e,R_{r_{k_0}}(\ldots R_{ 
r_{1}}(v_0)\ldots ))\|_{L^1(\Omega\times V)}-\|f_0\|_{L^1(\Omega\times V)} 
\nonumber \\ 
&&\hspace{0.5cm}\ge (2\omega^{k_0}-1)\|f_0\|_{L^1(\Omega\times V)}
\end{eqnarray}
where we note that by hypothesis (\ref{4.33}) we have $2\omega^{k_0}-1>0$. 
Thus, (\ref{4.50}) contradicts (\ref{4.45}). 
\qed 
\medskip

As a continuation of Lemma \ref{Lemma3.1}, we are also interested in 
boundedness from below and above along the group $S(t)$, $t\in\mathbb{R}$. 
\begin{corollary}\label{Corollary4.6} 
Let the conditions of Theorem \ref{Theorem4.5} be satisfied. \\ (a) Let 
$p_0\in L^1(\Omega\times V)$ and $t_0\in\mathbb{R}$ such that $S(t_0)p_0 
\in L^\infty(\Omega\times V)$. Then there are finite real numbers $p_{0, 
{\rm min}}$ and $p_{0,{\rm max}}$ such that 
\begin{eqnarray*}
p_{0,{\rm min}}\le S(t)\, p_0\le p_{0,{\rm max}}\quad\mbox{\rm a.e. on } 
\Omega\times V
\end{eqnarray*} 
for all $t\ge t_0$. In particular, if $S(t_0)p_0\ge 0$ a.e. and $\|1/ 
S(t_0)p_0\|_{L^\infty(\Omega\times V)}<\infty$ then there exists $p_{ 
0,{\rm min}}>0$. \\ 
(b) Let $p_0\in L^1(\Omega\times V)$. We have $\int_\Omega\int_V S(t) 
p_0(r,v)\, dv\, dr=\int_\Omega\int_V p_0(r,v)\, dv\, dr=:\langle p_0 
\rangle$ for all $t\in\mathbb{R}$. \\ 
(c) There exists $\alpha>0$ such that for all $p_0\in L^1(\Omega\times 
V)$ 
\begin{eqnarray}\label{4.51}
\|S(t)\, (p_0-\langle p_0\rangle\, \overline{g})\|_{L^1(\Omega\times V) 
}\le e^{-\alpha t}\|p_0-\langle p_0\rangle\, \overline{g}\|_{L^1(\Omega 
\times V)}\, ,\quad t\ge 0, 
\end{eqnarray} 
where $\overline{g}$ is the unique stationary probability density of 
$S(t)$, $t\in\mathbb{R}$, cf. (\ref{3.1}). In addition there is $M 
\ge 1$ such that with $m_1$ given by Proposition \ref{Proposition4.3}  
(b) we have for all $p_0\in L^1(\Omega\times V)$ 
\begin{eqnarray}\label{4.52}
\|S(t)\, p_0\|_{L^1(\Omega\times V)}\le Me^{m_1 t}\|p_0\|_{L^1(\Omega 
\times V)}\, ,\quad t\le 0.  
\end{eqnarray} 
\end{corollary}
Proof. Part (a) follows from Lemma \ref{Lemma3.1}, part (b) is a standard 
consequence of the non-negativity and linearity of the semigroup $S(t)$, 
$t\ge 0$. For the estimate (\ref{4.51}) we mention that it can be 
established  by (\ref{3.2}) and the calculations similar to those 
between (3.9) and (3.11) in \cite{CPW98}. Keeping in mind Corollary 
\ref{Corollary4.4} (b), (c), and Theorem \ref{Theorem4.5}, the estimate 
(\ref{4.52}) is standard, see e.g. \cite{EN00}, Proposition 2.2 of 
Chapter IV and its proof. Note here that $\sigma(-A)=\{-\mu:\mu\in\sigma 
(A)\}$. 
\qed

\section{Solutions to the Boltzmann Type Equation}\label{sec:5}
\setcounter{equation}{0} 

We are interested in global solutions to the Boltzmann type equations  
(\ref{2.2}) and (\ref{2.4}) for $t\ge 0$ and $t\in\mathbb{R}$. It is 
beneficial to construct local solutions in a first step. In particular, 
the proof of the existence of global solutions on $t\in\mathbb{R}$ uses 
crucially the Knudsen type group $S(t)$, $t\in\mathbb{R}$. The initial 
values $p_0$ at time zero are probability densities on $\Omega\times V$ 
satisfying $c\le p_0\le C$ for some $0<c\le C<\infty$. In this section 
we shall suppose the global conditions (i)-(vii). 

\subsection{Construction of Local Solutions to the Boltzmann Type 
Equation} 
\label{sec:5:1}

By the {\it normalization condition} in (ii), $S(t)$, $t\ge 0$, given 
in (iii) maps $L^1(\Omega\times V)$ linearly to $L^1(\Omega\times V)$ 
with operator norm one. We observe furthermore that by the definitions 
in (iv)-(vi), for fixed $t\ge 0$, $Q$ maps $(p(\cdot,\cdot,t),q(\cdot, 
\cdot,t))\in L^1(\Omega\times V)\times L^1(\Omega\times V)$ to $L^1 
(\Omega\times V)$ such that 
\begin{eqnarray}\label{5.1}
\int_\Omega\int_V Q(p(\cdot,\cdot,t),q(\cdot,\cdot,t))(r,v)\, dv\, dr=0
\end{eqnarray} 
and 
\begin{eqnarray}\label{5.2}
&&\hspace{-.5cm}\|Q(p(\cdot,\cdot,t),q(\cdot,\cdot,t))\|_{L^1(\Omega\times V)} 
\vphantom{\dot{f}}\nonumber \\ 
&&\hspace{.5cm}\le 2\|h_\gamma\|\|B\|\|p(\cdot,\cdot,t)\|_{L^1(\Omega\times V)} 
\|q(\cdot,\cdot,t)\|_{L^1(\Omega\times V)}
\end{eqnarray} 
where $\|\cdot\|$ denotes the sup-norm. 
\medskip 

Let $T>0$ and let $(L^1(\Omega\times V))^{[0,T]}$ be the space of all measurable 
real functions $f(r,v,t)$, $(r,v)\in\Omega\times V$, $t\in [0,T]$, such that $F_f 
(t):=f(\cdot,\cdot,t)\in L^1(\Omega\times V)$ for all $t\in [0,T]$ and $F_f\in C( 
[0,T];L^1(\Omega\times V))$. With the norm 
\begin{eqnarray*}
\|f\|_{1,T}:=\sup_{t\in [0,T]}\|f(\cdot,\cdot,t)\|_{L^1(\Omega\times V)}
\end{eqnarray*} 
$(L^1(\Omega\times V))^{[0,T]}$ is a Banach space. 

Let $p_0\in L^1(\Omega\times V)$ and let $p,q\in (L^1(\Omega\times V))^{[0, 
T]}$. For $(r,v,t)\in \Omega\times V\times [0,T]$ we set 
\begin{eqnarray}\label{5.3}
\Psi(p_0,p)(r,v,t):=S(t)\, p_0(r,v)+\lambda\int_0^t S(t-s)\, Q(p,p)\, (r,v,s) 
\, ds 
\end{eqnarray} 
where the integral converges in $L^1(\Omega\times V)$. We have $\Psi(p_0,p) 
\in (L^1(\Omega\times V))^{[0,T]}$ and according to (\ref{5.2})
\begin{eqnarray}\label{5.4} 
&&\hspace{-.5cm}\|\Psi(p_0,p)\|_{1,T}\le\|p_0\|_{L^1(\Omega\times V)}+\lambda 
T\|Q(p,p)\|_{1,T}\nonumber \\ 
&&\hspace{.5cm}\le\|p_0\|_{L^1(\Omega\times V)}+2\lambda T\|h_\gamma\|\|B\| 
\|p\|_{1,T}^2\, . 
\end{eqnarray} 
Moreover, if $\lambda\le 1/(16 T\|h_\gamma\|\|B\|)$ and $\|p_0\|_{L^1(\Omega 
\times V)}\le 3/2$ then (\ref{5.4}) leads to 
\begin{eqnarray}\label{5.5}
\|p\|_{1,T}\le 2\quad {\rm implies}\quad \|\Psi (p_0,p)\|_{1,T}\le 2\, .  
\end{eqnarray} 

Using symmetry and bilinearity of $Q(\cdot ,\cdot)$, for $p_0,p_0'\in L^1 
(\Omega\times V)$ and $p_1,p_2\in (L^1(\Omega\times V))^{[0,T]}$ we obtain 
from (\ref{5.2}), specified to $p=p_1+p_2$ and $q=p_1-p_2$, 
\begin{eqnarray}\label{5.6}
&&\hspace{-.5cm}\|\Psi(p_0,p_1)-\Psi(p_0',p_2)\|_{1,T}\nonumber \\ 
&&\hspace{.5cm}\le\max_{t\in [0,T]}\|S(t)\, (p_0-p_0')\|_{L^1(\Omega\times 
V)}+\lambda T\|Q(p_1+p_2,p_1-p_2)\|_{1,T} \nonumber \\ 
&&\hspace{.5cm}\le\|p_0-p_0'\|_{L^1(\Omega\times V)}+2\lambda T\|h_\gamma 
\|\|B\|\|p_1+p_2\|_{1,T}\|p_1-p_2\|_{1,T}\, .  
\end{eqnarray} 

It follows from (ii) and (iii) that $\frac{d}{dt}\int_\Omega\int_VS(t) 
\, p_0(r,v)\, dr\, dv$ $=0$, $t\ge 0$. Relation (\ref{5.1}) gives now  
\begin{eqnarray}\label{5.7}
\frac{d}{dt}\int_\Omega\int_V\Psi(p_0,p)(r,v,t)\, dr\, dv=0\, ,\quad t\in 
[0,T].   
\end{eqnarray} 

Let $\1$ denote the function constant to one on $\Omega\times V$. Furthermore, 
for $t\in [0,\infty]$ introduce 
\begin{eqnarray*}
c^{\1}_{t,\rm max}:=\sup\|S(\tau)\, \1\|_{L^\infty(\Omega\times V)} 
\end{eqnarray*} 
where, for $0\le t<\infty$ the supremum is taken over $\tau\in [0,t]$ and 
for $t=\infty$ the supremum is taken over $\tau\in\mathbb{R}_+$. 
Corollary \ref{Corollary4.6} yields $0<c^{\1}_{t,\rm max}<\infty$, $t\in [0, 
\infty]$. 

By (iv)-(vi) it holds for non-negative $p,q\in (L^1(\Omega\times V))^{[0,T]}$ 
that  
\begin{eqnarray}\label{5.8}
Q(p,q)(r,v,t)\ge -\|h_\gamma\|\|B\|\cdot\frac12\left(\|p\|_{1,T}\, q(r,v,t) 
+\|q\|_{1,T}\, p(r,v,t)\right)\, ,
\end{eqnarray} 
$r\in\Omega$, $v\in V$, $t\in [0,T]$. Let us assume 
\begin{eqnarray*}
\hat{p}_{t,{\rm max}}:=\, \mbox{\rm ess}\, \sup\{p(r',v',t'):r'\in\Omega, 
\ v'\in V,\ t'\in [0,t]\}<\infty 
\end{eqnarray*} 
and $\|p\|_{1,T}\le 1$. Relation (\ref{5.8}) together with definition 
(\ref{5.3}) and Lemma \ref{Lemma3.1} applied to $p_0$ as well as $p(\cdot, 
\cdot,s)$, $s\in [0,t]$, imply that $0\le\|S(\tau-s)p(\cdot,\cdot,s)\|_{L^1 
(\Omega\times V)}\le\|S(\tau-s)(\hat{p}_{t,{\rm max}}\cdot\1)\|\le\hat{p}_{ 
t,{\rm max}}\, c^{\1}_{t,\rm max}$ for all $0\le s\le\tau\le t$ and hence 
\begin{eqnarray}\label{5.9}
\Psi(p_0,p)(r,v,t)\ge p_{0,{\rm min}}-\lambda t\|h_\gamma\|\|B\|\cdot 
\hat{p}_{t,{\rm max}}\, c^{\1}_{t,\rm max}\, , 
\end{eqnarray} 
$r\in\Omega$, $v\in V$, $t\in [0,T]$. For the sake of clarity of the 
subsequent analysis we stress the different definitions of $p_{0,{\rm 
min}}$ as well as $p_{0,{\rm max}}$ where $p_0\in L^1(\Omega\times V)$, 
cf. Lemma \ref{Lemma3.1}, and $\hat{p}_{t,{\rm max}}$ where $p\in (L^1 
(\Omega\times V))^{[0,T]}$, cf. above (\ref{5.9}). 
\medskip

Let $\mathcal{N}$ denote the set of all non-negative $p_0\in L^1(\Omega 
\times V)$ with $\|1/p_0\|_{L^\infty(\Omega\times V)}<\infty$, $\|p_0 
\|_{L^\infty(\Omega\times V)}<\infty$, and $\|p_0\|_{L^1(\Omega\times 
V)}=1$. For all $p_0\in \mathcal{N}$, we may and do assume (\ref{3.3}) 
with $p_{0,{\rm min}}>0$ and $p_{0,{\rm max}}<\infty$. 

Furthermore for $p_0\in \mathcal{N}$, let $\mathcal{M}\equiv \mathcal{M}(p_0)$ 
be the set of all $p\in (L^1(\Omega\times V))^{[0,T]}$ such that $\|p 
(\cdot,\cdot,t)\|_{L^1(\Omega\times V)}=1$ and $\frac12p_{0,{\rm min 
}}\le p(\cdot,\cdot,t)\le p_{0,{\rm max}}+\frac12p_{0,{\rm min}}$ a.e. 
on $\Omega\times V$, $t\in [0,T]$. We note that $\mathcal{M}$ may depend 
on the choice of $p_{0,{\rm min}}$ and $p_{0,{\rm max}}$. 

\begin{lemma}\label{Lemma5.1} 
Fix $T>0$ as well as $c\ge 1$, and let $b>0$ be the number defined in 
condition (v). \\
(a) Let $p_0\in \mathcal{N}$ and let $p\in (L^1(\Omega\times V))^{[0,T]}$ 
be non-negative with $\|p\|_{1,T}\le 1$. Furthermore, let $0\le\beta< 
p_{0,{\rm min}}$. If $0<\lambda\le (p_{0,{\rm min}}-\beta)/(T\|h_\gamma\| 
\|B\|\cdot\hat{p}_{T,{\rm max}}\, c^{\1}_{T,\rm max})$ then 
\begin{eqnarray*}
\|\Psi(p_0,p)(\cdot,\cdot,t)\|_{L^1(\Omega\times V)}=1 
\end{eqnarray*} 
and, for all $t\in [0,T]$, 
\begin{eqnarray*}
\beta\le\Psi(p_0,p)(\cdot,\cdot,t)\le p_{0,{\rm max}}+\lambda T\|h_\gamma\|b
\cdot\hat{p}_{T,{\rm max}}\, c^{\1}_{T,\rm max}\quad\mbox{\rm a.e. on }\Omega 
\times V. 
\end{eqnarray*} 
(b) Let $p_0\in \mathcal{N}$ and let $p\in (L^1(\Omega\times V))^{[0,T]}$ 
be non-negative with $\|p\|_{1,T}\le 1$. Furthermore, let $\beta=\frac12 
p_{0,{\rm min}}$ and $\ \hat{p}_{T,{\rm max}}\le p_{0,{\rm max}}+\frac12 
p_{0,{\rm min}}$. If $\ 0<\lambda\le p_{0,{\rm min}}/(2T\|h_\gamma\|b\cdot 
(p_{0,{\rm max}}+\frac12p_{0,{\rm min}})\, c^{\1}_{T,\rm max})$ then, for all 
$t\in [0,T]$, we have 
\begin{eqnarray*} 
\frac12p_{0,{\rm min}}\le\Psi\left(p_0,p\right)(\cdot,\cdot,t)\le p_{0, 
{\rm max}}+\frac12p_{0,{\rm min}}\quad\mbox{\rm a.e. on }\Omega\times V. 
\end{eqnarray*} 
(c) Let $p_0\in L^1(\Omega\times V)$ and $p_1,p_2\in (L^1(\Omega\times 
V))^{[0,T]}$ and $\delta\in (0,1)$. For $0<\lambda\le\delta/(2T\|h_\gamma 
\|B\|(\|p_1\|_{1,T}+\|p_2\|_{1,T}))$ we have 
\begin{eqnarray}\label{5.10}
\left\|\Psi\left(p_0,p_1\right)-\Psi\left(p_0,p_2\right)\right\|_{1,T} 
\le\delta\left\|p_1-p_2\right\|_{1,T}\, .
\end{eqnarray} 
(d) In particular, let $p_0\in \mathcal{N}$, $p_1,p_2\in \mathcal{M}$ and $\delta 
\in (0,1)$. For 
\begin{eqnarray*}
0<\lambda\le\frac{\delta}{2}\, p_{0,{\rm min}}/(T\|h_\gamma\|b\cdot\left( 
p_{0,{\rm max}}+\textstyle{\frac12}p_{0,{\rm min}}\right)\, c^{\1}_{T,\rm max}) 
\, , 
\end{eqnarray*} 
we have $\Psi\left(p_0,p_1\right),\Psi\left(p_0,p_2\right)\in \mathcal{M}$ 
and (\ref{5.10}). 
\end{lemma}
Proof. Recalling $\|p\|_{1,T}\le 1$, $p\ge 0$, and relation (\ref{2.8}) we find 
\begin{eqnarray*}
\Psi(p_0,p)(\cdot,\cdot,t)\le p_{0,{\rm max}}+\lambda\|h_\gamma\|b\cdot\int_0^t 
S(t-s)(\hat{p}_{T,{\rm max}}\cdot\1)\,  ds\, ,\quad t\in [0,T]. 
\end{eqnarray*} 
Parts (a) and (b) of the Lemma are now a consequence of (\ref{5.3}), 
(\ref{5.7}), and (\ref{5.9}). For part (b) we note that $2\|B\|\le b$, cf. 
condition (v). Parts (c) and (d) follow from (\ref{5.6}), and parts 
(a) and (b). For part (d) we note that $c^{\1}_{T,\rm max}\ge 1$. 
\qed 
\medskip 

Let $T>0$. Let us iteratively construct a solution $p\equiv p(p_0)$ to 
(\ref{2.4}) restricted to $(r,v,t)\in\Omega\times V\times [0,T]$ By 
(i)-(iii) we have then the boundary conditions (\ref{2.3}). Set 
\begin{eqnarray*}
p^{(0)}(\cdot,\cdot,t):=p_0\, ,\quad p^{(n)}(\cdot,\cdot,t):=\Psi\left(p_0, 
p^{(n-1)}\right)(\cdot,\cdot,t)\, ,\quad t\in [0,T],\ n\in \mathbb{N}. 
\end{eqnarray*} 
We note that with $d(q_1,q_2):=\left\|q_1-q_2\right\|_{1,T}$, $q_1,q_2\in 
\mathcal{M}$, the pair $(\mathcal{M},d)$ is a complete metric space. Furthermore, 
Lemma \ref{Lemma5.1} (a) and (b) imply that, for $\lambda$ as in Lemma 
\ref{Lemma5.1} (b), $q\in \mathcal{M}\equiv \mathcal{M}(p_0)$ yields $\Psi\left( 
p_0,q\right)\in\mathcal{M}$. An immediate consequence of Lemma \ref{Lemma5.1} 
(d) and the Banach fixed point theorem is now part (b) of the following 
proposition. The first part of (a) is a consequence of Lemma \ref{Lemma5.1} 
(c), (\ref{5.5}), and again the Banach fixed point theorem. The continuity 
statement in the second part of (a) follows from (\ref{5.6}). 
\begin{proposition}\label{Proposition5.2} 
Fix $T>0$ and let $b>0$ be the number defined in condition (v). \\ 
(a) Fix $0<\lambda\le 1/(16 T\|h_\gamma\|\|B\|)$. Let $p_0\in L^1(\Omega 
\times V)$ with $\|p_0\|_{L^1(\Omega\times V)}\le 3/2$. There is a unique 
element $p\equiv p(p_0)\in (L^1(\Omega\times V))^{[0,T]}$ such that  
\begin{eqnarray*}
p=\Psi\left(p_0,p\right)\, . 
\end{eqnarray*} 
The map $\{q_0\in L^1(\Omega\times V):\|q_0\|_{L^1(\Omega\times V)}\le 
3/2\}\ni p_0\mapsto p(p_0)\in (L^1(\Omega\times V))^{[0,T]}$ is continuous. 
 \\ 
(b) Assume $p_0\in \mathcal{N}$ and $0<\lambda<\frac{1}{2}p_{0,{\rm min}}/ 
(T\|h_\gamma\|b\cdot (p_{0,{\rm max}}+\frac12p_{0,{\rm min}})\, c^{\1}_{T, 
{\rm max}})$. There is a unique element $p\equiv p(p_0)\in \mathcal{M}$ such 
that  
\begin{eqnarray*}
p=\Psi\left(p_0,p\right)\, . 
\end{eqnarray*} 
In other words, for $p_0\in \mathcal{N}$ the equation (\ref{2.4}), restricted to 
$(r,v,t)\in\Omega\times V\times [0,T]$, has a unique solution $p\equiv p(p_0)$. 
This solution satisfies  
\begin{eqnarray*}
\frac12p_{0,{\rm min}}\le p(p_0)(\cdot,\cdot,t)\le p_{0,{\rm max}}+ 
\frac12p_{0,{\rm min}}\quad\mbox{\rm a.e. on }\Omega\times V 
\end{eqnarray*} 
and $\|p(p_0)(\cdot,\cdot,t)\|_{L^1(\Omega\times V)}=1$, $t\in [0,T]$. 
Furthermore, 
\begin{eqnarray*}
\lim_{n\to\infty}\left\|p^{(n)}-p(p_0)\right\|_{1,T}=0\, . 
\end{eqnarray*} 
\end{proposition}

\subsection{Global Solutions to the Boltzmann Type Equation}\label{sec:5:2}

Lemma \ref{Lemma3.3} and Theorem \ref{Theorem4.5} will now be applied 
to investigate the existence and uniqueness of solutions to the 
integrated (mild) version of the Boltzmann equation (\ref{2.4}). As 
above, let $(A,D(A))$ denote the infinitesimal operator of the strongly 
continuous semigroup $S(t)$, $t\ge 0$, in $L^1(\Omega\times V)$. Obviously, 
$(A,D(A))$ is also the infinitesimal operator of the strongly continuous 
group $S(t)$, $t\in\mathbb{R}$, whenever the hypotheses of Theorem 
\ref{Theorem4.5} are satisfied. 

The following proposition prepares global existence and uniqueness of 
equation (\ref{2.4}). The extension of this result to $t\ge 0$ and 
eventually to $t\in\mathbb{R}$ will be possible as soon as we have 
verified non-negativity of the solution to (\ref{2.4}) whenever $p_0\ge 
0$. The latter will be a consequence of Lemma \ref{Lemma5.5} below. 
\begin{proposition}\label{Proposition5.3}
(a) Let $p_0\in L^1(\Omega\times V)$. There exists $T_{max}\equiv T_{max} 
(p_0)\in (0,\infty]$ such that the following holds. The equation (\ref{2.4}) 
with $p(\cdot ,\cdot ,0)=p_0$ has a unique solution $p(\cdot ,\cdot ,t)\in 
L^1(\Omega\times V)$, $t\in [0,T_{max})$, which is continuous in $t\in [0, 
T_{max})$ with respect to the topology in $L^1(\Omega\times V)$. Moreover, 
if $T_{max}<\infty$ then $\lim_{t\uparrow T_{max}}\|p(\cdot ,\cdot ,t)\|_{ 
L^1(\Omega\times V)}=\infty$. \\ 
(b) Let $p_0\in D(A)$. The solution $p(\cdot ,\cdot ,t)$, $t\in [0,T_{max})$, 
to the equation (\ref{2.4}) given in part (a) is continuous on $[0,T_{max})$ 
and continuously differentiable on $(0,T_{max})$ in the topology of $L^1 
(\Omega\times V)$. It satisfies $p(\cdot ,\cdot, t)\in D(A)$ and 
\begin{eqnarray}\label{5.11}
\frac{d}{dt}\, p(\cdot ,\cdot ,t)=Ap(\cdot ,\cdot ,t)+\lambda\, Q(p,p)\, 
(\cdot ,\cdot ,t)\, ,\quad t\in [0,T_{max}).  
\end{eqnarray} 
Here, $d/dt$ is the derivative in $L^1(\Omega\times V)$. At $t=0$ it is the 
right derivative. 
 \\ 
(c) Let $p_0\in L^1(\Omega\times V)$. For the solution $p(\cdot ,\cdot ,t)$, 
$t\in [0,T_{max})$, to the equation (\ref{2.4}) given in part (a) it holds that 
\begin{eqnarray*} 
\int_\Omega\int_Vp_0(r,v)\, dv\, dr=\int_\Omega\int_V p(r,v,t)\, dv\, dr\, , 
\quad t\in [0,T_{max}).
\end{eqnarray*} 
\end{proposition}
Proof. {\it Step 1 } We prove part (a). We recall that the equation (\ref{2.4}) 
is the integrated (mild) version of (\ref{5.11}) with boundary conditions 
(\ref{2.3}) induced by the boundary conditions of $S(t)$, $t\ge 0$, see (iii). 
For $q_1,q_2\in L^1(\Omega\times V)$ we have 
\begin{eqnarray*}
&&\hspace{-.5cm}\left\|Q(q_1,q_1)-Q(q_2,q_2)\right\|_{L^1(\Omega\times V)}= 
\left\|Q(q_1+q_2,q_1-q_2)\right\|_{L^1(\Omega\times V)} \\ 
&&\hspace{.5cm}\le\|B\|\, \|h_\gamma\|\, \|q_1+q_2\|_{L^1(\Omega\times V)}\, 
\|q_1-q_2\|_{L^1(\Omega\times V)}\, , 
\end{eqnarray*}  
i.e. $L^1(\Omega\times V)\ni q\mapsto Q(q,q)\in L^1(\Omega\times V)$ is 
locally Lipschitz continuous with constant $2C\|B\|\, \|h_\gamma\|$ on 
$\{q\in L^1(\Omega\times V):\|q\|_{L^1(\Omega\times V)}\le C\}$ for any 
$C>0$. Let us keep in mind the strong continuity of $S(t)$, $t\ge 0$, in 
$L^1(\Omega\times V)$, cf. Lemma \ref{Lemma3.3}. The standard reference 
\cite{Pa83} Theorem 1.4 of Chapter 6, says now that there is a unique 
solution $p$ to (\ref{2.4}) with $p(\cdot ,\cdot ,0)=p_0$ on some time 
interval $t\in [0,T_{max})$ such that 
\begin{eqnarray*}
T_{max}<\infty\quad\mbox{\rm implies}\quad\lim_{t\uparrow T_{max}}\|p 
(\cdot ,\cdot ,t)\|_{L^1(\Omega\times V)}=\infty\, . 
\end{eqnarray*}  
Furthermore, the just quoted source says also that $p(\cdot ,\cdot ,t)$ 
is continuous on $t\in [0,T_{max})$ with respect to the topology in $L^1 
(\Omega\times V)$. 
\medskip 

\noindent
{\it Step 2 } We prove part (b). For $h\in L^1(\Omega\times V)$ we have  
\begin{eqnarray*}
&&\hspace{-.5cm}\frac{\left\|Q(q+h,q+h)-Q(q,q)-2Q(q,h)\right\|_{L^1(\Omega 
\times V)}}{\|h\|_{L^1(\Omega\times V)}} \\ 
&&\hspace{.5cm}=\left\|Q(h,h)\right\|_{L^1(\Omega\times V)}/\|h\|_{L^1( 
\Omega\times V)} 
\end{eqnarray*} 
which converges to zero as $\|h\|_{L^1(\Omega\times V)}\to 0$ by (\ref{5.2}).
Therefore the map $L^1(\Omega\times V)\ni q\to Q(q,q)\in L^1(\Omega\times V)$ 
is Fr\'echet differentiable. The Fr\'echet derivative at $q\in L^1(\Omega\times 
V)$ has the representation 
\begin{eqnarray*}
\nabla Q(q,q)(\cdot)=2 Q(\cdot,q)\, .  
\end{eqnarray*} 
Next we remind of the notation $B\equiv B(L^1(\Omega\times V),L^1(\Omega\times 
V))$ of the space of all bounded linear operators $L^1(\Omega\times V)\mapsto 
L^1(\Omega\times V)$ endowed with the operator norm. We observe that  
\begin{eqnarray*}
&&\hspace{-.5cm}\|\nabla Q(p,p)(\cdot)-\nabla Q(q,q)(\cdot)\|_B=2\|Q(\cdot,p)- 
Q(\cdot,q)\|_B \\ 
&&\hspace{.5cm}=\sup_{\|h\|_{L^1(\Omega\times V)}=1}\|Q(h,p-q)\|_{L^1(\Omega\times 
V)} 
\end{eqnarray*} 
converges to zero as $\|p-q\|_{L^1(\Omega\times V)}\to 0$ by (\ref{5.2}). Thus 
the Fr\'echet derivative $\nabla Q$ as a map $L^1(\Omega\times V)\mapsto B(L^1 
(\Omega\times V),L^1(\Omega\times V))$ is continuous on $L^1(\Omega\times V)$. 
Now \cite{Pa83} Theorem 1.5 of Chapter 6, says that if $p_0\in D(A)$ then the 
solution $p(\cdot ,\cdot, t)$ to (\ref{2.4}) satisfies $p(\cdot,\cdot, t)\in 
D(A)$ as well as (\ref{5.11}). 
\medskip

\noindent
{\it Step 3 } Part (c) is an immediate consequence of (\ref{5.7}). 
\qed 
\medskip 

We mention that, even if $p(\cdot ,\cdot ,0)=p_0$ is positively bounded from 
below, it is not yet clear whether or not $p(\cdot ,\cdot ,t)$ is a.e. 
non-negative for all $t\in [0,T_{\max})$. We will address this problem 
in Lemma \ref{Lemma5.5} below. Together with Proposition 
\ref{Proposition5.3} this will lead to global solutions on $t\in [0, 
\infty)$ of (\ref{2.4}). See Theorem \ref{Theorem5.6} and Corollary 
\ref{Corollary5.7} below. 
\medskip

Let $p_0$ be a probability density and let $0<T<T_{max}(p_0)$. Denote 
by $p$ the solution to (\ref{2.4}) on $\Omega\times V\times [0,T]$ with 
$p(\cdot ,\cdot ,0)=p_0$ given in Proposition \ref{Proposition5.3} (a). The 
continuity with respect to $L^1(\Omega\times V)$  of the map $[0,T]\ni 
t\mapsto p(\cdot ,\cdot ,t)$, stated in Proposition \ref{Proposition5.3} 
(a), implies that 
\begin{eqnarray}\label{5.12}
\|p\|_{1,T}\le C_T\quad\mbox{\rm for some }0<C_T<\infty\, . 
\end{eqnarray} 

Let us recall the notation of Section \ref{sec:2}. In the present 
subsection we shall use the decomposition $Q(p,p)=Q^+(p,p)-Q^-(p,p)$ of 
the collision operator specified by 
\begin{eqnarray*}
Q^+(p,p)(r,v,t)=\int_V\int_{S_+^{d-1}}B(v,v_1,e)p(r,v^\ast,t)p_\gamma(r, 
v_1^\ast,t)\chi_{\{(v^\ast,v_1^\ast)\in V\times V\}}\, de\, dv_1 
\end{eqnarray*} 
in the sense of an element in $L^1(\Omega\times V)$. The expression 
\begin{eqnarray*}
\hat{B}(v,v_1):=\int_{S_+^{d-1}}B(v,v_1,e)\cdot\chi_{\{(v^\ast,v_1^\ast)\in V 
\times V\}}(v,v_1,e)\, de 
\end{eqnarray*} 
is by condition (v) well-defined for $(v,v_1)\in V\times V$. Keeping in mind 
(vi) and the definition of $p_\gamma$ in Section \ref{sec:2}, the term  
\begin{eqnarray*}
\hat{B}_p(r,v,t):=\lambda\int_V\hat{B}(v,v_1)\, p_\gamma(r,v_1,t)\, dv_1  
\end{eqnarray*} 
is well-defined and bounded on $(r,v,t)\in\overline{\Omega}\times V\times 
[0,T]$. Here 
\begin{eqnarray}\label{5.13}
\left|\hat{B}_p(r,v,t)\right|\le\lambda\|h_\gamma\|\|B\|\|p(\cdot,\cdot,t) 
\|_{L^1(\Omega\times V)}\le C_T\lambda\|h_\gamma\|\|B\|\, , 
\end{eqnarray} 
cf. also (\ref{5.12}). In fact, we have $Q^-(p,p)(\cdot ,\cdot ,t)=p(\cdot , 
\cdot ,t)\hat{B}_p(\cdot ,\cdot ,t)\in L^1(\Omega\times V)$. Moreover, by 
(vi), the map $\overline{\Omega}\ni r\mapsto p_\gamma(r,\cdot ,t)$ is 
bounded and uniformly continuous with respect to the topology of $L^1(V)$ for 
any $t\in [0,T]$. Thus by (v), $\hat{B}_p(\cdot ,\cdot ,t)$ is bounded and 
continuous on $\overline{\Omega}\times V$ for any $t\in [0,T]$. 
\medskip

\noindent
\br{13} Let $t\in [0,T]$. According to (\ref{2.4}) and Proposition 
\ref{Proposition5.3} (a) $[0,t]\ni s\mapsto S(t-s)\, Q(p,p)\, (\cdot ,\cdot ,s)$ 
is Bochner integrable, i.e. 
\begin{eqnarray*}
\int_0^t\|S(t-s)\, Q(p,p)\, (\cdot ,\cdot ,s)\|_{L^1(\Omega\times V)}\, ds<\infty 
\, .  
\end{eqnarray*} 
By (\ref{5.12}) we have $\int_0^t\|p(\cdot ,\cdot ,s)\|_{L^1(\Omega\times V)}\, 
ds<\infty$. Together with (\ref{5.13}), $\|S(t-s)\|_B=1$ as the operator norm 
in $B(L^1(\Omega\times V),L^1(\Omega\times V))$, and $Q^-(p,p)(\cdot ,\cdot ,s) 
=p(\cdot,\cdot ,s)\hat{B}_p(\cdot ,\cdot ,s)$, $s\in [0,t]$, this leads 
to 
\begin{eqnarray*}
\int_0^t\|S(t-s)\, Q^-(p,p)\, (\cdot ,\cdot ,s)\|_{L^1(\Omega\times V)}\, ds< 
\infty 
\end{eqnarray*} 
as well as
\begin{eqnarray*}
\int_0^t\|S(t-s)\, Q^+(p,p)\, (\cdot ,\cdot ,s)\|_{L^1(\Omega\times V)}\, ds< 
\infty\, .  
\end{eqnarray*} 
In other words, $S(t-\cdot)\, Q^-(p,p)\in L^1([0,t];L^1(\Omega\times V))$ as 
well as $S(t-\cdot)\, Q^+(p,p)\in L^1([0,t];L^1(\Omega\times V))$. According 
to \cite{Ja97}, Appendix C, both integrals 
\begin{eqnarray*}
\int_0^tS(t-s)\, Q^-(p,p)\, (\cdot ,\cdot ,s)\, ds\quad\mbox{\rm as well as}
\quad\int_0^tS(t-s)\, Q^+(p,p)\, (\cdot ,\cdot ,s)\, ds
\end{eqnarray*} 
can be evaluated a.e. on $\Omega\times V$ as ordinary Lebesgue integrals. 
Replacing the bounds, these integrals can also be evaluated over arbitrary 
Borel subsets of $[0,t]$ instead of over the whole interval $[0,t]$. 
\er
 
\br{14} From Remark \ref{13} and (\ref{2.4}) we obtain for a.e. 
$(r,v)\in\Omega\times V$ and all $0\le\tau\le T_{\Omega}(r,v)\wedge t$ 
\begin{eqnarray}\label{5.14}
&&\hspace{-.5cm}p(r-\tau v,v,t-\tau)-p(r,v,t)\vphantom{\int}\nonumber \\ 
&&\hspace{.5cm}=-\int_0^\tau\lambda Q(p,p)(r-sv,v,t-s)\, ds\nonumber \\ 
&&\hspace{.5cm}=\int_0^\tau\lambda Q^-(p,p)(r-sv,v,t-s)\, ds-\int_0^\tau\lambda 
Q^+(p,p)(r-sv,v,t-s)\, ds  \qquad
\end{eqnarray} 
which includes $\int_0^\tau|\lambda Q^-(p,p)(r-sv,v,t-s)|\, ds<\infty$ as well 
as $\int_0^\tau|\lambda Q^+(p,p)(r-sv,v,t-s)|\, ds<\infty$. 
By approximation on $\tau\in\{ta:a\in [0,1]$ is rational$\}$ along rays 
$\tau\mapsto r-\tau v$ relation (\ref{5.14}) holds also for a.e. $(r,v) 
\in\partial\Omega\times V$ with $v\circ n(r)\ge 0$. 
\er

For the next lemma we stress that $\hat{B}_p$ is well-defined and bounded 
on $\overline{\Omega}\times V\times [0,T]$ for any  $0<T<T_{max}$, cf. 
(\ref{5.13}). We remind of the notation diam$(\Omega)=\sup\{|r_1-r_2|: 
r_1,r_2\in\Omega\}$ and denote $T_m:={\rm diam}(\Omega)/v_{min}$. 
\begin{lemma}\label{Lemma5.4} 
Let $p_0\in L^1(\Omega\times V)$ and $0<T<T_{max}$. Let $p(\cdot ,\cdot 
,s)$, $s\in [0,T]$, be the solution to (\ref{2.4}) with $p(\cdot ,\cdot 
,0)=p_0$ given in Proposition \ref{Proposition5.3} (a). Then for all 
$t\in [0,T]$, all $(r,v)\in\Omega\times V$ as well as all $(r,v)\in 
\partial^{(1)}\Omega\times V$ such that $v\circ n(r)\ge 0$, and all $0\le 
\tau\le T_{\Omega}(r,v)\wedge t$ we have 
\begin{eqnarray*}
&&\hspace{-.5cm}0<\exp\left\{-\lambda (T_m\wedge T)\|h_\gamma\|\|B\|\right\} \\ 
&&\hspace{.5cm}\le\psi(r,v,t;\tau):=\exp\left\{\int_0^\tau\hat{B}_p(r-sv,v,t-s) 
\, ds\right\} \\ 
&&\hspace{.5cm}\le\exp\left\{\lambda (T_m\wedge T)\|h_\gamma\|\|B\|\right\}< 
\infty\, .  
\end{eqnarray*} 
Suppose (\ref{5.14}) for some $r,v,t,\tau$ as above. Then  
\begin{eqnarray}\label{5.15}
&&\hspace{-.5cm}p(r-\tau v,v,t-\tau)\nonumber \\ 
&&\hspace{.5cm}=\psi(r,v,t;\tau)\left(-\int_0^\tau\frac{\lambda Q^+(p,p)(r-sv,v, 
t-s)}{\psi(r,v,t;s)}\, ds+p(r,v,t)\right)\, . 
\end{eqnarray} 
\end{lemma}
Proof. As already noted below (\ref{5.13}), $\hat{B}_p(\cdot ,\cdot ,t)$ is 
bounded and continuous on $\overline{\Omega}\times V$ for any $t\in [0,T]$. In 
particular, $[0,T_\Omega\wedge t]\ni\tau \mapsto\hat{B}_p(r-\tau v,v,t)$ is 
continuous for all $(r,v)\in\Omega\times V$ or $(r,v)\in\partial^{(1)}\Omega\times V$ 
with $v\circ n(r)\ge 0$. We note also that by (\ref{5.13}) we have 
\begin{eqnarray*} 
0<\exp\left\{-\lambda (T_m\wedge T)\|h_\gamma\|\|B\|\right\}\le\inf\psi  
\end{eqnarray*} 
as well as 
\begin{eqnarray*} 
\sup\psi\le\exp\left\{\lambda (T_m\wedge T)\|h_\gamma\|\|B\|\right\}<\infty 
\end{eqnarray*} 
where both, the infimum as well as the supremum, are taken over the set $\{(r,v,t, 
\tau):(r,v)\in\Omega\times V,\ t\in [0,T],\ \tau\in [0,T_{\Omega}(r,v)\wedge t]\}$. 

For the rest of the proof, let $(r,v)\in\Omega\times V$ or $(r,v)\in\partial^{(1)}\Omega 
\times V$ with $v\circ n(r)\ge 0$ such that we have (\ref{5.14}). By the just 
mentioned boundedness of $\psi$, for 
\begin{eqnarray}\label{5.16}
&&\hspace{-.5cm}\varphi(\tau,t):=p(r,v,t+\tau)\, \psi(r,v,t+\tau;\tau)\nonumber \\ 
&&\hspace{.5cm}\equiv p(r,v,t+\tau)\exp\left\{\int_0^\tau\hat{B}_p(r-sv,v,t+\tau 
-s)\, ds\right\}\, ,\quad \tau\in [0,T_{\Omega}\wedge t],   
\end{eqnarray} 
the map $[0,T_{\Omega}\wedge t]\ni\tau\mapsto\varphi(\tau,t-\tau)$, $\tau\in [0,T_{ 
\Omega}\wedge t]$, is well-defined whenever (\ref{5.14}). From (\ref{5.16}) we 
get 
\begin{eqnarray}\label{5.17}
&&\hspace{-.5cm}\varphi(\tau,t-\tau)-\varphi(0,t)=p(r,v,t)\left(\exp 
\left\{\int_0^\tau\hat{B}_p(r-sv,v,t-s)\, ds\right\}-1\right)\nonumber \\ 
&&\hspace{.5cm}=\int_0^\tau\hat{B}_p(r-sv,v,t-s)p(r,v,t)\exp\left\{\int_0^s 
\hat{B}_p(r-uv,v,t-u)\, du\right\}\, ds\nonumber \\ 
&&\hspace{.5cm}=\int_0^\tau\hat{B}_p(r-s v,v,t-s)\varphi(s,t-s)\, ds\, .
\end{eqnarray} 
Recalling Remarks \ref{13} and \ref{14} we note that 
\begin{eqnarray}\label{5.18}
f(\tau,t):=\psi(r,v,t+\tau;\tau)\left(-\int_0^\tau\frac{\lambda Q^+(p,p) 
(r-sv,v,t+\tau-s)}{\psi(r,v,t+\tau;s)}\, ds+p(r,v,t+\tau)\right),\quad  
\end{eqnarray} 
$\tau\in [0,T_{\Omega}\wedge t]$, is well-defined for $(r,v)$ as chosen for 
this proof. Moreover we observe that by (\ref{5.16}) 
\begin{eqnarray}\label{5.19}
f(\tau,t)=-\psi(r,v,t+\tau;\tau)\int_0^\tau\frac{\lambda Q^+(p,p)(r-sv,v,t+\tau-s)} 
{\psi(r,v,t+\tau;s)}\, ds+\varphi(\tau ,t)\, .
\end{eqnarray} 
Using (\ref{5.16})-(\ref{5.19}) the following is a straightforward calculation, 
\begin{eqnarray}\label{5.20}
&&\hspace{-.5cm}f(\tau,t-\tau)-f(0,t)=-\psi(r,v,t;\tau)\int_0^\tau\frac{\lambda 
Q^+(p,p)(r-sv,v,t-s)}{\psi(r,v,t;s)}\, ds\nonumber \\ 
&&\hspace{1.0cm}+\varphi(\tau,t-\tau)-\varphi(0,t)\vphantom{\int_0^\tau}\nonumber 
 \\ 
&&\hspace{.5cm}=-\exp\left\{\int_0^\tau\hat{B}_p(r-uv,v,t-u)\, du\right\} 
\int_0^\tau\frac{\lambda Q^+(p,p)(r-sv,v,t-s)}{\psi(r,v,t;s)}\, ds\nonumber \\ 
&&\hspace{1.0cm}+\int_0^\tau\hat{B}_p(r-s v,v,t-s)\varphi(s,t-s)\, ds\nonumber \\
&&\hspace{.5cm}=-\int_0^\tau\hat{B}_p(r-sv,v,t-s)\exp\left\{\int_0^s\hat{B}_p 
(r-uv, v,t-u)\, du\right\}\times\nonumber \\ 
&&\hspace{1.0cm}\times\int_0^s\frac{\lambda Q^+(p,p)(r-uv,v,t-u)}{\psi(r,v,t;u)} 
\, du\, ds-\int_0^\tau\lambda Q^+(p,p)(r-uv,v,t-u)\, du\nonumber \\ 
&&\hspace{1.0cm}+\int_0^\tau\hat{B}_p(r-s v,v,t-s)\varphi(s,t-s)\, ds\nonumber \\
&&\hspace{.5cm}=\int_0^\tau\hat{B}_p(r-sv,v,t-s)\left(-\psi(r,v,t;s)\int_0^s\frac{ 
\lambda Q^+(p,p)(r-uv,v,t-u)}{\psi(r,v,t;u)}\, du\right.\nonumber \\ 
&&\hspace{1.0cm}\left.\vphantom{\int_0^s}+\varphi(s,t-s)\right)\, ds-\int_0^\tau 
\lambda Q^+(p,p)(r-uv,v,t-u)\, du\nonumber \\
&&\hspace{.5cm}=\int_0^\tau\left(\hat{B}_p(r-sv,v,t-s)f(s,t-s)-\lambda Q^+(p,p) 
(r-sv,v,t-s)\vphantom{\dot{f}}\right)\, ds\, . 
\end{eqnarray} 
Let us look at (\ref{5.20}) as an equation for $[0,T_{\Omega}\wedge t]\ni\tau\mapsto 
f(\tau ,t-\tau)$. Keeping in mind uniqueness of the related homogeneous equation 
(\ref{5.17}) we establish uniqueness of equation (\ref{5.20}) under the initial 
condition $f(0,t)=p(r,v,t)$. 

One representation of the solution to the equation (\ref{5.20}) under $f(0,t)= 
p(r,v,t)$ can be obtained from (\ref{5.18}) replacing there $t$ with $t-\tau$. 
In addition, it follows from (\ref{5.14}) that $f(\tau ,t-\tau):=p(r-\tau v,v,t 
-\tau)$, $0\le\tau\le T_{\Omega}\wedge t$, is a second representation. We have 
verified (\ref{5.15}). 
\qed
\begin{lemma}\label{Lemma5.5} 
Let $p_0\in L^1(\Omega\times V)$ and $p$ with $p(\cdot ,\cdot ,0)=p_0$ be 
the solution to (\ref{2.4}) given in Proposition \ref{Proposition5.3} (a). 
If there exist constants $0<c\le C<\infty$ with $c\le p_0\le C$ a.e. on 
$\Omega\times V$ then there exists a strictly decreasing positive function 
$[0,T_{max})\ni t\mapsto c_t$ such that $p(\cdot ,\cdot ,t)\ge c_t$ a.e. on 
$\Omega\times V$. 
\end{lemma}
Proof. Without loss of generality suppose $\|p_0\|_{L^1(\Omega\times V)} 
=1$. According to Proposition \ref{Proposition5.2} (b), for given $\lambda 
>0$ there exists $0<T<T_{max}$ and a unique solution $p$ on $\Omega\times 
V\times [0,T]$ to (\ref{2.4}) with $p(\cdot ,\cdot ,0)=p_0$ such that 
$0<c^\ast\le p(\cdot ,\cdot ,s)$ a.e. on $\Omega\times V$ for all $s\in [0,T]$ 
and some $c^\ast>0$. 

Let $\l$ denote the Lebesgue measure on $(\Omega\times V,\mathcal{B}(\Omega 
\times V))$ and assume now that there is $t^\ast\in [T,T_{max})$ such that 
\begin{eqnarray}\label{5.21}
\l\left(\{(r,v)\in\Omega\times V:p(r,v,s)\le 0\}\vphantom{l^1}\right)=0 
\quad\mbox{\rm for}\quad 0\le s<t^\ast 
\end{eqnarray}
and 
\begin{eqnarray}\label{5.22}
\l\left(\{(r,v)\in\Omega\times V:p(r,v,s)\le 0\}\vphantom{l^1}\right)>0 
\quad\mbox{\rm for}\quad s=t^\ast\quad\mbox{\rm or}\quad t^\ast<s<s_1   
\end{eqnarray}
for some $t^\ast<s_1<T_{max}$. Note that the opposite to this assumption is 
$p(\cdot ,\cdot ,s)>0$ a.e. on $\Omega\times V$ for every $s\in [0,T_{max})$. 
Therefore it is our objective to show that the assumption (\ref{5.21}), 
(\ref{5.22}) is false. 

Let $t\in [0,T_{max})$ and recall the times  
\begin{eqnarray*} 
t_k\equiv t_{\Omega,k}(r_0,v_0,\ldots ,v_{k-1}):=\sum_{l=1}^kT_{\Omega} 
(r_{l-1},v_{l-1})\, ,\quad 1\le k\le m\quad\mbox{\rm for }m\in\mathbb{N}, 
\end{eqnarray*} 
from Subsection \ref{sec:3:2} and set $\t t_k:=t_k\wedge t$. For a.e. 
$(r,v)\equiv (r_0,v_0)\in\Omega\times V$ and a.e. $(r,v)\equiv (r_0,v_0) 
\in\partial^{(1)}\Omega\times V$ with $v\circ n(r)\ge 0$ an iteration of 
(\ref{5.15}) leads to 
\begin{eqnarray}\label{5.23}
&&\hspace{-.5cm}p(r,v,t)=\left(\int_0^{{\t t}_1}\frac{\lambda Q^+(p,p) 
(r_0-s v_0,v_0,t-s)}{\psi(r_0,v_0,t;s)}\, ds\right.\nonumber \\ 
&&\hspace{2cm}+\left.\frac{\omega\chi_{[1,\infty)}(m)p\left(r_1,R_{r_1} 
(v_0),t-{t}_1\vphantom{l^1}\right)+\chi_{\{0\}}(m)p_0(r_e,v_0)}{\psi 
\left(r_0,v_0,t;{\t t}_1\right)}\vphantom{\int_0^{{\t t}_1}}\right) 
\nonumber \\
&&\hspace{.5cm}+\sum_{k=1}^\infty\, (1-\omega)M(r_1,v_0)\int_{v_1\circ 
n(r_1)\ge 0}(v_1\circ n(r_1))\ldots\times\nonumber \\ 
&&\hspace{1.0cm}\times (1-\omega)M(r_k,v_{k-1})\int_{v_k\circ n(r_k)\ge 
0}(v_k\circ n(r_k))\times\nonumber \\ 
&&\hspace{1.0cm}\times\left(\chi_{[k,\infty)}(m)\int_{t_k}^{{\t t}_{k+1}} 
\frac{\lambda Q^+(p,p)(r_k-(s-t_k)v_k,v_k,t-s)}{\psi\left(r_k,v_k,t-t_k; 
s-t_k\right)}\, ds\right.\nonumber \\ 
&&\hspace{1.5cm}\left.+\frac{\omega\, \chi_{[k+1,\infty)}(m)p\left(r_{k+1} 
,R_{r_{k+1}}(v_k),t-t_{k+1}\vphantom{l^1}\right)+\chi_{\{k\}}(m)p_0(r_e, 
v_k)}{\psi\left(r_k,v_k,t-t_k;{\t t}_{k+1}-t_k\right)}\vphantom{\int_{ 
{\t t}_{k-1}}^{{\t t}_k}}\right)\, dv_k\ldots\, dv_1\, \nonumber \\ 
\end{eqnarray} 
where we note that the structure of this sum is a slight reordering of 
the structure of the sums (\ref{3.20}) and (\ref{3.22}). By the 
assumption (\ref{5.21}), (\ref{5.22}) all summands are non-negative for 
$t\in [0,t^\ast)$. Thus, the infinite sum converges since the iteration 
of (\ref{5.15}) shows that any partial sum is bounded by the non-negative 
$p(r,v,t)$, $t\in [0,t^\ast)$. 

Let the function 
\begin{eqnarray*}
&&\hspace{-.5cm}g(r,v,t):=\chi_{\{0\}}(m)+\sum_{k=1}^\infty (1-\omega) 
M_{\rm min}\int_{v_1\circ n(r_1)\ge 0}(v_1\circ n(r_1))\ldots\times \\ 
&&\hspace{1.0cm}\times(1-\omega)M_{\rm min}\int_{v_k\circ n(r_k)\ge 0} 
(v_k\circ n(r_k))\chi_{\{k\}}(m)\, dv_k\ldots\, dv_1\, ,  
\end{eqnarray*} 
be defined for $t>0$ and all $(r,v)\in\left(\Omega\cup\partial^{(1)}\Omega 
\right)\times V$. It has the following properties.
\begin{itemize}
\item[(1)] We have $\lim_{t\to 0}g(r,v,t)=1$ for all $(r,v)\in\Omega\times 
V$ and all $(r,v)\in\partial^{(1)}\Omega\times V$ with $v\circ n(r)\ge 0$. 
\item[(2)] For fixed $t>0$ and $(r,v)\in\partial^{(1)}\Omega\times V$ with 
$v\circ n(r)\le 0$ the function $g$ satisfies the boundary conditions 
\begin{eqnarray*} 
g(r,v,t)=(1-\omega)M_{\rm min}\cdot J(r,t)(g)<J(r,t)(g)\cdot\left(\int_{w 
\circ n(r)\ge 0}w\circ n(r)\, dw\right)^{-1}\, .
\end{eqnarray*} 
\item[(3)] By (1) and (2) and the definition of $g$, the function $(0,\infty) 
\ni t\mapsto g(r,v,t)$ is non-increasing for all $(r,v)\in\left(\Omega\cup 
\partial^{(1)}\Omega\right)\times V$. 
\item[(4)] For $t>0$ and $(r,v)\in\partial^{(1)}\Omega\times V$ with $v\circ 
n(r)\ge 0$, the function $[0,T_\Omega(r,v)]\ni s\mapsto g(r-sv,v,t)$ is 
non-increasing. 
\end{itemize} 
As a consequence of properties (2)-(4), $\lim_{n\to 0}g(r_n,v_n,t)=0$ for 
fixed $t>0$ and some sequence $(r_n,v_n)\in\partial^{(1)}\Omega\times V$, 
$n\in\mathbb{N}$, implies $g(r,v,t)=0$ for a.e. $(r,v)\in\left(\Omega\cup 
\partial^{(1)}\Omega\right)\times V$. Thus 
\begin{eqnarray*}
\t c_t:=\mathop{\mathrm{ess~inf}}\limits_{(r,v)\in\Omega\times V} 
g(r,v,t)>0
\end{eqnarray*} 
and setting $\t c_0:=1$, the function $[0,\infty)\ni t\mapsto\t c_t$ is 
non-increasing. It follows from (\ref{5.23}) and Lemma \ref{Lemma5.4} 
that
\begin{eqnarray*}
&&\hspace{-.5cm}p (r,v,t)\ge\frac{\chi_{\{0\}}(m)p_0(r_e,v_0)}{\|\psi 
(\cdot,\cdot,t;t)\|_{L^\infty(\Omega\times V)}}+\sum_{k=1}^\infty(1- 
\omega)M(r_1,v_0)\int_{v_1\circ n(r_1)\ge 0}(v_1\circ n(r_1))\ldots 
\times \\ 
&&\hspace{1.0cm}\times(1-\omega)M(r_k,v_{k-1})\int_{v_k\circ n(r_k) 
\ge 0}(v_k\circ n(r_k))\frac{\chi_{\{k\}}(m)p_0(r_e,v_k)}{\|\psi( 
\cdot,\cdot,t;t)\|_{L^\infty(\Omega\times V)}}\, dv_k\ldots\, dv_1 \\ 
&&\hspace{.5cm}\ge\frac{c}{\exp\{\lambda t\|h_\gamma\|\|B\|\}}\, 
g(r,v,t)\, ,\quad (r,v,t)\in\Omega\times V\times (0,t^\ast), 
\end{eqnarray*} 
where we recall that for fixed $(r,v)$, the number $m$ is a function of 
$t$ and $v_1,v_2,\ldots ,v_e$. It follows that 
\begin{eqnarray}\label{5.24}
\mathop{\mathrm{ess~inf}}\limits_{(r,v)\in\Omega\times V}p(r,v,t)\ge 
\t c_t\cdot\frac{c}{\exp\{\lambda t\|h_\gamma\|\|B\|\}}=:c_t>0\, ,\quad
0\le t<t^\ast, 
\end{eqnarray} 
and that $[0,t^\ast)\ni t\mapsto c_t$ is strictly decreasing. Relations 
(\ref{5.24}) and (\ref{5.15}) imply that (\ref{5.22}) cannot hold. 

As already mentioned, this means $p(\cdot ,\cdot ,t)>0$ a.e. on $\Omega 
\times V$ for every $t\in [0,T_{max})$. Now we repeat the reasoning from  
(\ref{5.23}) on, for arbitrary $t\in [0,T_{max})$. We obtain (\ref{5.24}) 
for all $t\in [0,T_{max})$ including the definition of $c_t$ for all 
$t\in [0,T_{max})$ and that $[0,T_{max})\ni t\mapsto c_t$ is strictly 
decreasing.
\qed 
\medskip 

Let $b>0$ be the number defined in condition (v). Furthermore, let us 
remind of the notations $c^{\1}_{\infty,\rm max}=\sup_{t\ge 0}\|S(t)\, 
\1\|_{L^\infty(\Omega\times V)}$, which is finite by Lemma \ref{Lemma3.1}, 
and $p_{0,{\rm max}}$, cf. Lemma \ref{Lemma3.1}. We may now state the 
following result on the existence of global solutions on $t\in [0,\infty)$
to the equation (\ref{2.4}). 
\begin{theorem}\label{Theorem5.6} 
Let $p_0\in L^1(\Omega\times V)$ with $\|p_0\|_{L^1(\Omega\times V)}=1$ and 
suppose that there are constants $0<c\le C<\infty$ with $c\le p_0\le C$ a.e. 
on $\Omega\times V$. Then there is a unique map $[0,\infty)\ni t\mapsto p_t 
(p_0)\equiv p(\cdot ,\cdot ,t)\in L^1(\Omega\times V)$ such that (\ref{2.4})  
with $p(\cdot ,\cdot ,0)=p_0$. The solution $p\equiv p(p_0)$ to (\ref{2.4}) 
has the following properties. 
\begin{itemize} 
\item[(1)] The map $[0,\infty)\ni t\mapsto p(\cdot ,\cdot ,t)\in L^1(\Omega 
\times V)$ is continuous with respect to the topology in $L^1(\Omega\times 
V)$. 
\item[(2)] We have $\|p(\cdot ,\cdot ,t)\|_{L^1(\Omega\times V)}=1$, $t\ge 
0$. 
\item[(3)] There exists a strictly decreasing positive function $[0,\infty) 
\ni t\mapsto c_t$ such that $p(\cdot ,\cdot ,t)\ge c_t$ a.e. on $\Omega\times 
V$. 
\item[(4)] We have 
\begin{eqnarray*} 
\|p(\cdot ,\cdot ,t)\|_{L^\infty(\Omega\times V)}\le p_{0,{\rm max}}\cdot 
\exp\left\{\lambda\|h_\gamma\|\, b\, c^{\1}_{\infty,\rm max}\cdot t\vphantom 
{l^1}\right\}\, ,\quad t\ge 0.  
\end{eqnarray*}
\end{itemize} 
\end{theorem}
Proof. {\it Step 1 } The properties (1)-(3) are an immediate consequence of 
Proposition \ref{Proposition5.3}, parts (a) as well as (c), and Lemma 
\ref{Lemma5.5}. In particular we obtain $T_{max}=\infty$. 

Let us verify (4). It follows from Lemma \ref{Lemma3.1}, (\ref{2.4}), (\ref{2.8}), 
and property (3) of this theorem that 
\begin{eqnarray}\label{5.25}
&&\hspace{-.5cm}\|p(\cdot ,\cdot ,t)\|_{L^\infty(\Omega\times V)}\le p_{0, 
{\rm max}}+\lambda\|h_\gamma\|\, b\cdot\int_0^tS(t-s)(\|p(\cdot ,\cdot ,s) 
\|_{L^\infty(\Omega\times V)}\cdot\1)\, ds\nonumber \\ 
&&\hspace{.5cm}\le p_{0,{\rm max}}+\lambda\|h_\gamma\|\, b\, c^{\1}_{\infty, 
\rm max}\cdot\int_0^t\|p(\cdot ,\cdot ,s)\|_{L^\infty(\Omega\times V)}\, ds 
\, ,\quad t\ge 0,  
\end{eqnarray}
where we have not yet excluded that both sides are infinite. 

Recalling that $c^{\1}_{T,\rm max}\le c^{\1}_{\infty,\rm max}<\infty$ we 
may use Proposition \ref{Proposition5.2} (b) to claim that there is some 
$T\equiv T(\lambda,p_0)>0$ such that $\|p(\cdot ,\cdot ,t)\|_{L^\infty 
(\Omega\times V)}$ is bounded on $t\in [0,T]$. This allows to apply 
Gr\"onwall's inequality to (\ref{5.25}) in order to obtain (4) for $t\in 
[0,T]$. 
\medskip 

\noindent
{\it Step 2 } We keep in mind the particular form of (4) on every time 
interval $t\in [0,T_1]$ on which both sides of (\ref{5.25}) is finite. In 
order to apply Gr\"onwall's inequality on $t\ge 0$ it is sufficient to show 
that there is no $t_0>0$ such that 
\begin{eqnarray}\label{5.26}
\limsup_{t\uparrow t_0}\|p(\cdot ,\cdot ,t)\|_{L^\infty(\Omega\times V)}< 
\infty\quad\mbox{\rm and} \quad\|p(\cdot ,\cdot ,t_0)\|_{L^\infty (\Omega 
\times V)}=\infty. 
\end{eqnarray} 
This is what we are concerned with in the remainder of the proof. We therefore 
assume that there was a $t_0>0$ with (\ref{5.26}). It is our aim to lead this 
assumption to a contradiction. Let $n\in\mathbb{N}$ and 
\begin{eqnarray*}
B_n:=\left\{(r,v)\in\Omega\times V:\, p(r,v,t_0)>n\right\}\, . 
\end{eqnarray*} 
By the assumption we have $\l(B_n)>0$. Recalling Step 1 of the proof of Lemma 
\ref{Lemma3.3} let us choose $\ve_n>0$ such that $\l(B_n\cap(\Omega_{\ve_n} 
\times V))>0$ and $t_0\ge \ve_n/(2v_{max})$. For $s:=\ve_n/(2v_{max})$ we have  
\begin{eqnarray*}
r-\tau v\in\Omega\quad\mbox{\rm for all }\tau\in [0,s],\ r\in\Omega_{\ve_n},\ 
v\in V. 
\end{eqnarray*} 
This implies $\left\{r-\tau v:(r,v)\in B_n\cap(\Omega_{\ve_n}\times V) 
\right\}\subset\Omega$ and 
\begin{eqnarray}\label{5.27}
\l\left(\left\{(r-\tau v,v):(r,v)\in B_n\cap\left(\Omega_{\ve_n}\times V 
\right)\right\}\vphantom{l^1}\right)\ge c_n\, ,\quad\tau\in [0,s],  
\end{eqnarray} 
for some $c_n>0$. Choosing now $n>\limsup_{t\uparrow t_0}\|p(\cdot ,\cdot, 
t)\|_{L^\infty(\Omega\times V)}$ relations (\ref{5.26}) and (\ref{5.27}) 
contradict (\ref{5.14}).  

We have verified that there is no $t_0>0$ with (\ref{5.26}). Applying now 
Gr\"onwall's inequality to (\ref{5.25}) we obtain (4) for $t\ge 0$. 
\qed
\medskip

Recall that $(A,D(A))$ denotes the infinitesimal generator of the strongly 
continuous semigroup $S(t)$, $t\ge 0$, in $L^1(\Omega\times V)$. Together 
with Proposition \ref{Proposition5.3} (b) we obtain the following. 
\begin{corollary}\label{Corollary5.7}
Let $p_0\in D(A)$ with $\|p_0\|_{L^1(\Omega\times V)}=1$ and suppose that 
there are constants $0<c\le C<\infty$ with $c\le p_0\le C$ a.e. on $\Omega 
\times V$. Then there is a unique map $[0,\infty)\ni t\mapsto p_t(p_0) 
\equiv p(\cdot ,\cdot ,t)\in D(A)$, which is continuous on $[0,\infty)$ and 
continuously differentiable on $(0,\infty)$ with respect to the topology  
of $L^1(\Omega\times V)$, such that 
\begin{eqnarray*}
\frac{d}{dt}\, p(\cdot ,\cdot ,t)=Ap(\cdot ,\cdot ,t)+\lambda\, Q(p,p)\, 
(\cdot ,\cdot ,t)\, ,\quad t\ge 0, 
\end{eqnarray*} 
and $p(\cdot ,\cdot, 0)=p_0$. Here, $d/dt$ denotes the derivative in $L^1 
(\Omega\times V)$. At $t=0$ it is the right derivative. 

This map coincides with the solution to the equation (\ref{2.4}) for 
$p(\cdot ,\cdot, 0)=p_0\in D(A)$ of Theorem \ref{Theorem5.6} and has 
therefore the properties (2)-(4) of Theorem \ref{Theorem5.6}. 
\end{corollary} 

Introduce $a_1:=8\|h_\gamma\|\|B\|$, $-\tau_0\equiv -\tau_0(\lambda):= 
\max_{x\ge 1}\log(x)/(-m_1+a_1\lambda x)$ and note that $\tau_0(\lambda) 
\stack{\lambda\to 0}{\lra}$ $-\infty$. Let $K$ satisfy $-\tau_0=\log(K)/ 
(-m_1+a_1\lambda K)$. One main result of the paper is the following. 
\begin{theorem}\label{Theorem5.8} 
Suppose that the conditions of Theorem \ref{Theorem4.5} are satisfied. \\ 
(a) Let $p_0\in L^1(\Omega\times V)$ with $\|p_0\|_{L^1(\Omega\times V)}=1$ 
and suppose that there are constants $0<c\le C<\infty$ with $c\le p_0\le C$ 
a.e. on $\Omega\times V$. Then there is a unique map $[\tau_0,\infty)\ni 
t\mapsto p_t(p_0)\equiv p(\cdot,\cdot ,t)\in L^1(\Omega\times V)$ which is 
continuous with respect to the topology of $L^1(\Omega\times V)$ such that 
for every $\tau_0\le\tau\le 0$ 
\begin{eqnarray}\label{5.28} 
p(r,v,t)=S(t-\tau)\, p(r,v,\tau)+\lambda\int_\tau^t S(t-s)\, Q(p,p)\, (r,v,s) 
\, ds 
\end{eqnarray} 
a.e.~on $(r,v,t)\in\Omega\times V\times [\tau,\infty)$ and $p(\cdot ,\cdot , 
0)=p_0$. We have the following.
\begin{itemize} 
\item[(2')] We have $\|p(\cdot ,\cdot ,t)\|_{L^1(\Omega\times V)}\le K$, 
$t\ge\tau_0$ and $\|p(\cdot ,\cdot ,t)\|_{L^1(\Omega\times V)}=1$, $t\ge 0$. 
\item[(3')] If $p(\cdot,\cdot ,s)\ge c'$ a.e. on $\Omega\times V$ for some 
$c'>0$ and some $s\in [\tau_0,0]$ then there exists a strictly decreasing 
positive function $[s,\infty)\ni t\mapsto c_t\equiv c_t(c',s)$ such that 
$p(\cdot ,\cdot ,t)\ge c_t$ a.e. on $\Omega\times V$. 
\item[(4')] For $s$ introduced in (3') and $t\ge s$ we have 
\begin{eqnarray*}  
\|p(\cdot ,\cdot ,t)\|_{L^\infty(\Omega\times V)}\le p_{0,{\rm max}} 
\cdot\exp\left\{\lambda\|h_\gamma\|\, b\, c^{\1}_{\infty,\rm max}\cdot 
|t|\vphantom{l^1}\right\}\, .  
\end{eqnarray*}
\end{itemize} 
Conversely, for all $\tau_0<0$ there is a $\lambda>0$ such that 
(\ref{5.28}). If (\ref{5.28}) is modified by replacing $p$ with $p/(1 
\vee\|p/K\|_{L^1(\Omega\times V)})$, then the modified equation has a 
unique continuous solution $\mathbb{R}\ni t\mapsto \t p_t(p_0)\equiv 
\t p(\cdot,\cdot ,t)\in L^1(\Omega\times V)$ with $\t p_t(p_0)=p_t(p_0)$ 
for all $t\ge\tau_0$. The following holds.  
\begin{itemize} 
\item[(2'')] For all $t\in\mathbb{R}$ we have $\int_\Omega\int_V \t 
p(\cdot,\cdot ,t)\, dv\, dr =1$ and for all $t\le 0$ it holds that 
\begin{eqnarray*}
\left\|\t p_t(p_0)\right\|_{L^1(\Omega\times V)}\le\exp\left\{-t\cdot 
\left(-m_1+a_1\lambda K\vphantom{l^1}\right)\right\}\, . 
\end{eqnarray*}
\end{itemize} 
(b) Let $p_0\in D(A)$ with $\|p_0\|_{L^1(\Omega\times V)}=1$ and suppose 
that there are constants $0<c\le C<\infty$ with $c\le p_0\le C$ a.e. on 
$\Omega\times V$. Then there is a unique map $[\tau_0,\infty)\ni t\mapsto 
p_t(p_0)\equiv p(\cdot,\cdot ,t)\in D(A)$ which is continuous on $[\tau_0, 
\infty)$ and continuously differentiable on $(\tau_0,\infty)$ with respect 
to the topology of $L^1(\Omega\times V)$ such that 
\begin{eqnarray*}
\frac{d}{dt}\, p(\cdot ,\cdot ,t)=Ap(\cdot ,\cdot ,t)+\lambda\, Q(p,p)\, 
(\cdot ,\cdot ,t)\, ,\quad t\ge\tau_0, 
\end{eqnarray*} 
and $p(\cdot ,\cdot, 0)=p_0$. Here, $d/dt$ denotes the derivative in $L^1 
(\Omega\times V)$. At $t=\tau_0$ it is the right derivative. This map 
coincides with the solution to the equation (\ref{5.28}) if there $p(\cdot 
,\cdot, 0)=p_0\in D(A)$. We have properties (2')-(4') and (2'') of part (a).
\end{theorem} 
Proof. {\it Step 1 } In this step we establish a backward equation that 
in Step 2 will provide us with $p(\cdot ,\cdot,\tau)$ if $\tau_0\le\tau<0$. 
For a commented review of the calculus we use, see Remark \ref{15} below. 

Let $id$ denote the identity in $L^1(\Omega\times V)$. We recall that, 
according to Lemma \ref{Lemma3.3} and Corollary \ref{Corollary4.4}, 
the semigroup $e^{m_1t}\cdot S(-t)$, $t\ge 0$, is strongly continuous 
and contractive in $L^1(\Omega\times V)$. By Theorems 4.2 and 4.3 of 
Chapter 1 in \cite{Pa83}, it follows for its generator $(-A+m_1\, id, 
D(A))$ that $\left\|(\mu-m_1)f+Af\right\|_{L^1(\Omega\times V)}\ge\mu 
\left\|f\right\|_{L^1(\Omega\times V)}$, $f\in D(A)$, $\mu>0$, i.e. 
$(A-m_1\, id,D(A))$ is $m$-accretive. This implies in particular 
\begin{eqnarray*}
\|(f_1+\mu Af_1)-(f_2+\mu Af_2)\|_{L^1(\Omega\times V)}\ge(1+m_1\mu) 
\|f_1-f_2\|_{L^1(\Omega\times V)}\, ,
\end{eqnarray*}
for all $f_1,f_2\in D(A)$, $0<|m_1|\mu<1$, see also Remark \ref{15} 
below. Let $K\ge 1$ and introduce $\t f:=f/(1\vee\|f/K\|_{L^1(\Omega 
\times V)})$, $\t Q(f,f):=Q(\t f,\t f)$ and, accordingly, $\t Q^- 
(f,f):=Q^-(\t f,\t f)$ as well as $\t Q^+(f,f):=Q^+(\t f,\t f)$. 
For $0<|m_1|\mu<1$ it follows that
\begin{eqnarray*}
&&\hspace{-.5cm}\|f_1+\mu (Af_1+\lambda\t Q(f_1,f_1))-(f_2+\mu 
(Af_2+\lambda\t Q(f_2,f_2)))\|_{L^1(\Omega\times V)}\\ 
&&\hspace{.5cm}\ge\|(f_1+\mu Af_1)-(f_2+\mu Af_2)\|_{L^1(\Omega 
\times V)} \\ 
&&\hspace{1.0cm}-\mu\lambda\|\t Q^-(f_1,f_1)-\t Q^-(f_2,f_2)\|_{ 
L^1(\Omega\times V)}-\mu\lambda\|\t Q^+(f_1,f_1)-\t Q^+(f_2,f_2) 
\|_{L^1(\Omega\times V)} \\ 
&&\hspace{.5cm}= \|(f_1+\mu Af_1)-(f_2+\mu Af_2)\|_{L^1(\Omega 
\times V)}-2\mu\lambda\|\t Q^-(f_1,f_1)-\t Q^-(f_2,f_2)\|_{L^1 
(\Omega\times V)} \\ 
&&\hspace{.5cm}\ge (1+m_1\mu)\|f_1-f_2\|_{L^1(\Omega\times V)}-
2\mu\lambda\|Q^-(\t f_1-\t f_2,\t f_1+\t f_2)\|_{L^1(\Omega 
\times V)}\, ,
\end{eqnarray*}
where for the second last line we have applied a certain collision 
transformation, cf. \cite{RW05}, Lemma 1.10. Taking into consideration 
\begin{eqnarray*}
&&\hspace{-.5cm}\|\t f_1-\t f_2\|_{L^1(\Omega\times V)}=\left\|\frac 
{f_1}{1\vee\|f_1/K\|_{L^1}}-\frac{f_2}{1\vee\|f_2/K\|_{L^1}}\right\|_{ 
L^1(\Omega\times V)}  \\ 
&&\hspace{.5cm}\le\frac{\|f_1-f_2\|_{L^1}}{1\vee\|f_1/K\|_{L^1}}+\|f_2 
\|_{L^1}\left|\frac{1}{1\vee\|f_1/K\|_{L^1}}-\frac{1}{1\vee\|f_2/K\|_{ 
L^1}}\right|\le 2\|f_1-f_2\|_{L^1(\Omega\times V)}\, , 
\end{eqnarray*}
$\|\t f_1+\t f_2\|_{L^1(\Omega\times V)}\le\|\t f_1\|_{L^1(\Omega 
\times V)}+\|\t f_2\|_{L^1(\Omega\times V)}\le 2K$, and (iv) of Section 
\ref{sec:2} we find 
\begin{eqnarray}\label{5.29}
&&\hspace{-.5cm}\|f_1+\mu (Af_1+\lambda\t Q(f_1,f_1))-(f_2+\mu(Af_2+ 
\lambda\t Q(f_2,f_2)))\|_{L^1(\Omega\times V)}\nonumber \\ 
&&\hspace{.5cm}\ge(1+m_1\mu)\|f_1-f_2\|_{L^1(\Omega\times V)}-8\mu 
\lambda K\|h_\gamma\|\|B\|\|f_1-f_2\|_{L^1(\Omega\times V)}\nonumber \\ 
&&\hspace{.5cm}=(1-(-m_1+8\lambda K\|h_\gamma\|\|B\|)\mu)\|f_1-f_2\|_{ 
L^1(\Omega\times V)}\, . 
\end{eqnarray}
According to (3.5) in \cite{Ba10}, $A+\lambda\t Q(\cdot ,\cdot)$ is 
$(-m_1+8\lambda K\|h_\gamma\|\|B\|)$-accretive, cf. also Remark 
\ref{15} below.
\medskip

The existence of a unique solution to the equation 
\begin{eqnarray}\label{5.30}
&&\hspace{-.5cm}f+\mu(Af+\lambda\t Q(f,f))=g\, ,\quad g\in L^1(\Omega 
\times V), 
\end{eqnarray}
for $0<\mu<(-m_1+8\lambda K\|h_\gamma\|\|B\|)^{-1}$ is a consequence of 
the following two facts. Firstly, for fixed $\t g\in L^1(\Omega\times 
V)$ introduce $\mathcal{T}f:=\t g-(id+\mu A)^{-1}\mu\lambda\t Q(f,f)$ 
for any $f\in L^1(\Omega\times V)$, recall that $\|[id+\mu/(1+\mu\,  
m_1)(A-m_1 id)]^{-1}\|_{L^1(\Omega\times V)}\le 1$, and observe that as 
above 
\begin{eqnarray*}
&&\hspace{-.5cm}\|\mathcal{T}f_1-\mathcal{T}f_2\|_{L^1(\Omega\times V)} 
\le\mu(1+\mu\, m_1)^{-1}\times  \\ 
&&\hspace{1.0cm}\times\|[id+\mu/(1+\mu\, m_1) (A-m_1 id)]^{-1}\|_{L^1( 
\Omega\times V)}\cdot\lambda\|\t Q(f_1,f_1)-\t Q(f_2,f_2)\|_{L^1(\Omega 
\times V)} \\ 
&&\hspace{.5cm}\le 2\mu(1+\mu\, m_1)^{-1}\lambda\|\t Q^-(f_1,f_1)-\t 
Q^-(f_2,f_2)\|_{L^1(\Omega\times V)} \\ 
&&\hspace{.5cm}=2\mu(1+\mu\, m_1)^{-1}\lambda\|Q^-(\t f_1-\t f_2,\t f_1 
+\t f_2)\|_{L^1(\Omega\times V)} \\ 
&&\hspace{.5cm}\le 8\lambda K\mu(1+\mu\, m_1)^{-1}\|h_\gamma\|\|B\|\|f_1 
-f_2\|_{L^1(\Omega\times V)}\, . 
\end{eqnarray*}
Thus for $0<\mu<(-m_1+8\lambda K\|h_\gamma\|\|B\|)^{-1}$ the map $\mathcal 
{T}$ is a contraction. Secondly, 
\begin{eqnarray*}
&&\hspace{-.5cm}id+\mu(A+\lambda\t Q(\cdot ,\cdot))=(id+\mu A)\left(id+ 
(id+\mu A)^{-1}\mu\lambda\t Q(\cdot ,\cdot)\right)\, . 
\end{eqnarray*}
Since $(A-m_1\, id,D(A))$ is $m$-accretive, the range of $id+\mu A$ is 
$L^1(\Omega\times V)$ whenever $0<\mu<-1/m_1$. Banach's fixed point 
theorem implies now that the range of $id+\mu(A+\lambda\t Q(\cdot ,\cdot 
))$ is $L^1(\Omega\times V)$ for all $0<\mu<(-m_1+8\lambda K\|h_\gamma\| 
\|B\|)^{-1}$, i.e. we have (\ref{5.30}). Furthermore, choosing $f_2=0$ 
in (\ref{5.29}) we obtain 
\begin{eqnarray*}
\left\|\left(id+\mu(A+\lambda\t Q(\cdot ,\cdot))\right)^{-1}f\right 
\|_{L^1(\Omega\times V)}\le \left(1+\mu(m_1-8\lambda K\|h_\gamma\| 
\|B\|)\right)^{-1}\|f\|_{L^1(\Omega\times V)} 
\end{eqnarray*}
for $0<\mu<(-m_1+8\lambda K\|h_\gamma\|\|B\|)^{-1}$. 
\medskip

According to (\ref{5.29}), (\ref{5.30}), and the Crandall-Liggett 
theorem, as stated in Theorem 4.3 of \cite{Ba10}, for any $f\in L^1 
(\Omega\times V)$ there is a map $[0,\infty)\mapsto\t T(-t)f\in L^1 
(\Omega\times V)$ given by 
\begin{eqnarray}\label{5.31}
\t T(-t)f=\lim_{n\to\infty}\left(id +\frac{t}{n}(A+\lambda\t Q(\cdot, 
\cdot))\right)^{-n}f\, . 
\end{eqnarray}
for all $f\in L^1(\Omega\times V)$ and all $t\ge 0$. In addition, with 
$a_1=8\|h_\gamma\|\|B\|$ we get from the last estimate
\begin{eqnarray}\label{5.32}
\left\|\t T(-t)f\right\|_{L^1(\Omega\times V)}\le\exp\left\{t\cdot 
\left(-m_1+a_1\lambda K\vphantom{l^1}\right)\right\}\left\|f\right\|_{ 
L^1(\Omega\times V)}\, ,\quad t\ge 0. 
\end{eqnarray}

By \cite{We72}, Proposition 3.18, and Remark \ref{15} below, $\t T 
(-t)$, $t\ge 0$, is the unique non-linear strongly continuous semigroup 
on $L^1(\Omega \times V)$ satisfying  
\begin{eqnarray}\label{5.33}
&&\hspace{-.5cm}\t T(-t)f=e^{m_1t}S(-t)f-\int_0^t e^{m_1(t-u)}S(-t+u) 
\left(\lambda\t Q(\t T(-u)f,\t T(-u)f)+m_1 \t T(-u)f\right)\, du 
\nonumber \\ 
&&\hspace{.5cm}=S(-t)f-\int_0^t S(-t+u)\left(\lambda\t Q(\t T(-u)f,\t 
T(-u)f)\right)\, du\, . 
\end{eqnarray}
To verify the adjusted hypotheses for Proposition 3.18 of \cite{We72}, 
as they are formulated in Remark \ref{15} below, is no additional 
work. We remark that we already have mentioned that $(A-m_1\, id,D(A))$ 
is $m$-accretive in the beginning of the proof. Furthermore, we recall 
that in (\ref{5.29}) we have demonstrated that $A+\lambda\t Q(\cdot , 
\cdot )$ is $(-m_1+8\lambda K\|h_\gamma\|\|B\|)$-accretive. In addition, 
we have (\ref{5.30}). As part of (\ref{5.29}) we have already have 
shown that $\lambda\|\t Q(f_1,f_1)-\t Q(f_2,f_2)\|_{L^1(\Omega\times 
V)}\le 8\lambda K\|h_\gamma\|\|B\|\|f_1-f_2\|_{L^1(\Omega\times V)}$ 
for $f_1,f_2\in L^1(\Omega\times V)$. Thus, we have 
\begin{eqnarray*}
&&\hspace{-.5cm}\|(\lambda\t Q(f_1,f_1)+m_1f_1)-(\lambda\t Q(f_2,f_2) 
+m_1f_2)\|_{L^1(\Omega\times V)} \\ 
&&\hspace{.5cm}\le (-m_1+8\lambda K\|h_\gamma\|\|B\|)\|f_1-f_2\|_{L^1 
(\Omega\times V)}\, ,\quad f_1,f_2\in L^1(\Omega\times V). 
\end{eqnarray*}
In conclusion, the last paragraph is sufficient for Proposition 3.18 
in \cite{We72} under the adjusted hypotheses of Remark \ref{15} below, 
i.e. (\ref{5.31}) is the unique solution to the integral equation 
(\ref{5.33}). 
\medskip

For given $\lambda>0$, let $t_0=\max_{x\ge 1}\log(x)/(-m_1+a_1\lambda x)$ 
and $K\ge 1$ be the unique number satisfying $t_0=\log(K)/(-m_1+a_1\lambda 
K)$. It follows from (\ref{5.32}) that $\|\t T(-t)f\|_{L^1(\Omega\times V)} 
\le K\|f\|_{L^1(\Omega\times V)}$ for $-t\in [-t_0,0]$. Observing 
that for $F\in L^1(\Omega\times V)$ with $\|F\|_{L^1(\Omega\times V)}\le K$ 
we have $\t Q(F,F)=Q(F,F)$, for all $f$ with $\|f\|_{L^1(\Omega\times V)}\le 
1$ there is a unique solution $T(-t)f:=\t T(-t)f$, $-t\in [-t_0,0]$, to 
(\ref{5.33}) with $\t Q$ replaced by $Q$, that is 
\begin{eqnarray}\label{5.34}
T(-t)f=S(-t)f-\lambda\int_0^t S(-t+u)Q(T(-u)f,T(-u)f)\, du\, ,\quad 
-t\in [-t_0,0]. 
\end{eqnarray}
Moreover, for all $t_0>0$ there is a $\lambda>0$ such that (\ref{5.34}) 
since $\tau_0(\lambda)\stack{\lambda\to 0}{\lra}-\infty$. Recall from 
Corollary \ref{Corollary4.6} (b) that $S(-t)$, $t\ge 0$, is a semigroup 
for which $\int_\Omega\int_VS(-t)f\, dv\, dr=\int_\Omega\int_V f\, dv\, 
dr$, $t\ge 0$, and recall also the definition of $\t Q$. With this, 
representation (\ref{5.33}) implies $\int_\Omega\int_V\t T(-t)f\, dv\, dr= 
\int_\Omega\int_V f\, dv\, dr$, $t\ge 0$. 
\medskip

\nid
{\it Step 2 } Let $t\ge 0$. According to the result of Step 1, the equation 
\begin{eqnarray}\label{5.35}
p(r,v,t-u)=S(-u)p_0(r,v)-\lambda\int_0^u S(-u+s)Q(p,p)\, (r,v,t-s)\, ds 
\end{eqnarray} 
with $p_0=p(\cdot,\cdot,t)$ and $\|p_0\|_{L^1(\Omega\times V)}=1$ has a 
unique solution $p(\cdot,\cdot,t-u)$ on $u\in [0,-\tau_0]=[0,t_0]$ which 
is continuous in $L^1(\Omega\times V)$ 
by \cite{Pa83} Theorem 1.4 of Chapter 6. Applying $S(u)$ to both sides 
of (\ref{5.35}) and substituting $u$ with $u-w$ and $t$ with $t-w$, for 
any $w\in [0,u]$ we obtain 
\begin{eqnarray*}
\quad p(r,v,t-w)=S(u-w)p(r,v,t-u)+\lambda\int_0^{u-w}S(u-w-s)Q(p,p)\, 
(r,v,t-u+s)\, ds\, , 
\end{eqnarray*} 
i.e. equation (\ref{5.35}) constructs the solution to (\ref{2.4}) backward 
in time from $t$ to $t-u$ for any $u\in [0,-\tau_0]$. Furthermore, choosing 
$t=0$ in this equation and replacing after that $-w$ by $t$ and $u$ by $- 
\tau$ we obtain (\ref{5.28}). 
Now we can apply Gr\"onwall's inequality as in Step 1 of the proof of 
Theorem \ref{Theorem5.6}. Setting in (\ref{5.35}) $t=0$ and $u=-\tau$ 
we can conclude as in (\ref{5.25}) to obtain 
\begin{eqnarray*}
\|p(\cdot ,\cdot ,\tau)\|_{L^\infty(\Omega\times V)}\le p_{0,{\rm max}}+
\lambda\|h_\gamma\|\, b\, c^{\1}_{\infty,\rm max}\cdot\int_\tau^0\|p(\cdot , 
\cdot ,\tau-s)\|_{L^\infty(\Omega\times V)}\, ds\, ,\quad \tau\in [s,0],  
\end{eqnarray*}
which implies $\|p(\cdot ,\cdot ,\tau)\|_{L^\infty(\Omega\times V)}\le p_{0, 
{\rm max}}\cdot\exp\{\lambda\|h_\gamma\|\, b\, c^{\1}_{\infty,\rm max}\cdot 
|\tau|\}$. The latter and (4) of Theorem \ref{Theorem5.6} result in (4'). 

Considerations similar to Step 2 hold for $\t p$ which corresponds to $\t 
T(-\cdot)$. 
\medskip 

\noindent 
{\it Step 3 } For part (a) we refer to Theorem \ref{Theorem4.5}, Corollary 
\ref{Corollary4.6}, and Theorem \ref{Theorem5.6}, and in particular to 
Steps 1 and 2 of the present proof. Keeping Corollary \ref{Corollary5.7} 
in mind, for part (b) there is just to demonstrate that $p_0\in D(A)$ 
implies $p(\cdot ,\cdot ,\tau)\in D(A)$, $\tau\in [\tau_0,0]$. However, 
analyzing the equation 
\begin{eqnarray}\label{5.36}
\frac{d}{dt}\, q(\cdot ,\cdot ,t)=-Aq(\cdot ,\cdot ,t)-\lambda\, Q(q,q)\, 
(\cdot ,\cdot ,t)\, ,\quad t\in [0,-\tau_0], \quad q(\cdot ,\cdot ,0)=p_0 
\, ,
\end{eqnarray} 
as in the proof of Proposition \ref{Proposition5.3}, Step 2, and Corollary 
\ref{Corollary5.7}, $p(\cdot ,\cdot ,\tau)\in D(A)$ follows from \cite{Pa83}, 
Theorem 1.5 of Chapter 6. The latter states that if $p_0\in D(-A)=D(A)$ then 
the solution $p(\cdot ,\cdot, t)$ to (\ref{5.35}) satisfies $p(\cdot,\cdot, t)
\in D(A)$ as well as (\ref{5.36}). Here we note that by Theorem 
\ref{Theorem4.5}, $(-A,D(A))$ is the generator of the strongly continuous 
semigroup $S(-u)$, $u\ge 0$. Also here, similar considerations hold for $\t p$. 
\qed
\br{15} 
Let $(X,\|\cdot\|)$ be a real Banach space, $(X^\ast,\|\cdot\|_*)$ the 
dual space, and let $(\cdot ,\cdot)$ denote the pairing between $X$ and 
$X^\ast$. Furthermore, let $X\ni x\mapsto J(x)=\{x^\ast\in X^\ast: 
(x^\ast,x) = \|x\|^2=\|x^\ast\|_*^2\}$ denote the possibly multivalued 
duality map, for the existence we refer to \cite{Ba10}, Section 1.1. A 
possibly multivalued operator $A$ from $X$ to $X$ is called {\it 
accretive} if for every pair $[x_1,y_1],[x_2,y_2]\in A$, there is $z\in  
J(x_1-x_2)$ such that $(y_1-y_2,z)\ge 0$. Denoting by $id$ the identity 
in $X$, the operator $A$ is said to be {\it $m$-accretive} if the range 
of $id +A$ is $X$. The operator $A$ is said to be {\it $w$-accretive} 
for some $w\in\mathbb{R}$, if $w\, id+A$ is accretive, and $A$ is called 
{\it $w$-$m$-accretive} if $w\, id+A$ is $m$-accretive. 

Proposition 3.1 of \cite{Ba10} says now that $A$ is accretive if and only 
if the inequality $\|x_1-x_2\|\le\|x_1-x_2+\mu (y_1-y_2)\|$ holds for all 
$\mu>0$, or equivalently for some $\mu >0$, and all $[x_i,y_i]\in A$, $i= 
1,2$. In this case, the operator $A$ is single-valued. As a consequence, 
$A$ is $w$-accretive if and only if $(1-\mu w)\|x_1-x_2\|\le\|x_1-x_2+\mu 
(y_1-y_2)\|$ holds for all $0<\mu<1/w$, or equivalently for some $0<\mu 
<1/w$, and all $[x_i,y_i]\in A$, $i = 1,2$, see (3.5) in \cite{Ba10}. 
%\medskip
 
If a single valued operator $A$ in a Banach space $(X,\|\cdot\|)$ is accretive 
then $-A$ is said to be {\it dissipative}. The Lumer-Philips theorem, Theorem 
4.3 in Chapter 1 of \cite{Pa83}, says that if $-A$ is the infinitesimal 
generator of a strongly continuous contraction semigroup of linear operators 
in $X$ then, for all $\mu>0$, the range of $\mu\, id+A$ is $X$ and $-A$ is 
dissipative. Furthermore, Theorem 4.2 in Chapter 1 of \cite{Pa83} states that 
$-A$ is dissipative if and only if $\|(\mu\, id+A)f\|\ge\mu\|f\|$ for all 
$f\in X$ and all $\mu>0$. 
%\medskip

Let $f\in L^1([0,T];X)$, $f_1,\ldots ,f_N\in X$, $\varepsilon >0$, $0=t_0\le 
t_1\le\ldots \le t_N$ such that $t_i-t_{i-1}<\varepsilon$ for all $i\in\{1, 
\ldots,N\}$ as well as $T-\varepsilon <t_N\le T$, and $\sum_{i=1}^N\int_{ 
t_{i-1}}^{t_i}\|f(s)-f_i\|\, ds< \varepsilon$. Consider a  piecewise constant 
function $z:[0,t_N]\mapsto X$ whose values $z_i$ on $(t_{i-1},t_i]$ satisfy 
the finite difference equation $$\frac{z_i- z_{i-1}}{t_i-t_{i-1}}+Az_i\ni 
f_i\, ,\quad i\in\{1,\ldots,N\}.$$ Such a function $z = (z_i)^N_{i=1}$ is 
called an {\it $\varepsilon$-approximate} solution to the Cauchy problem 
\begin{eqnarray}\label{5.37}
y'(t)+Ay(t)\ni f(t)\, ,\quad t\in [0,T],\quad y(0) =y_0,
\end{eqnarray}
if it satisfies in addition $\|z(0)-y_0\|< \varepsilon$. A {\it mild 
solution in the sense of \cite{Ba10}} of the Cauchy problem (\ref{5.37}) 
is a function $y\in C([0,T];X)$ with the property that for each $\varepsilon 
>0$ there is an $\varepsilon$-approximate solution $z$ of $y'+Ay\ni f$ on 
$[0,T]$ such that $\|y(t)-z(t)\|\le\varepsilon$ for all $t\in [0,T]$ and 
$y(0) = x$.

Let $w\in\mathbb {R}$ and $A$ be a $w$-accretive operator for which  
$\overline{D(A)}$ is for small $\mu >0$ a subset of the range of $id+\mu A$.  
Let $y_0\in\overline{D(A)}$. The {\it Crandall-Liggett theorem} (Theorem 4.3 
in \cite{Ba10}) states now that the Cauchy problem (\ref{5.37}) has for 
$f=0$ a unique mild solution in the sense of \cite{Ba10} for which $$y(t)= 
\lim_{n\to\infty}\left(id +\frac{t}{n}A\right)^{-n}\, y_0$$ uniformly in 
$t$ on compact intervals.
%\medskip 

The proof of Proposition 3.18 of \cite{We72} is for our purposes the link 
between a {\it mild solution in the sense of \cite{Ba10}} and the usual 
understanding of a {\it mild solution}. This proof uses 
$\varepsilon$-approximate solutions in the sense of \cite{Ba10} in order 
to approximate a mild solution in the usual sense. Recall that the latter 
is a certain integrated version of a first order differential equation in 
Banach space as introduced e.g. in the standard reference \cite{Pa83}, 
Definition 1.1 of Chapter 6. 

Reviewing the proof of Proposition 3.18 of \cite{We72}, we observe that we 
can adjust the hypotheses of this proposition. Instead of requiring that 
$A$ is $m$-accretive and $B$ is continuous and accretive in the context of 
\cite{We72}, it is sufficient to require that $A$ is $w_1$-$m$-accretive 
for some $w_1\in\mathbb{R}$, $A+B$ is $w_2$-$m$-accretive for some $w_2\in 
\mathbb{R}$, and that $B$ is globally Lipschitz on $X$. Ignoring the first 
three lines of the proof and leaving the rest unchanged but keeping the 
Crandall-Liggett theorem in the form of Theorem 4.3 in \cite{Ba10} in mind, 
we obtain Proposition 3.18 of \cite{We72} under the adjusted hypotheses. 
Now, the non-linear strongly continuous semigroup (3.17) of \cite{We72} is 
a semigroup of contractions on $X$ if and only if $w_2=0$. In other words, 
the two notions of a {\it mild solution} coincide for our purposes under the 
here adjusted hypotheses of Proposition 3.18 of \cite{We72}. 
\er 

\begin{thebibliography}{1} 

\bibitem{AG11}{Aoki, K., Golse, F.:} On the speed of approach to equilibrium 
for a collisionless gas. Kinet. Relat. Models {\bf 4}(1), 87-107 (2011).

\bibitem{Ba10}{Barbu, V.:} Nonlinear differential equations of monotone 
types in Banach spaces. Springer Monographs in Mathematics. Springer, 
New York (2010)

\bibitem{Be20}{Bernou, A.:} A semigroup approach to the convergence rate of 
a collisionless gas. Kinet. Relat. Models {\bf 13}(6), 1071-1106 (2020) 

\bibitem{BC02}{Bobylev, A. V., Cercignani, C.:} Exact eternal solutions of 
the Boltzmann equation. J. Statist. Phys. {\bf 106}(5-6), 1019-1038 (2002)

\bibitem{BGSS18}{Bodineau, T., Gallagher, I., Saint-Raymond, L.,  
Simonella, S.:} One-sided convergence in the Boltzmann-Grad limit. {Ann. 
Fac. Sci. Toulouse Math.} {\bf 27}(5), 985-1022 (2018) 

\bibitem{Br15-1}{Briant, M.:} Instantaneous exponential lower bound for 
solutions to the Boltzmann equation with Maxwellian diffusion boundary 
conditions. Kinet. Relat. Models {\bf8}(2), 281-308 (2015)

\bibitem{Br15-2}{Briant, M.:} Instantaneous filling of the vacuum for the 
full Boltzmann equation in convex domains. Arch. Ration. Mech. Anal. 
{\bf 218}(2), 985-1041 (2015)

\bibitem{CPW98}{Caprino, S., Pulvirenti, M., Wagner, W.:} Stationary 
particle systems approximating stationary solutions to the Boltzmann 
equation. {SIAM J. Math. Analysis} {\bf 29}(4), 913-934 (1998) 

\bibitem{CIP94}{Cercignani, C., Illner, R., Pulvirenti, M.:} The 
mathematical theory of dilute gases. Applied Mathematical Sciences, 
{\bf 106}. Springer, New York, 1994. 

\bibitem{Ch20}{Chen, H.:} Cercignani-Lampis boundary in the Boltzmann 
theory. Kinet. Relat. Models {\bf 13}(3), 549-597 (2020)

\bibitem{CM03}{Chernov, N., Markarian, R.:} Introduction to the ergodic 
theory of chaotic billiards. Second edition. Facultad de Ingenier\'{i}a
Universidad de la Rep\'{u}blica - Uruguay (2003) \\
{\tt https://www.fing.edu.uy/imerl/grupos/ssd/publicaciones/pdfs/articulos/ 
2003/2003CheMarIntrod.pdf}

\bibitem{EN00}{Engel, K.-J., Nagel, R.:} One-parameter semigroups 
for linear evolution equations. With contributions by S. Brendle, 
M. Campiti, T. Hahn, G. Metafune, G. Nickel, D. Pallara, C. Perazzoli, 
A. Rhandi, S. Romanelli and R. Schnaubelt. Graduate Texts in Mathematics 
{\bf 194}. Springer, New York (2000)

\bibitem{EN06}{Engel, K.-J., Nagel, R.:} A short course on operator 
semigroups. Universitext. Springer, New York (2006) 

\bibitem{EM20}{Esposito, R., Marra, R.:} Stationary non equilibrium 
states in kinetic theory. J. Stat. Phys. {\bf 180}(1-6), 773-809 (2020) 

\bibitem{Fo01} {Fournier, N.:} Strict positivity of the solution to a 
2-dimensional spatially homogeneous Boltzmann equation without cutoff. 
Ann. Inst. H. Poincar\'{e} Probab. Statist. {\bf 37}(4), 481-502 (2001) 

\bibitem{GBV09}{Gamba, I. M., Panferov, V., Villani, C.:} Upper 
Maxwellian bounds for the spatially homogeneous Boltzmann equation. 
Arch. Ration. Mech. Anal. {\bf 194}(1), 253-282 (2009) 

\bibitem{GS69}{Gikhman, I. I., Skorokhod, A. V.:} {Introduction to the 
theory of random processes}. Inc. W. B. Saunders Co., Philadelphia 
London Toronto (1969) 

\bibitem{IMS20}{Imbert, C., Mouhot, C., Silvestre, L.:} Gaussian lower 
bounds for the Boltzmann equation without cutoff. SIAM J. Math. Anal. 
{\bf 52}(3), 2930-2944 (2020)

\bibitem{Ja97}{Janson, S.:} Gaussian Hilbert spaces. Cambridge Tracts 
in Mathematics {\bf 129}. Cambridge University Press, Cambridge (1997) 

\bibitem{KLT13}{Kuo, H.-W., Liu, T.-P., Tsai, L.-C.:} Free molecular flow 
with boundary effect. Comm. Math. Phys. {\bf 318}(2), 375-409 (2013) 

\bibitem{Lo18}{L\"obus, J.-U.:} Boundedness of the stationary solution 
to the Boltzmann equation with spatial smearing, diffusive boundary 
conditions, and Lions' collision kernel. SIAM J. Math. Anal. {\bf 50}(6), 
5761-5782 (2018) 

\bibitem{LMR20}{Lods, B., Mokhtar-Kharroubi, M., Rudnicki, R.:} Invariant 
density and time asymptotics for collisionless kinetic equations with 
partly diffuse boundary operators. Ann. Inst. H. Poincar\'{e} Anal. Non 
Lin\'{e}aire. {\bf 37}(4) 877-923 (2020). 

\bibitem{Ma1879}{Maxwell, C.:} Stresses in Rarefied Gases arising from 
Inequalities of Temperature, Phil. Trans. Roy. Soc. London {\bf 170}, 
231-256 (1879)

\bibitem{Mi10}{Mischler, S.:} Kinetic equations with Maxwell boundary 
conditions. Ann. Sci. \'{E}c. Norm. Sup\'{e}r. (4) {\bf 43}(5), 719-760 
(2010). 

\bibitem{Mo05}{Mouhot, C.:} Quantitative lower bounds for the full 
Boltzmann equation. I. Periodic boundary conditions. Comm. Partial 
Differential Equations {\bf 30}(4-6), 881-917 (2005).  

\bibitem{Pa83}{Pazy, A.:} Semigroups of linear operators and applications 
to partial differential equations. Springer, New York (1983) 

\bibitem{PW97}{Pulvirenti, A., Wennberg, B.:} A Maxwellian lower bound 
for solutions to the Boltzmann equation. Comm. Math. Phys. {\bf 183}(1), 
145-160 (1997)

\bibitem{RV08}{Rezakhanlou, F., Villani, C.:} Entropy methods for the 
Boltzmann equation. Lectures from a Special Semester on Hydrodynamic Limits 
held at the Universit\'e de Paris VI, Paris, 2001. Edited by F. Golse and 
S. Olla. Lecture Notes in Mathematics, {\bf 1916}. Springer, Berlin, 2008. 

\bibitem{RW05}{Rjasanow, S., Wagner, W.:} Stochastic numerics for the 
Boltzmann equation. Springer Series in Computational Mathematics, {\bf 37}. 
Springer, Berlin (2005)

\bibitem{S-R09} {Saint-Raymond, L.:} Hydrodynamic limits of the Boltzmann 
equation. Lecture Notes in Mathematics, {\bf 1971}. Springer, Berlin, 2009.

\bibitem{We72} {Webb, G. F.:} Continuous nonlinear perturbations of linear 
accretive operators in Banach spaces. {J. Functional Analysis} {\bf 10}, 
191-203 (1972) 

\bibitem{WZ06}{Wei, J., Zhang, X.:} Eternal solutions of the Boltzmann equation 
near travelling Maxwellians. J. Math. Anal. Appl. {\bf 314}(1), 219-232 (2006) 
\end{thebibliography}
\end{document}